\documentclass[12pt,nosumlimits,twoside,reqno]{amsart}
\usepackage{amssymb}
\usepackage{amsmath}
\usepackage[bookmarksnumbered,plainpages,hyperindex]{hyperref}
\usepackage{amsfonts}
\usepackage{color}
\usepackage{graphics}
\usepackage[all]{xy}
\usepackage{ifpdf}
\ifpdf  \fi
\newtheorem{theorem}{\sc Theorem}[section]

\newtheorem{lemma}[theorem]{\sc Lemma}
\newtheorem{corollary}[theorem]{\sc Corollary}
\theoremstyle{definition}
\newtheorem{definition}[theorem]{\sc Definition}
\newtheorem{example}[theorem]{\sc Example}
\theoremstyle{remark}

\newtheorem{claim}[theorem]{}
\numberwithin{equation}{section}
\DeclareSymbolFont{largesymbols}{OMX}{yhex}{m}{n}
\DeclareMathAccent{\widetriangle}{\mathord}{largesymbols}{"E6}
\def\M{\mathcal M}
\def\V{\vee}
\def\L{\wedge}
\def\ot{\otimes}
\def\Hom{\mathsf{Hom}}
\def\End{\mathsf{End}}
\def\1{_{\langle 1 \rangle}}
\def\2{_{\langle 2 \rangle}}
\def\stac#1{\raise-.2cm\hbox{$\,\stackrel{\displaystyle\otimes}{\scriptscriptstyle{#1}}\,$}}
\def\cten#1{\raise-.2cm\hbox{$\stackrel{\displaystyle\widehat{\otimes}}
{\scriptscriptstyle{#1}}$}}
\def\bigcheck{\widetriangle{\ \ }}
\def\bighat{\widehat{\ \ }}

\def\black{\color{black}}

\newcommand{\wt}[1]{\widetriangle{#1}}

\begin{document}

\title{A categorical approach to cyclic duality}
\author{Gabriella B\"ohm}
\author{Drago\c s \c Stefan}
\address{Research Institute for Particle and Nuclear Physics, Budapest,
H-1525 Budapest 114, P.O.B.49, Hungary.}
\email{G.Bohm@rmki.kfki.hu}
\address{University of Bucharest, Faculty of Mathematics and Informatics,
Bucharest, 14 Academiei Street, Ro-010014, Romania.}
\email{drgstf@gmail.com} 
\date{\today}

\begin{abstract}
The aim of this paper is to provide a unifying categorical framework for the
many examples of para-(co)cyclic modules arising from Hopf cyclic
theory. Functoriality of the coefficients is immediate in this approach.
A functor corresponding to Connes's cyclic duality is constructed.
Our methods allow, in particular, to extend Hopf cyclic theory to (Hopf)
bialgebroids.
\end{abstract}
\maketitle
\tableofcontents

\section*{Introduction}

Cyclic cohomology extends and unifies cohomology theories like de Rham
cohomology and Lie algebra cohomology of matrices. It has applications e.g. in
homological algebra, algebraic topology, Lie algebras, algebraic K-theory and
so non-commutative differential geometry.

The foot-stone in cyclic cohomology theory is a so called {\em cocyclic
  object}, i.e. a cosimplicial object equipped with an isomorphism at each
grade $n$, roughly implementing a cyclic permutation of the coface and
codegeneracy morphisms. In particular, the $n+1$st power of this cocyclic
morphism is required to be the identity.
The study of cocyclic objects, or a quest of their examples, can be
divided to two steps. First one can deal with a more general structure, called
a {\em para-cocyclic object}, obtained by relaxing the requirement about the
$n+1$st power of the cocyclic morphism at grade $n$ to be trivial. Truly
cocyclic subobjects or quotients of para-cocyclic objects are then studied as
a subsequent step.

There are many known examples of (para-)cocyclic objects, relevant for
various purposes. A large family of examples, occurring as symmetries in
non-commutative differential geometry, is associated to (co)module
(co)algebras of bialgebras. The first example of this kind appeared in
\cite{ConMos:HopfCyc} where it was used by Connes and Moscovici to give a
geometrical interpretation of the non-commutative Chern-character. Further
examples of para-cocyclic objects, associated to (co)module (co)algebras of
Hopf algebras were constructed by Hajac at al. in \cite{HaKhRaSo2}, where also
non-trivial coefficients provided by (co)modules of the Hopf algebra were
introduced. As a most important achievement, also criteria (on the
coefficients) for the existence of truly cocyclic subobjects and quotients
were found. These constructions were extended to bialgebras (extending Hopf
algebras) by Kaygun in \cite{Kay:BialgCycHom} and \cite{Kay:UniHCyc}. A new
type of coefficients, so-called contramodules, was proposed by Brzezi\'nski in
\cite{Brz:coef}. In this way, currently there are known eight families of
para-cocyclic objects associated to bialgebras: A cosimplicial object can be
constructed from a module algebra or a comodule algebra, or from a module
coalgebra or a comodule coalgebra $A$ (yielding four possibilities), using
either a functor of the form $A \otimes (-)$ or a functor of the form
$\Hom(A,-)$ (doubling the number of examples). In each case there turns out to
be an appropriate choice of the coefficients resulting in a para-cocyclic
structure.

Dually to (para-)cocyclic objects, one may consider (para-)cyclic objects,
i.e. (para-) cocyclic objects in the opposite category. Using bialgebras,
there can be constructed again eight families of examples.

As it was observed by Connes in \cite{Co:CycDual}, the category of cyclic
objects and the category of cocyclic objects in a given category are
isomorphic. This isomorphism, called {\em cyclic duality}, does not extend to
the categories of para-cyclic and para-cocyclic objects only to their
appropriate subcategories. These (full) subcategories have those objects whose
para-(co)cyclic morphisms are isomorphisms at each grade, cf. Khalkhali and
Rangipour's work \cite{KR}.

For para-(co)cyclic objects associated to (co)module (co)algebras of
bialgebras, the para-(co)cyclic morphisms are not isomorphisms in general. They
are isomorphisms, however, if the bialgebra in question is a Hopf algebra
with an invertible antipode. In this case the eight families of associated
para-cocyclic objects and the eight families of para-cyclic objects turn
out to be pairwise related by cyclic duality.
\smallskip

The aim of this paper is to provide a general construction of para-(co)cyclic
objects, including in particular existing constructions in terms of
bialgebras together with their generalizations to bialgebroids,
cf. \cite{ConMos:DiffCyc}. We do {\em not} investigate here, however, the
existence of truly (co)cyclic subobjects or quotients.

An important antecedent work of similar aims is Kaygun's paper
\cite{Kay:UniHCyc}, where a universal construction of para-(co)cyclic objects,
including examples from bialgebras, was presented. The construction in this
work is built on monoids and comonoids in symmetric monoidal
categories. Therefore, while it is perfectly suitable to describe (co)module
(co)algebras of bialgebras, it has to be generalized in order to be able to
cope with bialgebroids over non-commutative base algebras. Such a
generalization (under the names {\em admissible septuple} and its {\em
  transposition map}) was proposed in our previous work \cite{BS}. In Section
\ref{sec:cat_A} of the current paper we introduce a category in which
admissible septuples and their transposition maps are special objects. This
newly introduced category ${\mathcal A}$ comes equipped with a functor
${\mathcal Z}^*$ from ${\mathcal A}$ to the category of para-cocyclic objects
in the category of functors. As a consequence, any object in ${\mathcal A}$
induces a functor from a category (of coefficients) to the category of
para-(co)cyclic objects in another category (usually the category of modules
over a commutative ring). Thus the resulting construction of para-(co)cyclic
modules is functorial for the choice of coefficients.

Behind the construction of the above category ${\mathcal A}$ there are some
2-categorical considerations.
Consider an abstract 2-category ${\mathcal T}$ with a single
0-cell $o$, a 1-cell $t: o \to o$ and 2-cells $\eta:o \Rightarrow
t$ and $\mu: tt \Rightarrow t$, such that $\mu \circ \eta t = t =
\mu \circ t \eta$ and $\mu \circ \mu t = \mu \circ t \mu$ (i.e.
such that $(t,\mu,\eta)$ is a monad in ${\mathcal T}$). A monad
$(T,m,u)$ on a category $\M$ can be described then as a 2-functor
$F$ from ${\mathcal T}$ to an appropriate 2-category $\mathsf{Cat}$ of
(some) categories, functors and natural transformations, such that
$Fo=\M$, $Ft=T$, $F\mu=m$ and $F\eta=u$. A lax natural
transformation between 2-functors $F,F': {\mathcal T} \to
\mathsf{Cat}$ is precisely the same as a monad morphism
$(Ft,F\mu,F\eta) \to (F't,F'\mu,F'\eta)$ in the sense of
\cite{Street:monad} (for a review see Section \ref{sec:prelims} below).
Extending this picture, we may consider an abstract 2-category ${\mathcal S}$
of three 0-cells $o$, $d$ and $c$, with two monads on $o$ related by a
distributive law $\phi$, together with a so-called ${\mathcal
S}(\phi,c)$-algebra in ${\mathcal S}(o,c)$ and a ${\mathcal
S}(d,\phi)$-algebra in ${\mathcal S}(d,o)$. (For a review of algebras over
distributive laws in \cite{Bur:D-alg} see Section \ref{sec:prelims} below.)
The objects in our category ${\mathcal A}$ are 2-functors ${\mathcal
S}\to \mathsf{Cat}$ and the morphisms are lax natural transformations between
them.
In Section \ref{sec:ex.A} we collect eight families of examples of objects in
${\mathcal A}$ associated to bialgebroids.

A dual construction in Section \ref{sec:cat_B} provides us with another
category ${\mathcal B}$, admitting a functor ${\mathcal Z}_*$ from ${\mathcal
B}$ to the category of para-cyclic objects in the category of functors. In
Section \ref{sec:ex.B} we construct eight families of examples of objects in
${\mathcal B}$ in terms of bialgebroids.

In Section \ref{sec:duality} we investigate the natural question how Connes'
cyclic duality functor (more precisely its extension $\bighat$ to the
categories of para-(co)cyclic objects with invertible para-(co)cyclic
morphisms at each grade) lifts to a functor $\bigcheck$ between appropriate
subcategories of ${\mathcal A}$ and ${\mathcal B}$, such that ${\mathcal Z}_*
\circ {\bigcheck {\ }}\cong {\bighat {\ }} \circ {\mathcal Z}^*$. Behind this
lifting there are liftings of (co)monads.

In Section \ref{sec:ex_duals} we show that the examples associated to
(co)module (co)algebras of a Hopf algebroid with a bijective antipode, belong
to the subcategories of ${\mathcal A}$ and ${\mathcal B}$ on which the lifted
cyclic duality functor  $\bigcheck$ is defined. We also check that these
examples are pairwise related by the functor $\bigcheck$.

The paper contains an Appendix, summarizing some facts about modules, comodules
and contramodules of bialgebroids and Hopf algebroids, used to construct the
examples in the paper.
\bigskip

{\bf Notation.}
Composition of functors and the corresponding composition of natural
transformations is denoted by juxtaposition. That is, for consecutive functors
$F:{\mathcal C}\to {\mathcal D}$ and $G:{\mathcal D}\to {\mathcal E}$, the
composite is denoted by $GF$, with object map $X \mapsto GFX$. For
endofunctors $T:{\mathcal C}\to {\mathcal C}$, we also write $TT=T^2$. For
natural transformations $ v :F \to F'$ and $ w :G \to
G'$, the value of $ v $ at an object $X$ is denoted by $ v
 X$. Moreover, $G v $ is a natural transformation $GF \to
GF'$ whose value at $X$ is given by $G( v  X)$ -- written simply as
$G v  X$. Similarly, $ w  F$ is a natural transformation
$GF \to G'F$, with value at $X$ given by $ w (FX)= w  FX$.
The composition of consecutive natural transformations is denoted by $\circ$.
The identity morphism of any object $X$ is denoted by the same symbol $X$.
 In order to simplify some computations, we shall use the string
representation of functors and their natural transformations. Throughout this
paper, the composition of functors is represented by horizontal juxtaposition
of strings, the functor acting first being represented by the rightmost
string. A natural transformation $F _{1} \cdots F _{n}\longrightarrow G
_{1}\cdots G _{m}$ will be represented as a `microchip' with $n$ inputs
$F_1,\dots,F_n$ at its top and $m$ outputs $G_1,\dots,G_m$ at its bottom. The
identity natural transformation of a functor $F$ will be represented just as a
string, without any box, as in the first diagram. For a natural transformation
$v:F\rightarrow F'$, we draw $Gv:GF\rightarrow GF'$ and $vH:FH\rightarrow F'H$
as in the second and the third pictures.
\[
\includegraphics{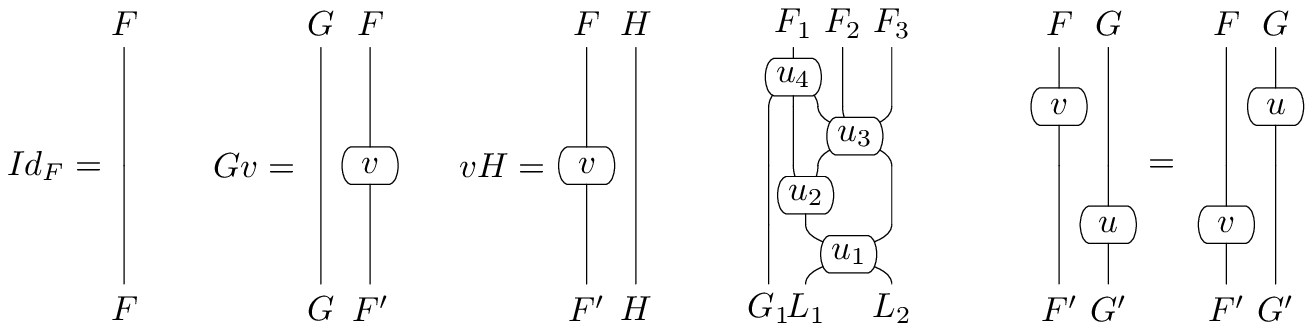}
\]
The composition of natural transformations is represented by vertical
juxtaposition of the corresponding `layers'. The natural transformation acting
first is represented by the top layer. Hence if translating diagrams to
formulae, we have to read our diagrams from the bottom to the top and from
left to right: For example, let us consider natural transformations
\[
u_1:K_1H_2\rightarrow L_1L_2, \quad
u_2:G_2H_1\rightarrow K_1,\quad
u_3:G_3F_2F_3\rightarrow H_1H_2,\quad
u_4:F_1\rightarrow G_1G_2G_3.
\]
Then the fourth diagram of the above picture is the string representation of
the natural transformation
\begin{equation*}
G_1u_1\circ G_1u_2H_2\circ G_1G_2u_3\circ u_4 F_2F_3.
\end{equation*}
In this notation, naturality of morphisms is visualized by their behaviour
as `pearls' on the strings. That is, those boxes which do not have common
ingoing or outgoing strings, can be freely moved above or below each other,
cf. the last equality in the above figure.

The above diagrammatic notation is used more generally in any
2-category: 1-cells are represented by vertical strings, their domains
corresponding to the surfaces on their right and codomains corresponding to
the surfaces on their left. 2-cells are represented by boxes, with domains
represented by `incoming legs' at their top and codomains represented by
`outgoing legs' at their bottom. Horizontal and vertical compositions in a
2-category are represented by horizontal and vertical juxtapositions of such
diagrams. Diagrams like the rightmost one above, come from the middle four
interchange law.

Note that in the literature the dual diagrammatic notation is used equally
frequently. In that case, our strings representing 1-cells are replaced by the
orthogonal lines -- hence surfaces on the sides of the original lines are
replaced by source and target points of the new orthogonal lines; and source
and target points of the original lines are replaced by surfaces on the
sides of the the new orthogonal lines. In this notation 2-cells are
represented by labels of the faces surrounded by their domain and codomain
1-cells. E.g. the 2-cell in the second figure above is represented as
$$
\xymatrix@=10pt{
\ar@{<-}[rr]^-G&&
\ar @{<-} @/^1pc/[rr]^-F \ar @{<-} @/_1pc/[rr]_-{F'}&
\Downarrow_v&
}
$$
Though both diagrammatic notations contain precisely the same
information, in this paper we prefer to work with string diagrams.
\black


\section{$\Phi$-module functors and their morphisms}\label{sec:prelims}

In this section we recall some notions from category theory, the
constructions of the later sections are built on.

The following notions are introduced in \cite{Street:monad}.

\begin{definition}\label{def:monad_morphism}
A {\em monad} on a category $\M$ is a triple $(T,m,u)$, where
$T:\M\to \M$ is a functor and $m:T^2\to T$ and $u:{\mathcal M} \to T$ are
natural transformations, called the {\em multiplication} and {\em unit},
respectively, which satisfy the last two ({\em associativity} and {\em
  unitality}) conditions in the following figure. (Our string representations
for the multiplication and the unit of a monad are introduced in the first two
equalities in the same figure.)
\[
\includegraphics{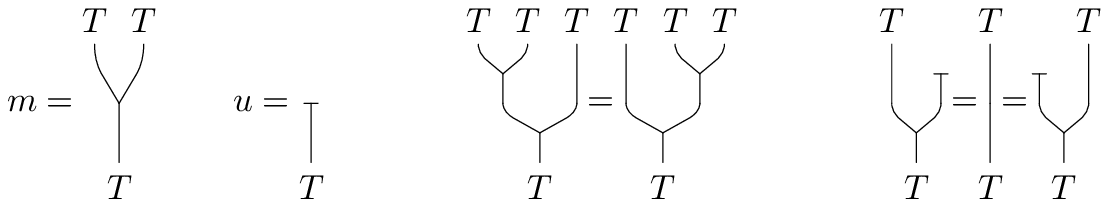}
\]
A {\em morphism} from a monad $(T',m',u')$ on $\mathcal M'$ to a monad
$(T,m,u)$ on $\mathcal M$ is a pair $(F,f)$, where  $F:\M'\to \M$ is a functor
and $f:TF \to FT'$ is a natural transformation which satisfies the following
two relations.
\[
\includegraphics{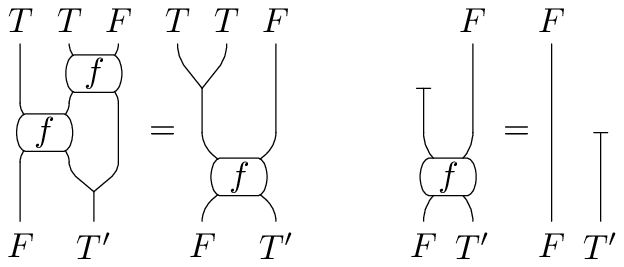}
\]
Monads and their morphisms constitute a category, with identity morphism $(\M,
T):(T,m,u)\to (T,m,u)$ and composition law $(G,g)\circ (F,f) =(GF,Gf\circ
gF)$, for monad morphisms $(F,f):(T'',m'',u'')\to (T',m',u')$ and
$(G,g):(T',m',u')\to (T,m, u)$.
\end{definition}

The category of monads and their morphisms is in fact isomorphic to a
full subcategory in the category of 2-functors and lax natural
transformations.
Let us denote in any 2-category the horizontal composition by juxtaposition
and the vertical composition by $\circ$.
Consider the 2-category $\mathsf{Cat}$ whose 0-cells are
some categories (whose monads we are aiming to describe), 1-cells are functors
and 2-cells are natural transformations.
(In order to avoid set theoretical problems arising from the paradox of
{\em ``the category of all categories''}, some restrictions on the 0-cells are
needed. Since our most important examples in Sections \ref{sec:ex.A} and
\ref{sec:ex.B} are (co)module categories, for our purposes allowing only for
small categories would be too restrictive. All of our examples are included,
for example, if 0-cells are all admissible or all locally presentable
categories, cf. \cite{porst}. Readers interested in other examples might
choose other classes (or in some cases even finite sets) of 0-cells to define
an appropriate (large) 2-category $\mathsf{Cat}$.) \black
On the other hand, consider the 2-category freely generated by a monad. That
is, the 2-category ${\mathcal T}$ with a single 0-cell $o$, a non-identity
1-cell $t: o \to o$ and its iterated horizontal composites, and 2-cells given
by composites of the non-identity 2-cells $\eta:o \Rightarrow t$ and $\mu: tt
\Rightarrow t$, modulo the relations $\mu \circ \eta t = t = \mu \circ t \eta$
and $\mu \circ \mu t = \mu \circ t \mu$ (for a diagrammatic representation see Definition
\ref{def:monad_morphism}); \black
that is, such that $(t,\mu,\eta)$ is a monad in ${\mathcal T}$.
A 2-functor $K:{\mathcal T}\to \mathsf{Cat}$ is precisely the same as a monad
$(Kt,K\mu,K\eta)$ on the category $Ko$.

A {\em lax natural transformation} $K' \to K$, for 2-functors $K,K':{\mathcal
  C}\to {\mathcal D}$ between any 2-categories ${\mathcal C}$ and ${\mathcal
  D}$, is given by 1-cells $F_C:K'C \to KC$ in ${\mathcal D}$, labelled by the
0-cells $C$ of ${\mathcal C}$, and 2-cells $f_\gamma:(K
\gamma)F_{C'} \to F_{C}(K'\gamma)$ in ${\mathcal D}$, labelled by the 1-cells
$\gamma:C' \to C$ in ${\mathcal C}$. These data obey the following conditions.
\begin{itemize}
\item[{(i)}] Naturality of $f$; that is, for any 2-cell $\Gamma:\gamma
  \to \delta$ in ${\mathcal C}$, the following diagram commutes.
$$
\xymatrix{
(K\gamma)F_{C'} \ar[r]^-{f_\gamma}\ar[d]_-{(K\Gamma)F_{C'}}&
F_{C}(K'\gamma) \ar[d]^-{F_{C}(K'\Gamma)}\\
(K\delta)F_{C'} \ar[r]^-{f_\delta}&
F_{C}(K'\delta)\ .
}
$$
\item[{(ii)}] Compatibility of $f$ with the horizontal composition; that is,
for any 1-cells \break
$\xymatrix@=1.8pc{C'' \ar[r]|{\,\gamma'\,}& C'\ar[r]|{\,\gamma\,}& C}$,
the identity
$f_{\gamma\gamma'}= f_{\gamma}(K'\gamma')\circ (K\gamma)f_{\gamma'}$ holds.
\item[{(iii)}] Compatibility of $f$ with the identity 1-cells; i.e., for any
  identity 1-cell $\xymatrix@=1.8pc{C\ar[r]|{\,C\,} & C}$, $f_C$ is equal to
  the identity 2-cell $\xymatrix@=2pc{F_C \ar[r]|(.47){\,F_C\,}& F_C}$.
\end{itemize}
A lax natural transformation between 2-functors $K,K': {\mathcal T}
\to \mathsf{Cat}$ is given then by a functor $F=F_o:K'o \to Ko$
and a natural transformation $f=f_t:(Kt)F \to F(K't)$. (The value of $f$ on
the other 1-cells $t^n$ for $n\neq 1$ is determined by the compatibility
conditions (ii) and (iii) with the horizontal composition and with the
identity 1-cell.) The naturality condition (i), applied to $\Gamma=\mu$ and
$\Gamma = \eta$, respectively, yields precisely the same conditions in
Definition \ref{def:monad_morphism} on a monad morphism
$(F,f):(K't,K'\mu,K'\eta) \to (Kt,K\mu,K\eta)$.

\begin{definition}
An {\em algebra} for a monad $(T,m,u)$ on a category $\M$, is a monad morphism
 from the identity monad on the terminal category (of a single object and
its identity morphism) to $(T,m,u)$.
That is, a pair $(M,\varrho)$, where $M$ is an object in ${\mathcal
  M}$ and $\varrho:TM\to M$ is a morphism in $\M$ such that
 the first two diagrams in
$$
\xymatrix{
T^2M \ar[r]^-{T\varrho}\ar[d]_-{mM}&
TM \ar[d]^-{\varrho}\\
TM\ar[r]^-{\varrho}&
M
}\qquad\qquad
\xymatrix{
M\ar[d]_-{uM}\ar@{=}[rd]&\\
TM \ar[r]^-{\varrho}&
M
}\qquad\qquad
\xymatrix{
TM' \ar[r]^-{T\varphi}\ar[d]_-{\varrho'}&
TM\ar[d]^-{\varrho}\\
M'\ar[r]^-{\varphi}&
M
}
$$
commute.
A {\em morphism} of $T$-algebras $(M',\varrho')\to (M,\varrho)$ is a morphism
$\varphi:M'\to M$ in $\M$ such that the third diagram above commutes. Algebras
of a monad $T$ and their morphisms constitute the so-called Eilenberg-Moore
category $\M^T$.
\end{definition}

Note that via composition on the right, a monad $T:\M \to {\mathcal
M}$ induces a monad $\mathsf{Cat}(T,-)$ on the category $\mathsf{Cat}(\M,-)$,
whose objects are functors of domain $\M$ and whose  morphisms are natural
transformations. Symmetrically, there is a monad $\mathsf{Cat}(-,T)$ acting by
composition on the left on the category $\mathsf{Cat}(-,\M)$,  whose objects
are functors of codomain $\M$ and whose morphisms are natural transformations.
In order to distinguish algebras of these induced monads from $T$-algebras, we
call an algebra of the monad $\mathsf{Cat}(T,-)$  a {\em right $T$-module
functor} and we term an algebra of the monad $\mathsf{Cat}(-,T)$  a {\em left
$T$-module functor}.

\begin{definition}\label{def:distr_law}
Consider two monads $(T_l,m_l,u_l)$ and $(T_r,m_r,u_r)$ on the same category
$\M$. A (monad) {\em distributive law} is a natural transformation
$\Phi:T_r T_l \to T_l T_r$ such that the following four relations hold
true.
\[
\includegraphics{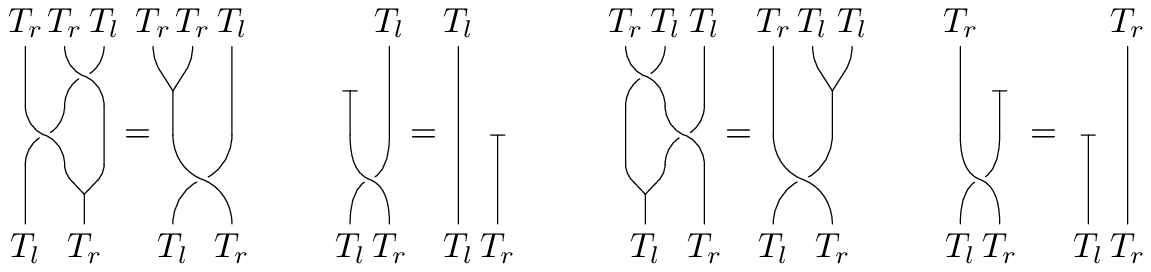}
\]
The last two defining relations are equivalent to the fact that
$(T_l,\Phi)$ is a monad endomorphism of $(T_r,m_r,u_r)$.

Note that in the representation of $\Phi$, the string corresponding to
$T_r$ crosses over the other one. If $\Phi$ is an isomorphism, then the string
representation of $\Phi^{-1}$ is obtained from that of $\Phi$ by an up-down
reflection.
\end{definition}

As it was proven in \cite{Be}, a distributive law $\Phi:T_r T_l \to T_l T_r$ as
in Definition \ref{def:distr_law} induces a monad structure on the composite
functor $T_lT_r$, with the multiplication $m$ and unit $u$, whose string
representations are
\begin{equation}\label{eq:composite_monad}
\includegraphics{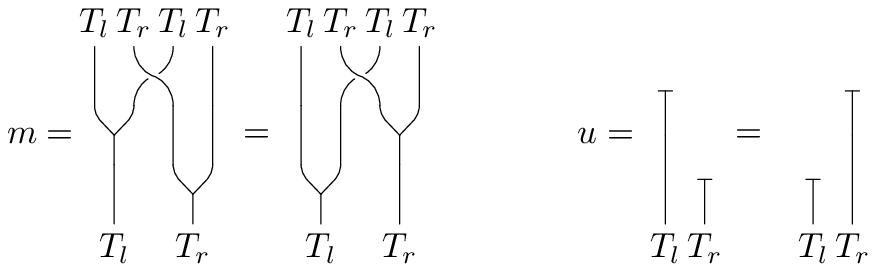}
\end{equation}
The following well-known lemma describes morphisms between such
composite monads.

\begin{lemma}\label{lem:comp_d_law}
Consider two monads $(T_l,m_l,u_l)$ and $(T_r,m_r,u_r)$ on the same category
$\M$ and two monads $(T'_l,m'_l,u'_l)$ and $(T'_r,m'_r,u'_r)$ on
$\M'$. Let $\Phi:T_r T_l \to T_l T_r$ and $\Phi':T'_r T'_l \to T'_l
T'_r$ be distributive laws. For monad morphisms
$(G,q_l): (T'_l,m'_l,u'_l)\to (T_l,m_l,u_l)$ and
$(G,q_r): (T'_r,m'_r,u'_r)\to (T_r,m_r,u_r)$, the
following assertions are equivalent.
\begin{itemize}
\item[{(i)}]
 $(G, q_l T'_r \circ T_l q_r)$
is a monad morphism
 $T'_l T'_r \to T_l T_r$ , cf. \eqref{eq:composite_monad};
\item[{(ii)}]  The following identity holds true.
\begin{equation}\label{eq:d-law_comp}
\includegraphics{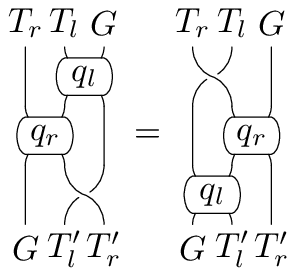}
\end{equation}
\end{itemize}
\end{lemma}

Distributive laws $(T_l,T_r,\Phi)$ as objects, with triples
$(G,q_l,q_r)$ satisfying the equivalent conditions in Lemma
\ref{lem:comp_d_law} as morphisms between them, constitute a category which
can be regarded again as a full subcategory in the category of
2-functors and lax natural transformations.
To this end, consider now a 2-category ${\mathcal L}$ of a single 0-cell $o$,
generated by two monads $(t_l,\mu_l,\eta_l)$ and $(t_r,\mu_r,\eta_r)$ on $o$
and a 2-cell $\phi:t_rt_l \to t_l t_r$ such that
\begin{eqnarray*}
&\phi \circ \mu_r t_l = t_l \mu_r \circ \phi t_r \circ t_r \phi
\qquad\qquad
&\phi \circ \eta_r t_l = t_l \eta_r \qquad\textrm{and}\\
&\phi \circ t_r\mu_l = \mu_l t_r \circ t_l \phi \circ\phi t_l
\qquad\qquad
&\phi \circ t_r \eta_l = \eta_l t_r
\end{eqnarray*}
(for a diagrammatic representation see Definition \ref{def:distr_law})
\black i.e. such that $\phi$ is a distributive law.
A 2-functor $K:{\mathcal L}\to \mathsf{Cat}$ is the same as a pair of monads
on the category $Ko$ related by a distributive law $K\phi$. A lax natural
transformation between such 2-functors is the same as a pair of monad
morphisms as in Lemma \ref{lem:comp_d_law}.

The following definition is quoted from \cite{Bur:D-alg}.

\begin{definition}\label{def:t-module}
Consider two monads $(T_l,m_l,u_l)$ and $(T_r,m_r,u_r)$ on the same category
$\M$ and a distributive law $\Phi:T_r T_l \to T_l T_r$.
A {\em $\Phi$-algebra} is a pair, consisting of an object $X$ of $\M$
and a morphism $\xi:T_r X \to T_l X$, such that
 the following diagrams commute.
$$
\xymatrix{
T_r^2X\ar[r]^-{T_r\xi}\ar[d]_-{m_rX}&
T_rT_lX\ar[r]^-{\Phi X}&
T_lT_rX\ar[r]^-{T_l\xi}&
T_l^2X\ar[d]^-{m_lX}\\
T_rX\ar[rrr]^-{\xi}&&&
T_lX
}\qquad \qquad
\xymatrix{
X\ar@{=}[r]\ar[d]_-{u_rX}&
X\ar[d]^-{u_lX}\\
T_rX\ar[r]^-{\xi}&
T_lX
}
$$
\end{definition}

If $T_l$ is equal to the identity functor $\M$ then there is a
trivial distributive law $\Phi=T_r$. In this case $\Phi$-algebras
are the same as  $T_r$-algebras.

A distributive law $\Phi:T_r T_l \to T_l T_r$, between monads on
$\M$, induces a distributive law
 $\mathsf{Cat}(\Phi,-): \mathsf{Cat}(T_l,-)\mathsf{Cat}(T_r,-)\to
\mathsf{Cat}(T_r,-)\mathsf{Cat}(T_l,-)$ between monads on $\mathsf{Cat}(\M,-)$
 and a distributive law
$\mathsf{Cat}(-,\Phi): \mathsf{Cat}(-,T_r)\mathsf{Cat}(-,T_l)\to
\mathsf{Cat}(-,T_l) \mathsf{Cat}(-,T_r)$ between monads on
$\mathsf{Cat}(-,\M)$ .
Algebras for these induced
distributive laws are called {\em right and left $\Phi$-module functors},
respectively.
Explicitly, a right $\Phi$-module functor is a pair, consisting of a functor
$\sqcap:\M \to {\mathcal C}$  (where ${\mathcal C}$ is any category)
and a natural transformation $i:\sqcap T_l \to \sqcap T_r$, such that
the following relations hold true.
\begin{equation}\label{eq:right_t-module}
\includegraphics{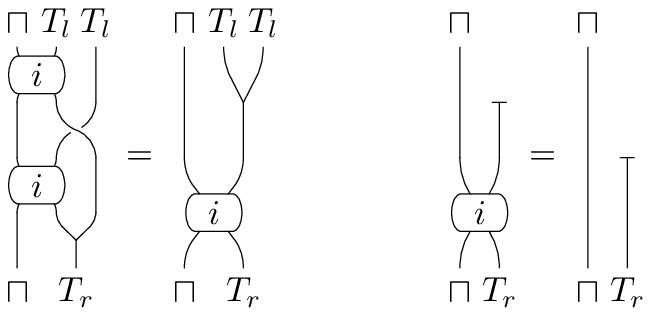}
\end{equation}
A left $\Phi$-module functor is a pair, consisting of a functor
$\sqcup:{\mathcal D}\to \M$  (where ${\mathcal D}$ is any category)
 and a natural transformation $w:T_r \sqcup \to T_l \sqcup$, such that
 the following relations hold true.
\begin{equation}\label{eq:left_t-module}
\includegraphics{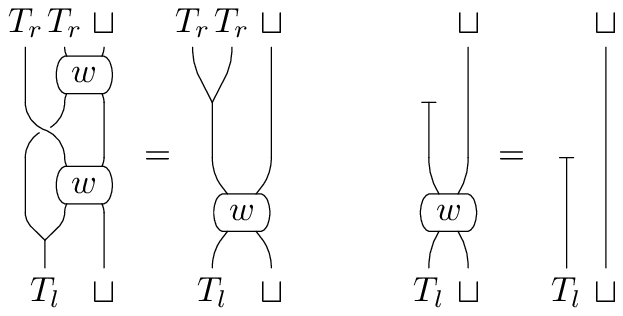}
\end{equation}

Once again, right or left $\Phi$-module functors can be described as
2-functors from an appropriately chosen 2-category to $\mathsf{Cat}$. Consider
a 2-category ${\mathcal R}$, with two 0-cells $o$ and $c$, generated by two
monads $(t_l,\mu_l,\eta_l)$ and $(t_r,\mu_r,\eta_r)$ on $o$ related by a
distributive law $\phi:t_r t_l \to t_l t_r$, and a further 1-cell $p:o\to c$
together with a 2-cell $\iota:p t_l \to p t_r$ such that
$$
\iota \circ p \mu_l = p \mu_r \circ \iota t_r\circ p \phi \circ \iota t_l
\qquad \textrm{and}\qquad
\iota \circ p \eta_l = p \eta_r,
$$
(for a diagrammatic representation see \eqref{eq:right_t-module}) \black
i.e. such that $(p,\iota)$ is an ${\mathcal R}(\phi,c)$-algebra
in ${\mathcal R}(o,c)$. A 2-functor $K$ from ${\mathcal R}$ to $\mathsf{Cat}$
can be described then as a pair of monads on the category $Ko$, related by a
distributive law $K\phi$, together with a right $K\phi$-module functor
$(Kp,K\iota)$. A 2-functor from the horizontal opposite of ${\mathcal R}$ to
$\mathsf{Cat}$ corresponds to a pair of monads related by a
distributive law and a left module functor for it.

Making use of the above observations, we can define morphisms between right
or left $\Phi$-module functors as lax natural transformations between the
corresponding 2-functors. Explicitly, this yields the following.

\begin{definition}\label{def:t-module_morphism}
Consider two monads $(T_l,m_l,u_l)$ and $(T_r,m_r,u_r)$ on the same category
$\M$ and two monads $(T'_l,m'_l,u'_l)$ and $(T'_r,m'_r,u'_r)$ on
$\M'$. Let $\Phi:T_r T_l \to T_l T_r$ and $\Phi':T'_r T'_l \to T'_l
T'_r$ be distributive laws.

A {\em morphism}
 from a right $\Phi'$-module functor
$(\sqcap':\M'\to {\mathcal C}', i':\sqcap' T'_l \to \sqcap' T'_r)$ to a
right $\Phi$-module functor
$(\sqcap:\M\to {\mathcal C},i:\sqcap T_l \to \sqcap T_r)$
is a quintuple
$(G,q_l,q_r,\L,\pi)$,
where
$(G,q_l): (T'_l,m'_l,u'_l)\to (T_l,m_l,u_l) $ and
$(G,q_r): (T'_r,m'_r,u'_r)\to (T_r,m_r,u_r) $
are monad morphisms,
$\L:{\mathcal C}'\to {\mathcal C} $
is a functor and
$\pi:  \sqcap G \to \L \sqcap' $
is a natural transformation, such that \eqref{eq:d-law_comp}  and the
following relation hold.
\begin{equation}\label{eq:right_t-mod_morphism}
\includegraphics{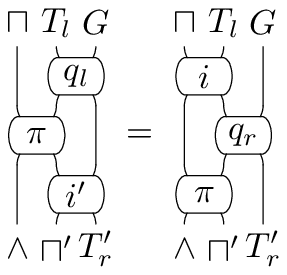}
\end{equation}

Similarly, a {\em morphism}
 from a left $\Phi'$-module functor
$(\sqcup':{\mathcal D}'\to \M', w':T'_r \sqcup'\to T'_l \sqcup')$
to a left $\Phi$-module functor
$(\sqcup:{\mathcal D}\to \M,w:T_r \sqcup\to T_l \sqcup)$
is a quintuple
$(G,q_l,q_r,\V,\omega)$,
where
$(G,q_l): (T'_l,m'_l,u'_l)\to (T_l,m_l,u_l) $ and
$(G,q_r): (T'_r,m'_r,u'_r)\to (T_r,m_r,u_r) $
are monad morphisms,
$\V: {\mathcal D}'\to {\mathcal D} $
is a functor and
$\omega: \sqcup \V \to G \sqcup' $
is a natural transformation, such that
\eqref{eq:d-law_comp} and the following relation hold.
\begin{equation}\label{eq:left_t-mod_morphism}
\includegraphics{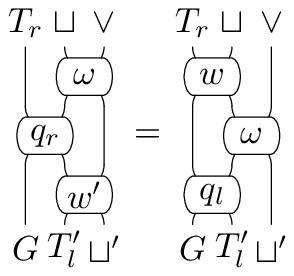}
\end{equation}
\end{definition}

 It is immediately clear by their definition as lax natural
transformations that
morphisms of (left or right) $\Phi$-module functors can be composed in the
appropriate sense.


\section{Para-cocyclic objects and $\Phi$-module functors}\label{sec:cat_A}

In this section we construct a category ${\mathcal A}$,  which  comes
equipped with a functor to a category of para-cocyclic objects in
the category of functors. This implies  that any object of
${\mathcal A}$ induces a functor from a category ${\mathcal D}$
(of coefficients) to the category of para-cocyclic objects in some
category ${\mathcal C}$.

Motivated by the constructions in Section \ref{sec:prelims}, consider a
2-category ${\mathcal S}$ with three 0-cells $o$, $c$ and $d$, generated by
two monads on $o$ related by a distributive law $\phi$ and an
 ${\mathcal S}(\phi,c)$-algebra  in ${\mathcal S}(o,c)$ and an
 ${\mathcal S}(d,\phi)$-algebra  in ${\mathcal S}(d,o)$.
Explicitly, ${\mathcal S}$ has 1-cells depicted in
$$
\includegraphics{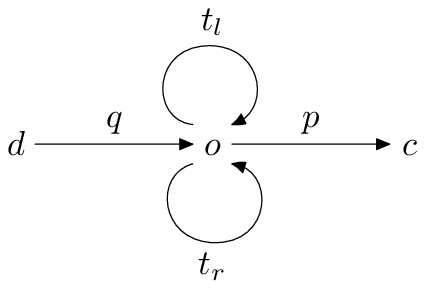}
$$
2-cells are horizontal and vertical composites of identity 2-cells and the
following 2-cells.
\begin{gather*}
\mu_l:t_l^2 \Rightarrow t_l\qquad
\eta_l:o \Rightarrow t_l \qquad
\mu_r:t_r^2 \Rightarrow t_r\qquad
\eta_r:o \Rightarrow t_r
\\
\phi:t_rt_l \Rightarrow t_lt_r\qquad
\iota:pt_l \Rightarrow pt_r \qquad
\varpi: t_r q \Rightarrow t_l q
\end{gather*}
On these generating 2-cells one imposes three types of relations. The first one
\[
\includegraphics{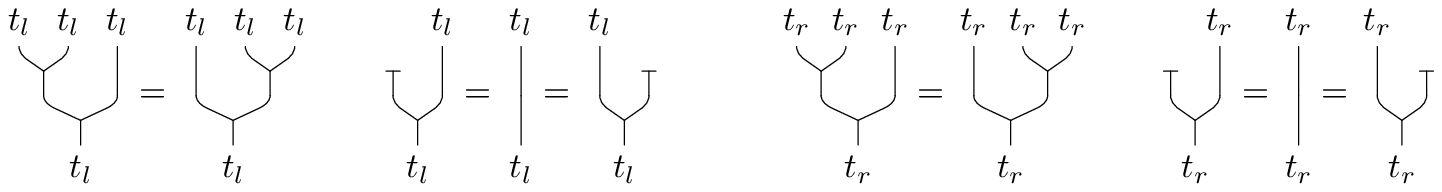}
\]
means that $(t_l,\mu_l,\eta_l)$ and $(t_r,\mu_t,\eta_t)$ are monads in
$\mathcal S$. One also imposes the relations
\[
\includegraphics{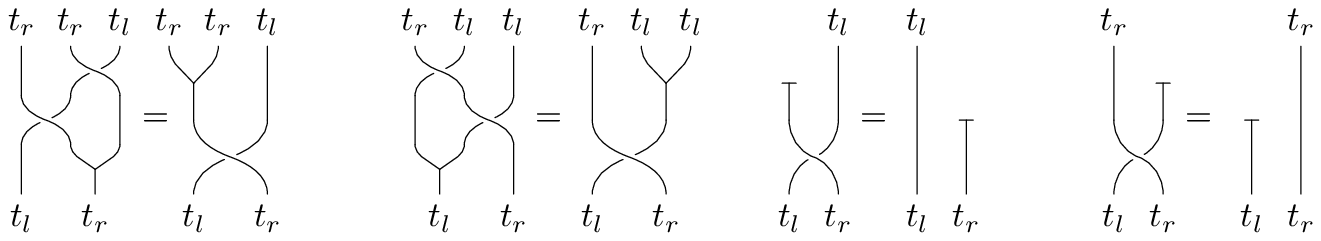}
\]
so that $\phi$ is a distributive law. Finally, the relations
\[
\includegraphics{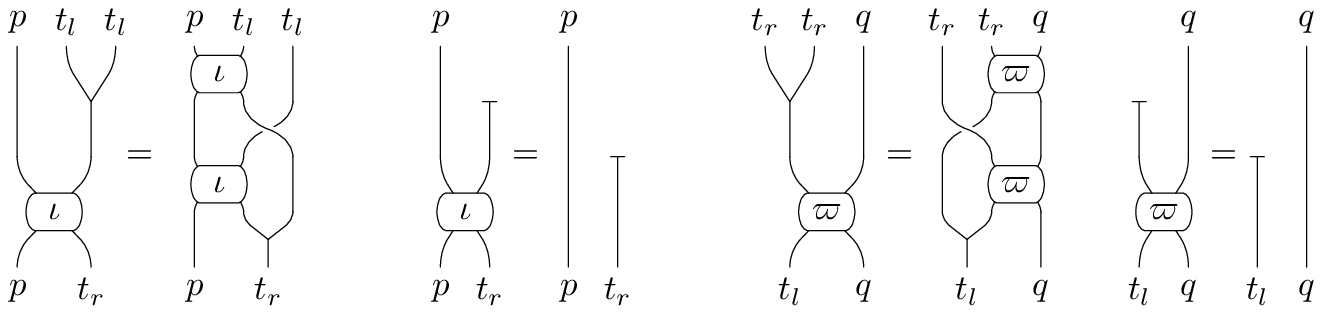}
\]
mean that $(p,\iota)$ is ${\mathcal S}(\phi,c)$-algebra  in ${\mathcal
S}(o,c)$ and $(q,\varpi)$ is an ${\mathcal S}(d,\phi)$-algebra  in ${\mathcal
S}(d,o)$.\black

We define a category ${\mathcal A}$ as the {\em opposite} of the category of
2-functors from ${\mathcal S}$ to $\mathsf{Cat}$ and lax natural
transformations between them. Explicitly, this means the following.

\begin{definition} \label{def:cat_A}
The category ${\mathcal A}$ is defined to have {\em objects}
$(T_l,T_r,\Phi,\sqcap,i,\sqcup,w)$, where
\begin{itemize}
\item $T_l$ and $T_r$ are monads on the same category $\M$;
\item $\Phi:T_r T_l \to T_l T_r$ is a distributive law;
\item $(\sqcap:\M\to {\mathcal C}, i:\sqcap T_l \to \sqcap T_r)$ is
  a right $\Phi$-module functor;
\item $(\sqcup:{\mathcal D} \to \M,w:T_r \sqcup \to T_l \sqcup)$ is a
  left $\Phi$-module functor.
\end{itemize}
A {\em morphism}
$(T_l,T_r,\Phi,\sqcap,i,\sqcup,w) \to (T'_l,T'_r,\Phi',\sqcap',i',\sqcup',w')$
is a datum $(G,q_l,q_r,\L,\pi,$ $\V,\omega)$, such that
\begin{itemize}
\item $(G,q_l,q_r,\L,\pi)$ is a morphism from the right
 $\Phi'$-module functor $(\sqcap', i')$
to the right
 $\Phi$-module functor $(\sqcap,i)$ (cf. Definition
\ref{def:t-module_morphism}); \black
\item $(G,q_l,q_r,\V,\omega)$ is a morphism from the left
 $\Phi'$-module functor $(\sqcup', w')$
to the left
 $\Phi$-module functor $(\sqcup,w)$ (cf. Definition
\ref{def:t-module_morphism}). \black
\end{itemize}
\end{definition}

Recall that a {\em para-cocyclic object} in a category ${\mathcal C}$
consists of a family $\{Z^n\}_n$ of objects in ${\mathcal C}$, for all
non-negative integers $n$, and morphisms
$$
d^k:Z^{n-1} \to Z^n,\qquad s^k:Z^{n+1}\to Z^n, \qquad \textrm{for } k=0,\dots, n,
$$
called {\em coface} and {\em codegeneracy} morphisms, respectively, satisfying
cosimplicial relations, together with so called {\em para-cocyclic morphisms}
$t^n:Z^n \to Z^n$, for all $n\geq 0$, which satisfy, for all $k=1,\dots,n$,
$$
t^n \circ d^0=d^n,\quad
t^n \circ d^k = d^{k-1}\circ t^{n-1},\quad
t^n \circ s^0 =s^n\circ t^{n+1}\circ t^{n+1},\quad
t^n \circ s^k = s^{k-1}\circ t^{n+1}.
$$
A morphism $(Z^*,d^*,s^*,t^*)\to (Z^{\prime *},d'^*,s'^*,t'^*)$ is a family of
morphisms $\{f^n:Z^n \to Z^{\prime n}\}_{n\geq 0}$ in ${\mathcal C}$,
compatible with the coface, codegeneracy and para-cocyclic morphisms in the
evident sense.

\begin{definition}\label{def:cat_cocyc}
The category ${\overline{\mathcal P}}$ is defined to have {\em
objects} that are para-cocyclic objects in the category of
functors. That is, for any non-negative integer $n$, a functor
$Z^{n}:{\mathcal D}\to {\mathcal C}$ together with natural
transformations $d^k:Z^{n-1}\to Z^n$, $s^k:Z^{n+1}\to Z^n$,
$t^n:Z^n \to Z^n$, for $0\leq k \leq n$, satisfying the defining
relations of a para-cocyclic object.

{\em Morphisms} from $(Z^*:{\mathcal D}\to {\mathcal
C},d^*,s^*,t^*)$ to $(Z^{\prime *}:{\mathcal D}'\to {\mathcal
C}',d'^*,s'^*,t'^*)$ are triples $(\L,\V,\xi^*)$, where $ \L
:{\mathcal C}'\to {\mathcal C}$ and $  \V :{\mathcal D}'\to
{\mathcal D}$ are functors and $\xi^* : (Z^* \V,d^* \V,s^* \V,t^*
\V) \to (\L Z^{\prime *},\L d'^*, \L s'^*,\L t'^*)$  is a morphism
of para-cocyclic objects.
\end{definition}
In terms of $\Phi:T_rT_l\rightarrow T_lT_r$ and $q_l:T_lG\rightarrow
GT'_l$, we define inductively some new natural transformations:
Let $q_l^0$ be the identity natural transformation $G\to G$ and $\Phi^0$ be
the identity natural transformation $T_r\to T_r$. Put $q_l^1:=q_l$ and
$\Phi^1:=\Phi$. For every $n>1$ we now define $\Phi^n:T_rT_l^n\rightarrow
T_l^nT_r$ and $q_l^n:T_l^nG\rightarrow G{T'_l}^n$ by
$\Phi^{n}:={T_l}^{n-1}\Phi\circ \Phi^{n-1}T_l$ and
$q_l^{n}:=q_l{T'_l}^{n-1}\circ T_lq_l^{n-1}$, respectively. For these natural
transformation we will use the string diagrams in the figure below.
\[
\includegraphics{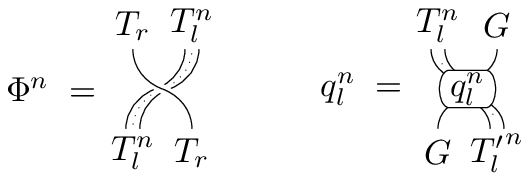}
\]
Note that in this figure each stripe represents a bunch of $n$ strings,
each one representing the functor $T_l$.
\begin{theorem}\label{thm:Z^*}
There is a functor ${\mathcal Z}^*:{\mathcal A} \to {\overline {\mathcal P}}$,
with object map
$$
(T_l, T_r, \Phi,\sqcap:\M\to {\mathcal C},i,\sqcup:{\mathcal D}\to {\mathcal
  M},w)\mapsto  \big(\ \sqcap T_l^{*+1} \sqcup:{\mathcal D}\to
  {\mathcal C},d^*,s^*,t^*
\ \big).
$$
The functor ${\mathcal Z}^*$ takes a morphism
$$
(G,q_l,q_r,\L,\pi,\V,\omega):(T_l,T_r,\Phi,\sqcap,i,\sqcup,w)\to
(T'_l,T'_r,\Phi',\sqcap',i',\sqcup',w')
$$
to the triple $(\L,\V,\xi^*)$. At every degree $n\geq 0$ and for $0\leq k\leq
n$,  the coface morphisms $d^k$, the codegeneracy morphisms $s^k$, the
para-cocyclic morphism $t^n$ and the morphism $\xi^n$ are given by the natural
transformations below.
\[
\includegraphics{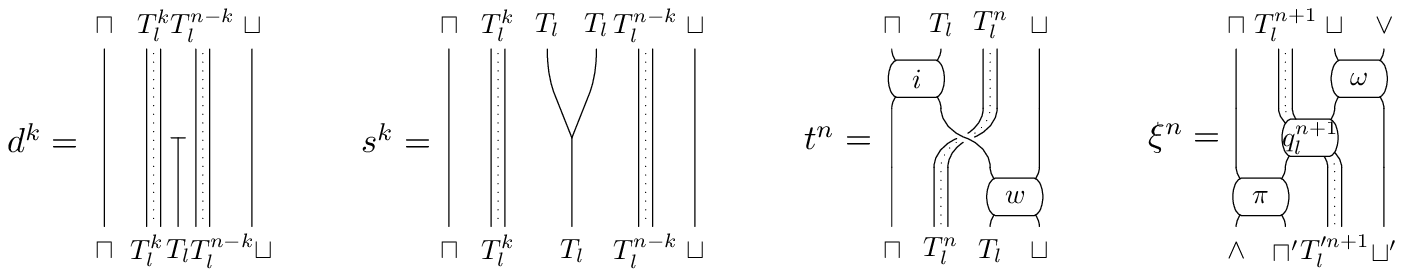}
\]
\end{theorem}

\begin{proof}
The datum $(\sqcap T_l^{*+1} \sqcup,d^*,s^*)$ is obviously a
cosimplex in the category of functors, cf. \cite{We}. Its
para-cocyclicity is checked with the same steps in \cite[Theorem
  1.10]{BS}. It remains to show that $\xi^*$ is a morphism of para-cocyclic
objects. Its compatibility with the coface and codegeneracy morphisms follows
by naturality and monad morphism property of $q_l$,  see the following
diagrammatic computations.
\[
\includegraphics{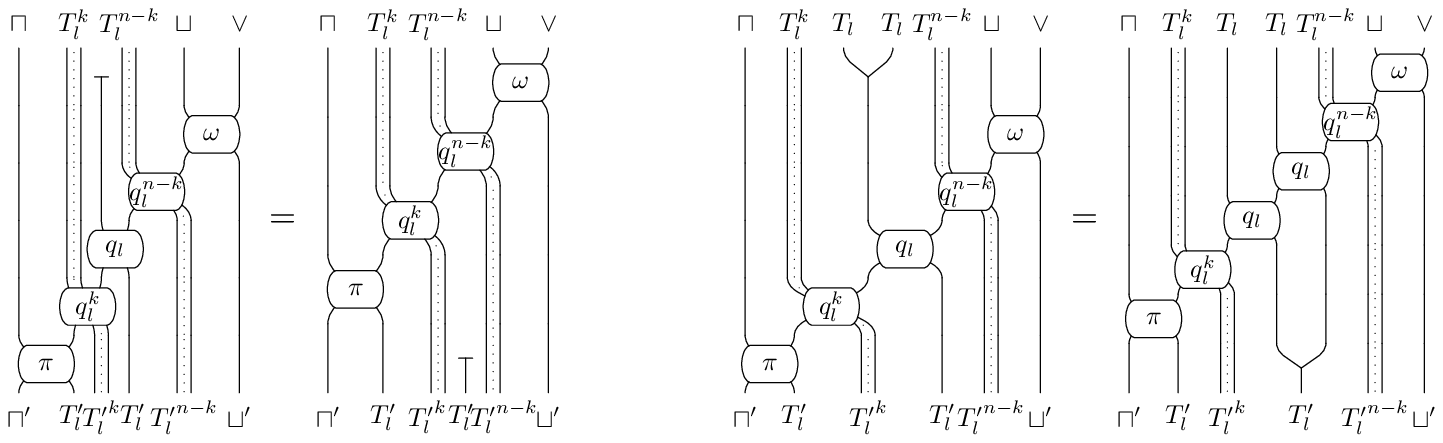}
\]
Compatibility with the para-cocyclic morphisms is proved in the figure below.
\[
\includegraphics{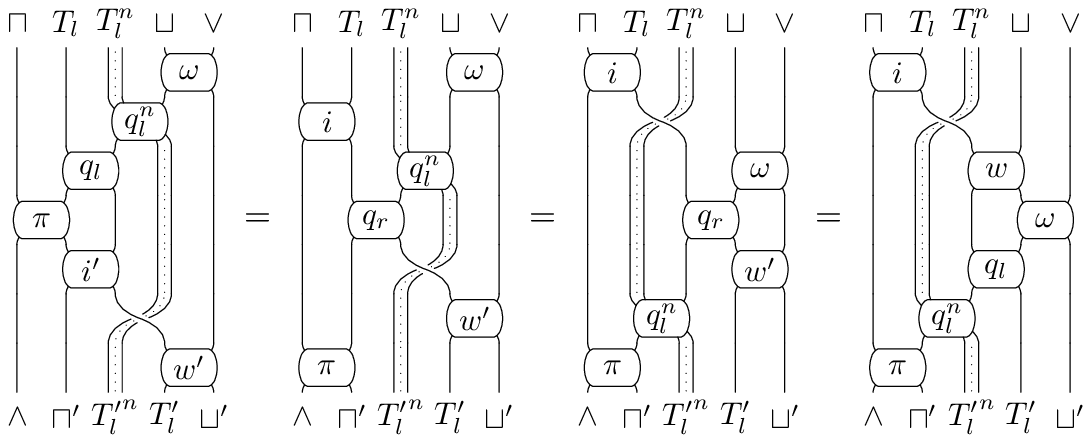}
\]
It follows by using naturality and
\eqref{eq:right_t-mod_morphism} in the first equality, then applying
\eqref{eq:d-law_comp} repeatedly in the second equality, and using in
the last equality \eqref{eq:left_t-mod_morphism} together with
naturality.
\end{proof}

\begin{corollary}
Any object $(T_l,T_r,\Phi,\sqcap:\M\to {\mathcal
  C},i,\sqcup:{\mathcal D}\to \M,w)$ of the category ${\mathcal A}$
determines a functor from ${\mathcal D}$ to the category of para-cocyclic
objects in ${\mathcal C}$. The objects of ${\mathcal D}$ play the role of
coefficients for the resulting para-cocyclic object in ${\mathcal C}$.
\end{corollary}


\section{Examples from Hopf cyclic theory}\label{sec:ex.A}

In this section we list some examples of objects in the category
${\mathcal A}$ in Definition \ref{def:cat_A}, arising from Hopf
cyclic theory (of bialgebroids, hence in particular of
bialgebras). They give rise to families of para-cocyclic objects
in the category $\mathsf{Mod}$-$k$ of modules over a commutative
ring $k$. They extend examples in \cite{JaSt:CycHom},
\cite{HaKhRaSo2}, \cite{Kay:BialgCycHom} and \cite{Brz:coef}.

Regarding the regular $R$-bimodule $R$ as a right $R^e:= R \ot R^{op}$-module,
and regarding any $R$-bimodule as a left $R^e$-module, we can define a
functor $R\ot_{R^e}(-):R$-$\mathsf{Mod}$-$R\to \mathsf{Mod}$-$k$. Applying it
to the $R$-module tensor product of two $R$-bimodules $M$ and $N$, it yields
the so called {\em cyclic $R$-module tensor product}. Throughout the paper, it
will be denoted by
$$
M {\widehat \ot}_R N := R \ot_{R^e} (M\ot_R N) \cong M\ot_{R^e} N.
$$
For finitely many $R$-bimodules $\{M_i\}_{i=1,\dots,n}$, we put
$$
M_1 {\widehat \ot}_R M_2 {\widehat \ot}_R \dots {\widehat \ot}_R M_n
:= (M_1 \ot_R \dots \ot_R M_{i-1}) \ot_{R^e} (M_i \ot_R \dots \ot_R
M_n),
$$
where the right hand side yields the same $k$-module for any $i=1,\dots, n$
(defining the $0$-fold tensor product to be equal to $R$).

For a short review of modules, comodules and contramodules of a
bialgebroid, we refer to the Appendix. Throughout, actions in
modules are denoted by juxtaposition and for coactions in
comodules we use a Sweedler type index notation, with implicit
summation understood.

The first example of an object in the category $\mathcal A$ in Definition
\ref{def:cat_A} arises from \cite[Theorem 2.4]{BS}.

\begin{example}\label{ex:1}
Let $B$ be a left bialgebroid over a $k$-algebra $L$ and $A$ be a left
$B$-module
algebra. Then $A$ is in particular an $L$-ring, with multiplication
$\mu:A\ot_L A \to A$ and unit $\eta:L \to A$.
An object in ${\mathcal A}$ is given by the following data.
\begin{itemize}
\item The monads $T_l= A \ot_L (-)$ and $T_r = (-) \ot_L A$ on
$L$-$\mathsf{Mod}$-$L$, with monad structures
\begin{eqnarray*}
&\mu\ot_L (-):T_l^2 \to T_l, \qquad
&\eta\ot_L (-) :L\textrm{-}\mathsf{Mod}\textrm{-}L \to T_l
\qquad \textrm{and}\\
&(-)\ot_L \mu:T_r^2 \to T_r, \qquad
&(-)\ot_L \eta :L\textrm{-}\mathsf{Mod}\textrm{-}L \to T_r,
\end{eqnarray*}
respectively;
\item The trivial distributive law $\Phi=A\ot_L (-)\ot_L A:T_r T_l \to T_l
  T_r$ ;
\item The right $\Phi$-module functor $(\sqcap,i)$, where
  $\sqcap=L\ot_{L^e} (-):L$-$\mathsf{Mod}$-$L \to \mathsf{Mod}$-$k$ and for
  any $L$-bimodule $P$,
$$
i_P: A {\widehat \ot}_L P \to P  {\widehat \ot}_L  A, \qquad a  {\widehat
  \ot}_L  p\mapsto p {\widehat \ot}_L  a;
$$
\item The left $\Phi$-module functor $(\sqcup,w)$, where
  $\sqcup:B$-$\mathsf{Comod} \to L$-$\mathsf{Mod}$-$L$ is the forgetful functor
  and for any left $B$-comodule $M$, with coaction $m\mapsto m_{[-1]} \ot_L
  m_{[0]}$,
$$
w_M: M \ot_L A \to A \ot_L M, \qquad m\ot_L a\mapsto m_{[-1]} a
\ot_L m_{[0]}.
$$
\end{itemize}
Applying the functor ${\mathcal Z}^*$ in Theorem \ref{thm:Z^*}, we obtain a
para-cocyclic object in $\mathsf{Mod}$-$k$, for any left $B$-comodule $M$. At
degree $n$, it is the $k$-module $A^{\otimes_L\, n+1} {\widehat\ot}_L
M$. Coface and codegeneracy maps are
\begin{eqnarray*}
&&d^k(a_0\cten L a_1 \cten L \dots \cten L a_{n-1}
\cten L m)= a_0 \cten L \dots \cten L a_{k-1}
\cten L 1_A \cten L a_k \cten L \dots {\widehat
  \ot}_L a_{n-1} \cten L m,\\
&&s^k(a_0\cten L a_1 \cten L \dots \cten L a_{n+1}
\cten L m)= a_0 \cten L \dots \cten L a_{k-1}
\cten L a_k a_{k+1} \cten L a_{k+2} \cten L \dots
\cten L a_{n+1} \cten L m,
\end{eqnarray*}
for $k=0,\dots,n$. The para-cocyclic operator comes out as
$$
t^n(a_0{\widehat \ot}_L a_1 {\widehat \ot}_L \dots {\widehat \ot}_L a_n
{\widehat \ot}_L m)=a_1 {\widehat \ot}_L \dots {\widehat \ot}_L a_n  {\widehat
  \ot}_L m_{[-1]} a_0 {\widehat \ot}_L m_{[0]}.
$$
\end{example}

The next example is obtained from \cite[Theorem 2.7]{BS}.

\begin{example}\label{ex:2}
Let $B$ be a right bialgebroid over a $k$-algebra $R$ and $(A,\mu,\eta)$ be a
right $B$-comodule algebra (hence in particular an $R$-ring). For the
$B$-coaction on $A$, use the notation $a\mapsto a^{[0]}\ot_R a^{[1]}$. An
object in ${\mathcal A}$ is given by the following data.
\begin{itemize}
\item The same monads $T_l= A \ot_R (-)$ and $T_r = (-) \ot_R A$ on
$R$-$\mathsf{Mod}$-$R$, introduced in Example \ref{ex:1} (replacing $L$ by $R$);
\item The trivial distributive law $\Phi=A\ot_R (-)\ot_R A:T_r T_l \to T_l
  T_r$ ;
\item The same right $\Phi$-module functor $(\sqcap,i)$, introduced in Example
  \ref{ex:1} (replacing $L$ by $R$);
\item The left $\Phi$-module functor $(\sqcup,w)$, where
  $\sqcup:\mathsf{Mod}$-$B \to R$-$\mathsf{Mod}$-$R$ is the forgetful functor
  and for any right $B$-module $N$,
$$
w_N: N \ot_R A \to A \ot_R N, \qquad  m \ot_R a\mapsto
a^{[0]}\ot_R  m  a^{[1]}.
$$
\end{itemize}
For any right $B$-module $N$, the corresponding para-cocyclic object in
$\mathsf{Mod}$-$k$ is, at degree $n$, $A^{\otimes_R\, n+1} {\widehat\ot}_R
N$. Coface and codegeneracy maps are given by the same formulae in Example
\ref{ex:1} (replacing $L$ by $R$ and $M$ by $N$). The para-cocyclic operator
has the form
$$
t^n(a_0{\widehat \ot}_R a_1 {\widehat \ot}_R \dots {\widehat
\ot}_R a_n {\widehat \ot}_R  m )=a_1 {\widehat \ot}_R \dots
{\widehat \ot}_R a_n  {\widehat
  \ot}_R a_0^{[0]} {\widehat \ot}_R   m   a_0^{[1]}.
$$
\end{example}

In the following example, for a left $R$-module $P$, a right $R$-module $Q$
and $R$-bimodules $C$ and $D$, the canonical isomorphisms
\begin{eqnarray*}
&&\Hom_{-,R}(C,\Hom_{-,R}(D,Q))\cong \Hom_{-,R}(C\ot_R D, Q)
\qquad \textrm{and} \\
&&
\Hom_{R,-}(C,\Hom_{R,-}(D,P))\cong \Hom_{R,-}(D\ot_R C, P)
\end{eqnarray*}
are suppressed.

\begin{example}\label{ex:3}
Let $B$ be a right bialgebroid over a $k$-algebra $R$ and
$(C,\Delta,\epsilon)$ be a right
$B$-module coring (hence in particular an $R$-coring).
An object in ${\mathcal A}$ is given by the following data.
\begin{itemize}
\item The monads $T_l:=\Hom_{-,R}(C,-)$ and $T_r:= \Hom_{R,-}(C,-)$ on
  $R$-$\mathsf{Mod}$-$R$, with monad structures
\begin{eqnarray*}
&\Hom_{-,R}(\Delta,-):T_l^2 \to T_l \ \qquad
&\Hom_{-,R}(\epsilon,-): R \textrm{-}\mathsf{Mod}\textrm{-}R \to T_l
\qquad \textrm{and}\\
&\Hom_{R,-}(\Delta,-):T_r^2 \to T_r \qquad
&\Hom_{R,-}(\epsilon,-): R \textrm{-}\mathsf{Mod}\textrm{-}R \to T_r,
\end{eqnarray*}
respectively;
\item The distributive law
$$
\Phi:\Hom_{R,-}(C,\Hom_{-,R}(C,-)) \cong
\Hom_{R,R}(C \ot_k C,-)\cong
\Hom_{-,R}(C,\Hom_{R,-}(C,-)),
$$
given by switching the arguments;
\item The right $\Phi$-module functor $(\sqcap,i)$, where
  $\sqcap=\Hom_{R,R}(R,-):R$-$\mathsf{Mod}$-$R\to \mathsf{Mod}$-$k$ and
$$
i_P: \Hom_{R,R}(R,\Hom_{-,R}(C,P))\cong \Hom_{R,R}(C,P)\cong
\Hom_{R,R}(R,\Hom_{R,-}(C,P))
$$
is the  hom-tensor  adjunction natural isomorphism, for any
$R$-bimodule $P$; \item The left $\Phi$-module functor
$(\sqcup,w)$, where
  $\sqcup:B$-$\mathsf{Ctrmod}\to R$-$\mathsf{Mod}$-$R$ is the forgetful
  functor and for a left $B$-contramodule $(Q,\alpha)$,
$$
w_Q: \Hom_{R,-}(C,Q) \to \Hom_{-,R}(C,Q),\qquad
f \mapsto \big(\ c \mapsto \alpha(f(c-)) \ \big).
$$
\end{itemize}
For any left $B$-contramodule $(Q,\alpha)$, this yields a para-cocyclic object
in $\mathsf{Mod}$-$k$. It is given by $\Hom_{R,R}(C^{\ot_R\,
  n+1},Q)$, at degree $n$. Coface and codegeneracy maps are
\begin{eqnarray*}
&&(d^k\varphi^{(n-1)})(c_0\stac R c_1 \stac R \dots \stac R c_n)=
\varphi^{(n-1)}(c_0\stac R \dots \stac R c_{k-1} \epsilon(c_k)\stac R
c_{k+1}\stac R \dots \stac R c_n),\\
&&(s^k\varphi^{(n+1)})(c_0\stac R c_1 \stac R \dots \stac R c_n)=
\varphi^{(n+1)}(c_0\stac R \dots \stac R c_{k-1} \stac R \Delta(c_k)\stac R
c_{k+1}\stac R \dots \stac R c_n),
\end{eqnarray*}
for $\varphi^{(j)}\in \Hom_{R,R}(C^{\ot_R j+1},Q)$ and $k=0,\dots,n$.
The para-cocyclic map is equal to
$$
(t^n\varphi^{(n)})(c_0\ot_R c_1 \ot_R \dots \ot_R c_n)=
\alpha \big(\varphi^{(n)}(c_n(-)\ot_R c_0\ot_R c_1\ot_R \dots \ot_R c_{n-1})\big).
$$
\end{example}

\begin{example}\label{ex:4}
Let $B$ be a left bialgebroid over a $k$-algebra $L$ and $(C,\Delta,\epsilon)$
be a left
$B$-comodule coring (hence in particular an $L$-coring), with $B$-coaction
$c\mapsto c_{[-1]}\ot_L c_{[0]}$.
An object in ${\mathcal A}$ is given by the following data.
\begin{itemize}
\item The same monads $T_l:=\Hom_{-,L}(C,-)$ and $T_r:= \Hom_{L,-}(C,-)$ on
  $L$-$\mathsf{Mod}$-$L$ introduced in Example \ref{ex:3} (replacing $R$ by
  $L$);
\item The same distributive law $\Phi$ introduced in Example \ref{ex:3}
  (replacing $R$ by $L$);
\item The same right $\Phi$-module functor $(\sqcap,i)$ introduced in Example
  \ref{ex:3} (replacing $R$ by $L$);
\item The left $\Phi$-module functor $(\sqcup,w)$, where
  $\sqcup:\mathsf{Mod}$-$B\to L$-$\mathsf{Mod}$-$L$ is the forgetful functor
  and for any right $B$-module $N$,
$$
w_N:\Hom_{L,-}(C,N)\to \Hom_{-,L}(C,N),\qquad f \mapsto \big(\ c\mapsto
f(c_{[0]})c_{[-1]}\ \big).
$$
\end{itemize}
The cosimplicial structure of the para-cocyclic object in $\mathsf{Mod}$-$k$,
corresponding to a right $B$-module $N$, is the same in Example
\ref{ex:3} (replacing $R$ by $L$ and $Q$ by $N$). The para-cocyclic map comes
out as
$$
(t^n\varphi ^{(n)})(c_0\ot_L  c_1 \ot_L  \dots \ot_L c_n)=
\varphi^{(n)}(c_{n[0]}\ot_L c_0\ot_L c_1 \ot_L \dots \ot_L
c_{n-1}) c_{n[-1]}.
$$
\end{example}
Specializing the above four examples to {\em bialgebras} instead of
bialgebroids, in all of them the functors $\sqcap$ become identity functors.

For an $R$-coring $C$, a {\em $C$-bicomodule} is an $R$-bimodule $M$, together
with a right $C$-coaction $\varrho^M:M \to M\ot_R C$ and a left $C$-coaction
${}^M\! \varrho:M\to C \ot_R M$, such that $\varrho^M$ is a left
$R$-module map, ${}^M\!\varrho$ is a right $R$-module map and $({}^M\!
\varrho \ot_R C) \circ \varrho^M = (C\ot_R \varrho^M) \circ {}^M\!
\varrho$. Morphisms of bicomodules are right $C$-comodule maps as well as left
$C$-comodule maps. The category of $C$-bicomodules is denoted by
$C$-$\mathsf{Comod}$-$C$.

\begin{example}\label{ex:5.D}
Let $B$ be a right bialgebroid over a $k$-algebra $R$ and
$(C,\Delta,\epsilon)$ be a left $B$-comodule
coring, hence in particular an $L:=R^{op}$-coring. An object in
${\mathcal A}$ is given by the following data.
\begin{itemize}
\item The monads $T_l=C \ot_L (-)$ and $T_r=(-)\ot_L C$ on
$C$-$\mathsf{Comod}$-$C$. For a $C$-bicomodule $(M,{}^M \varrho,\varrho^M)$,
$T_l  M = C \ot_L M$ is a $C$-bicomodule via the left and right
coactions
$$
c\ot_L m \mapsto \Delta(c) \ot_L m
\qquad \textrm{and} \qquad
c\ot_L m \mapsto c \ot_L \varrho^M(m).
$$
The monad structure of $T_l$ is given by the multiplication and unit
$$
C \ot_L \epsilon \ot_L M  : T_l^2 M \to T_l M
\qquad \textrm{and}\qquad
{}^M \varrho: M \to T_l M.
$$
Symmetrically, $T_r M = M \ot_L C$ is a $C$-bicomodule via the left and right
coactions
$$
m \ot_L c \mapsto {}^M \varrho(m) \ot_L c
\qquad \textrm{and} \qquad
m\ot_L c \mapsto m\ot_L \Delta(c).
$$
The monad structure of $T_r$ is given by the multiplication and unit
$$
M\ot_L \epsilon \ot_L C : T_r^2 M \to T_r M
\qquad \textrm{and} \qquad
\varrho^M: M \to T_r M.
$$
\item The trivial distributive law $\Phi=C \ot_L (-)\ot_L C$;
\item The right $\Phi$-module functor $\sqcap: C
  \textrm{-}\mathsf{Comod}\textrm{-} C\to \mathsf{Mod}\textrm{-}k$, given by
  the equalizer
$$
\xymatrix{
\sqcap M \ar[r] &
L {\widehat \otimes}_L M
\ar@<2pt>[rr]^-{L {\hat{\otimes}}_L \varrho^M}
\ar@<-2pt>[rr]_-{L {\hat{\otimes}}_L {}^M \varrho}
&&
C {\widehat \otimes}_L M \cong M {\widehat \otimes}_L C
}
$$
for any $C$-bicomodule $M$. The natural transformation $i$ is given by the
isomorphism
$$
i:\sqcap T_l \cong L {\widehat{\otimes}}_L (-) \cong (-) {\widehat{\otimes}}_L
L \cong \sqcap T_r.
$$
\item The left $\Phi$-module functor $\sqcup=(-)\ot_L C
  :\mathsf{Mod}\textrm{-}B\to C
  \textrm{-}\mathsf{Comod}\textrm{-} C$. For the left $B$-coaction on $C$,
  introduce the index notation $c \mapsto c^{[-1]}\ot_R c^{[0]}$ and for the
  comultiplication in $C$ write $\Delta(c) = c\1 \ot_L c\2$.
For any right $B$-module $N$, $\sqcup N = N \ot_L C$ is a $C$-bicomodule via the
left and right coactions
$$
m \ot_L c \mapsto {c\1}^{[0]} \ot_L m {c\1}^{[-1]} \ot_L c\2
\qquad \textrm{and} \qquad
m \ot_L c \mapsto m \ot_L c\1\ot_L c\2.
$$
The natural transformation $w:(-)\ot_L C \ot_L C \to C\ot_L (-)\ot_L C$ is
given by
$$
w_N(m\ot_L c \ot_L d)= c^{[0]} \ot_L m c^{[-1]} \ot_L d.
$$
\end{itemize}
For any right $B$-module $N$, this yields a para-cocyclic object
in $\mathsf{Mod}\textrm{-}k$. At degree $n$, it is given by
$\sqcap (C^{\ot_L n+1}\ot_L N \ot_L C) \cong C^{\ot_L n+1}
{\widehat{\otimes}}_L N$. For every $0\leq k\leq n-1$, the
corresponding coface map is:
\begin{eqnarray*}
&&d^k(c_0 {\widehat{\otimes}}_L \dots {\widehat{\otimes}}_L
c_{n-1} {\widehat{\otimes}}_L m) = c_0 {\widehat{\otimes}}_L \dots
c_{k-1} {\widehat{\otimes}}_L \Delta(c_k) {\widehat{\otimes}}_L
c_{k+1} {\widehat{\otimes}}_L \dots {\widehat{\otimes}}_L c_{n-1}
{\widehat{\otimes}}_L m,
\end{eqnarray*}
while
\begin{eqnarray*}
&&d^n(c_0 {\widehat{\otimes}}_L \dots {\widehat{\otimes}}_L
c_{n-1} {\widehat{\otimes}}_L m) = {c_0} \2 {\widehat{\otimes}}_L
c_1 {\widehat{\otimes}}_L \dots {\widehat{\otimes}}_L c_{n-1}
{\widehat{\otimes}}_L {{c_0}\1}^{[0]} {\widehat{\otimes}}_L m
{{c_0}\1}^{[-1]}.
\end{eqnarray*}
If $0\leq k\leq n$, then the codegeneracy map $s^k$ is given
by
\begin{eqnarray*}
&&s^k(c_0 {\widehat{\otimes}}_L \dots {\widehat{\otimes}}_L
c_{n+1} {\widehat{\otimes}}_L m) = c_0 {\widehat{\otimes}}_L \dots
c_{k-1} {\widehat{\otimes}}_L c_k \epsilon(c_{k+1})
{\widehat{\otimes}}_L c_{k+2} {\widehat{\otimes}}_L \dots
{\widehat{\otimes}}_L c_{n+1} {\widehat{\otimes}}_L m.
\end{eqnarray*}
The para-cocyclic map is
$$
t^n(c_0 {\widehat{\otimes}}_L \dots {\widehat{\otimes}}_L c_n
{\widehat{\otimes}}_L m) =
c_1 {\widehat{\otimes}}_L \dots {\widehat{\otimes}}_L c_{n}
{\widehat{\otimes}}_L c_0^{[0]} {\widehat{\otimes}}_L m c_0^{[-1]}.
$$
\end{example}

\begin{example}\label{ex:6.D}
Let $B$ be a right bialgebroid over a $k$-algebra $R$ and
$(C,\Delta,\epsilon)$ be a right $B$-module coring, hence in particular an
$R$-coring. An object in ${\mathcal A}$ is given by the following data.
\begin{itemize}
\item The same monads $T_l= C \ot_R (-)$ and $T_r= (-)\ot_R C$ on $C
  \textrm{-}\mathsf{Comod}\textrm{-} C$, introduced in Example \ref{ex:5.D}
  (replacing $L$ by $R$);
\item The trivial distributive law $\Phi=C \ot_R (-) \ot_R C$;
\item The same right $\Phi$-module functor $(\sqcap, i)$, introduced in
  Example \ref{ex:5.D} (replacing $L$ by $R$);
\item The left $\Phi$-module functor $\sqcup=(-)\ot_R
  C:B\textrm{-}\mathsf{Comod}\to C \textrm{-}\mathsf{Comod}\textrm{-} C$. For
  a left $B$-comodule $M$, with coaction $m \mapsto m^{[-1]}\ot_R m^{[0]}$,
  using the notation $\Delta(c)= c\1\ot_R c\2$, $\sqcup M= M \ot_R C$ is a
  $C$-bicomodule with left and right coactions
$$
m\ot_R c \mapsto c\1 m^{[-1]} \ot_R m^{[0]} \ot_R c\2
\qquad \textrm{and} \qquad
m\ot_R c \mapsto m \ot_R c\1 \ot_R c\2.
$$
The natural transformation $w:(-)\ot_R C \ot_R C \to C \ot_R (-)\ot_R C$ is
given by
$$
w_N(m\ot_R c \ot_R d)= c m^{[-1]} \ot_R m^{[0]} \ot_R d.
$$
\end{itemize}
For any left $B$-comodule $M$ this determines a para-cocyclic
object in $\mathsf{Mod}\textrm{-}k$. At degree $n$, it is given by
$\sqcap(C^{\ot_R n+1}\ot_R M \ot _R C) \cong C^{\ot_R n+1}
{\widehat{\otimes}}_R M$. For every $0\leq k\leq n-1$, the
corresponding coface map is:
\begin{eqnarray*}
&&d^k(c_0 {\widehat{\otimes}}_R \dots {\widehat{\otimes}}_R
c_{n-1} {\widehat{\otimes}}_R m) = c_0 {\widehat{\otimes}}_R \dots
c_{k-1} {\widehat{\otimes}}_R \Delta(c_k) {\widehat{\otimes}}_R
c_{k+1} {\widehat{\otimes}}_R \dots {\widehat{\otimes}}_R c_{n-1}
{\widehat{\otimes}}_R m,
\end{eqnarray*}
while
\begin{eqnarray*}
&&d^n(c_0 {\widehat{\otimes}}_R \dots {\widehat{\otimes}}_R
c_{n-1} {\widehat{\otimes}}_R m) = {c_0} \2 {\widehat{\otimes}}_R
c_1 {\widehat{\otimes}}_R \dots {\widehat{\otimes}}_R c_{n-1}
{\widehat{\otimes}}_R {{c_0}\1} m^{[-1]} {\widehat{\otimes}}_R
m^{[0]}.
\end{eqnarray*}
If $0\leq k\leq n$, then the codegeneracy map $s^k$ is given
by
\begin{eqnarray*}
&&s^k(c_0 {\widehat{\otimes}}_R \dots {\widehat{\otimes}}_R
c_{n+1} {\widehat{\otimes}}_R m) = c_0 {\widehat{\otimes}}_R \dots
c_{k-1} {\widehat{\otimes}}_R c_k \epsilon(c_{k+1})
{\widehat{\otimes}}_R c_{k+2} {\widehat{\otimes}}_R \dots
{\widehat{\otimes}}_R c_{n+1} {\widehat{\otimes}}_R m.
\end{eqnarray*}
The para-cocyclic map is
$$
t^n(c_0 {\widehat{\otimes}}_R \dots {\widehat{\otimes}}_R c_{n}
{\widehat{\otimes}}_R m) =
c_1 {\widehat{\otimes}}_R \dots {\widehat{\otimes}}_R c_{n}
{\widehat{\otimes}}_R c_0 m^{[-1]}{\widehat{\otimes}}_R m^{[0]}.
$$
Restricting to the case when $B$ is a bialgebra over $k$, this para-cocyclic
module yields a symmetrical version of \cite[(2.1)-(2.4)]{HaKhRaSo2}
(note the minor difference of using a left or a right module coalgebra). In
\cite{HaKhRaSo2} additional assumptions are made on the left comodule $M$
under which a truly cocyclic quotient exists.
\end{example}

\begin{example}\label{ex:7.D}
Let $B$ be a left bialgebroid over a $k$-algebra $L$ and $(A,\mu,\eta)$ be a
left $B$-module algebra, hence in particular an $L$-ring. The following data
define an object in ${\mathcal A}$.
\begin{itemize}
\item The monads $T_l=\Hom_{-,L}(A,-)$ and $T_r =\Hom_{L,-}(A,-)$ on
$A\textrm{-}\mathsf{Mod}\textrm{-}A$. For any $A$-bimodule $X$,
$T_lX=\Hom_{-,L}(A,X)$ is an $A$-bimodule via
$$
(a_1 f a_2)(a) = a_1 f(a_2 a),\qquad \textrm{for }a,a_1,a_2\in A,\ f\in
  \Hom_{-,L}(A,X).
$$
The monad structure is given by
$$
\Hom_{-,L}(A\ot_L \eta,X):T_l^2 X \to T_l X
\qquad \textrm{and}\qquad
X \to T_l X,\quad x\mapsto [a\mapsto xa].
$$
Symmetrically, $T_r X =\Hom_{L,-}(A,X)$ is an $A$-bimodule via
$$
(a_1 g a_2)(a) = g(a a_1)a_2,\qquad \textrm{for }a,a_1,a_2\in A,\ g\in
  \Hom_{L,-}(A,X).
$$
The monad structure is given by
$$
\Hom_{L,-}(\eta \ot_L A,X):T_r^2 X \to T_r X
\qquad \textrm{and}\qquad
X \to T_r X,\quad x\mapsto [a\mapsto ax].
$$
\item The distributive law
$$
\Phi: \Hom_{-,L}(A,\Hom_{L,-}(A,-))\cong \Hom_{L,L}(A\ot_k A,-)\cong
\Hom_{L,-}(A,\Hom_{-,L}(A,-)).
$$
\item The right $\Phi$-module functor
$\sqcap:A\textrm{-}\mathsf{Mod}\textrm{-}A\to
\mathsf{Mod}\textrm{-}k$, given for any $A$-bimodule $X$ by the
equalizer
$$
\xymatrix{
\sqcap X \ar[r] &
\Hom_{L,L}(L,X)
\ar@<2pt>[rrr]^-{f\mapsto [a\mapsto af(1)]}
\ar@<-2pt>[rrr]_-{f\mapsto [a\mapsto f(1)a]} &&&
\Hom_{L,L}(A,X).
}
$$
That is, $\sqcap X$ is the center of the $A$-bimodule $X$. The
natural transformation $i$ is given by the isomorphism
$$
\sqcap \Hom_{-,L}(A,-)\cong \Hom_{L,L}(L,-)\cong \sqcap \Hom_{L,-}(A,-).
$$
\item The left $\Phi$-module functor
  $\sqcup=\Hom_{L,-}(A,-):B\textrm{-}\mathsf{Ctrmod} \to
  A\textrm{-}\mathsf{Mod}\textrm{-}A$. For any left $B$-contramodule
  $(Q,\alpha)$, $\sqcup Q = \Hom_{L,-}(A,Q)$ is an $A$-bimodule, via
$$
(a_1ga_2)(a) = \alpha \big(g((-a_2)a a_1)\big),
\qquad \textrm{for }a,a_1,a_2\in A,\ g\in \Hom_{L,-}(A,Q).
$$
$w:\Hom_{L,-}(A,\Hom_{L,-}(A,-))\to \Hom_{-,L}(A,\Hom_{L,-}(A,-))$
is given by
$$
\big(w_Q(h)\big)(a)(b)=\alpha\big(h(b)(-a)\big),
\quad \textrm{for }a,b\in A, \ h\in \Hom_{L,-}(A,\Hom_{L,-}(A,Q)).
$$
\end{itemize}
For any left $B$-contramodule $(Q,\alpha)$, this determines a
para-cocyclic object in $\mathsf{Mod}\textrm{-}k$. At degree $n$,
it is given by $\sqcap \Hom_{-,L}(A^{\ot_L n+1},\Hom_{L,-}(A,Q))
\cong \Hom_{L,L}(A^{\otimes_L n+1},Q)$. For every $0\leq k\leq
n-1$, the corresponding coface map is
 $$(d^k\varphi^{(n-1)})(a_0\stac L \dots \stac L a_n)=
\varphi^{(n-1)}(a_0\stac L \dots \stac L a_{k-1} \stac L a_k
a_{k+1}\stac L a_{k+2}\stac L \dots \stac L a_n),
 $$
while
 $$
(d^n \varphi^{(n-1)})(a_0\stac L \dots \stac L a_n)=
\alpha\big(\varphi^{(n-1)}((-a_{n})a_0 \stac L a_1 \stac L \dots
\stac L a_{n-1})\big).
 $$
If $0\leq k\leq n$, then the  codegeneracy map $s^k$ is given
by
\begin{eqnarray*}
&&(s^k \varphi^{(n+1)})(a_0\stac L \dots \stac L a_n)=
\varphi^{(n+1)}(a_0\stac L \dots \stac L a_k \stac L 1_A \stac L
a_{k+1} \stac L \dots \stac L a_n).
\end{eqnarray*}
The para-cocyclic map is
$$
(t^n \varphi^{(n)})(a_0\ot_L \dots \ot_L a_n)=
\alpha\big( \varphi^{(n)}(-a_n \ot_L a_0 \ot_L \dots\ot_L a_{n-1})\big),
$$
where $\varphi^{(j)}\in \Hom_{L,L}(A^{\ot_L j+1},Q)$.

This yields a non-commutative base version of the para-cocyclic module
in \cite[p. 6]{Brz:coef}.
(In \cite{Brz:coef} additional assumptions are made on the left contramodule
$(Q,\alpha)$ under which a truly cocyclic subobject exists.)
\end{example}

\begin{example}\label{ex:8.D}
Let $B$ be a left bialgebroid over a $k$-algebra $L$ and $(A,\mu,\eta)$ be a
right $B$-comodule algebra with $B$-coaction $a\mapsto a_{[0]}\ot_L
a_{[1]}$. Then $(A,\mu,\eta)$ is in particular an $R:=L^{op}$-ring. The
following data determine an object in ${\mathcal A}$.
\begin{itemize}
\item The same monads $T_l=\Hom_{-,R}(A,-)$ and $T_r=\Hom_{R,-}(A,-)$ on
  $A\textrm{-}\mathsf{Mod}\textrm{-}A$ introduced in Example \ref{ex:7.D}
  (replacing $L$ by $R$);
\item The same distributive law $\Phi$ introduced in Example \ref{ex:7.D}
  (replacing $L$ by $R$);
\item The same right $\Phi$-module functor $(\sqcap,i)$ introduced in Example
  \ref{ex:7.D} (replacing $L$ by $R$);
\item The left $\Phi$-module functor $\sqcup=\Hom_{R,-}(A,-):\mathsf{Mod}
  \textrm{-} B\to A\textrm{-}\mathsf{Mod} \textrm{-} A$. For any right
  $B$-module $N$, $\sqcup  N = \Hom_{R,-}(A,N)$ is an $A$-bimodule via
$$
(a_1 g a_2)(a) = g(a_{2[0]} a a_1) a_{2[1]},
\qquad \textrm{for }a,a_1,a_2 \in A,\ g\in \Hom_{R,-}(A,N).
$$
$w: \Hom_{R,-}(A,\Hom_{R,-}(A,-)) \to
\Hom_{-,R}(A,\Hom_{R,-}(A,-))$ is given by
$$
\big(w_N(h)\big)(a)(b)= h(b)(a_{[0]}) a_{[1]}.
$$
\end{itemize}
For any right $B$-module $N$, this determines a para-cocyclic
object in $\mathsf{Mod} \textrm{-}k$. At degree $n$, it is given
by $\sqcap \Hom_{-,R}(A^{\ot_R n+1},\Hom_{R,-}(A,N)) \cong
\Hom_{R,R}(A^{\otimes_R n+1},N)$. For every $0\leq k\leq n-1$,
the corresponding coface map is
 $$
 (d^k \varphi^{(n-1)})(a_0\stac R \dots \stac R a_n)=
\varphi^{(n-1)}(a_0\stac R \dots \stac R a_{k-1} \stac R a_k
a_{k+1}\stac R a_{k+2}\stac R  \dots \stac R a_{n}),
 $$
while
 $$(d^n
\varphi^{(n-1)})(a_0\stac R \dots \stac R a_n)=
\varphi^{(n-1)}(a_{n[0]}a_0 \stac R a_1 \stac R \dots \stac R
a_{n-1})a_{n[1]}.
 $$
If $0\leq k\leq n$, then the codegeneracy map $s^k$ is given
by
\begin{eqnarray*}
&&(s^k \varphi^{(n+1)})(a_0\stac R \dots \stac R a_n)=
\varphi^{(n+1)}(a_0\stac R \dots \stac R a_k \stac R 1_A \stac R
a_{k+1} \stac R \dots \stac R a_n).
\end{eqnarray*}
The para-cocyclic map is
$$
(t^n \varphi^{(n)})(a_0\ot_R \dots \ot_R a_n)=
\varphi^{(n)}(a_{n[0]} \ot_R a_0 \ot_R \dots\ot_R a_{n-1})a_{n[1]},
$$
where $\varphi^{(j)}\in \Hom_{R,R}(A^{\ot_R j+1},N)$.

This yields a non-commutative base version of the para-cocyclic module
in \cite[(3.1)--(3.4)]{HaKhRaSo2}.
(In \cite{HaKhRaSo2} additional assumptions are made on the right
module $N$ under which a truly cocyclic submodule exists.)
\end{example}


\section{Para-cyclic objects} \label{sec:cat_B}

Symmetrically to the considerations in Section \ref{sec:prelims} and Section
\ref{sec:cat_A}, one can obtain another category ${\mathcal B}$ together with
a functor ${\mathcal Z}_*$ from ${\mathcal B}$ to a category of para-cyclic
objects in the category of functors. Without repeating the details,  in
this section  we summarize the main steps.

\begin{definition}\label{def:comonad_morphism}
A {\em comonad} on a category $\M$ is a triple $(S,d,e)$, where
$S:\M\to \M$ is a functor and $d:S\to S^2$ and $e:S \to
\M$ are natural transformations called the {\em comultiplication} and
{\em counit}, respectively. Their string representations are shown in the
first two pictures of the figure below.
They satisfy the coassociativity and counitality constraints expressed by the
third and the fourth equalities in the same figure.

A {\em morphism} from a comonad $(S',d',e')$ on $\mathcal M'$ to a comonad
$(S,d,e)$ on $\mathcal M$ is a pair $(F,f)$, where  $F:\M'\to \M$
is a functor and  $f:SF \to FS'$ is a natural transformation which is
compatible with the comultiplications and the counits in the sense of
the last two relations of the following figure.
\[
\includegraphics{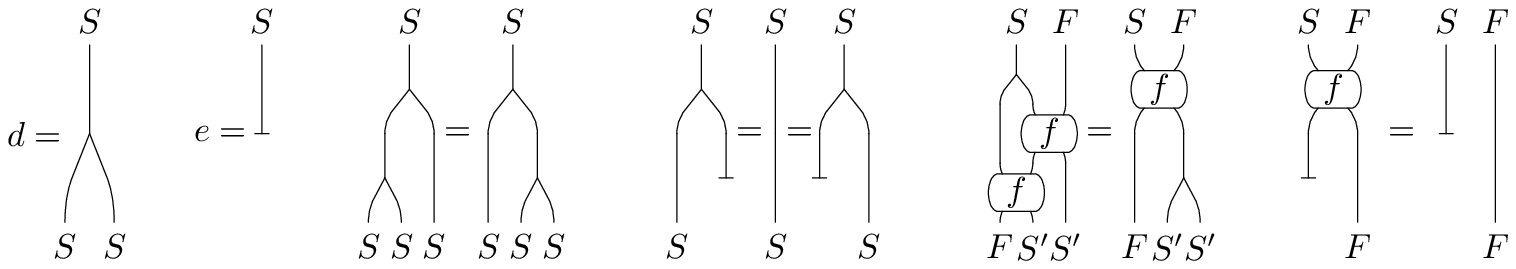}
\]
Comonads and their morphisms constitute a category which is isomorphic
to the
category of 2-functors and lax natural transformations, from
the vertical opposite of the 2-category ${\mathcal T}$ in Section
\ref{sec:prelims} to $\mathsf{Cat}$.

\end{definition}

\begin{definition}
A {\em coalgebra} for a comonad $S$ on a category $\M$ is a
pair $(M,\varrho)$, where $M$ is an object in $\M$ and
$\varrho:M\to SM$ is a morphism in $\M$ which is coassociative and
counital in the evident sense.

A {\em morphism} of $S$-coalgebras  $(M',\varrho')\to (M,\varrho)$
is a morphism  $\varphi:M'\to M$  in $\M$, such that
 $S\varphi \circ \varrho' = \varrho\circ \varphi$.
Coalgebras of a comonad $S$ and their morphisms
constitute the so-called Eilenberg-Moore category $\M^S$.
\end{definition}

Via composition on the right, a comonad $S:\M \to {\mathcal M}$
induces a comonad $\mathsf{Cat}(S,-)$  on the category $\mathsf{Cat}(\M,-)$.
Symmetrically, there is a comonad $\mathsf{Cat}(-,S)$ on the category
$\mathsf{Cat}(-,\M)$. We call a coalgebra of the comonad $\mathsf{Cat}(S,-)$
a {\em right $S$-comodule functor} and we term a coalgebra of the comonad
$\mathsf{Cat}(-,S)$ a {\em left $S$-comodule functor}.

\begin{definition}\label{def:comonad_distr_law}
Consider two comonads $(S_l,d_l,e_l)$ and $(S_r,d_r,e_r)$ on the same category
$\M$. A comonad {\em distributive law} is a natural transformation
$\Psi:S_l S_r \to S_r S_l$, such that the following equalities hold.
\[
\includegraphics{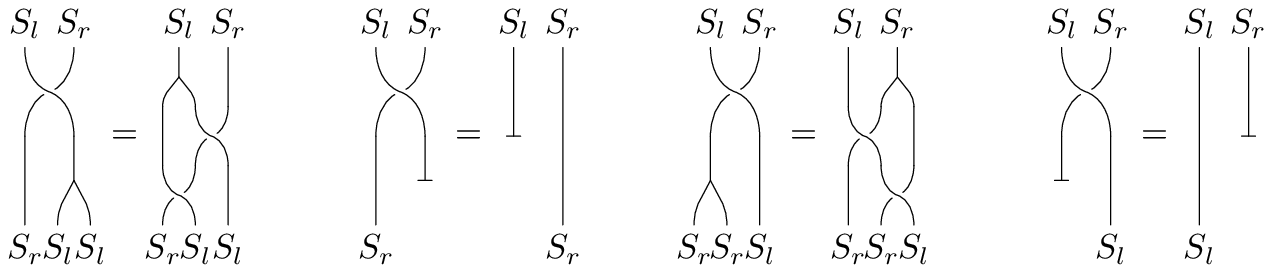}
\]
\end{definition}
A comonad distributive law $\Psi:S_l S_r \to S_r S_l$ as
in Definition \ref{def:comonad_distr_law} induces a comonad structure on the
composite functor $S_l S_r$, with comultiplication $d$ and counit $e$ whose
string representations are given in the figure below.
\[
\includegraphics{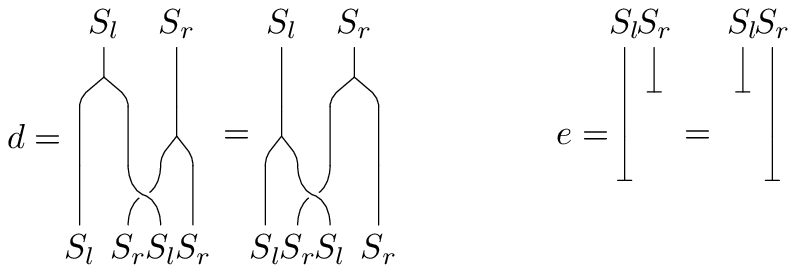}
\]
\begin{definition}\label{def:t-comodule}
Consider two comonads $(S_l,d_l,e_l)$ and $(S_r,d_r,e_r)$ on the same category
$\M$ and a comonad distributive law $\Psi:S_l S_r \to S_r S_l$.
A {\em $\Psi$-coalgebra} is a pair consisting of an object $X$ in ${\mathcal
  M}$ and a morphism $\xi:S_l X \to S_r X$ rendering commutative the following
diagrams.
\begin{equation}\label{eq:comod.i}
\xymatrix{
S_lX\ar[rrr]^-{\xi}\ar[d]_-{d_lX}&&&
S_rX\ar[d]^-{d_rX}\\
S_l^2X\ar[r]_-{S_l\xi}&
S_lS_rX\ar[r]_-{\Psi X}&
S_rS_lX\ar[r]_-{S_r\xi}&
S_r^2X
}
\qquad \textrm{and}\qquad
\xymatrix{
S_l X \ar[r]^-{\xi}\ar[d]_-{e_l X}&
S_r X \ar[d]^-{e_r X}\\
X\ar@{=}[r]&
X\ .
}
\end{equation}
\end{definition}

Coalgebras of the comonad distributive law $\mathsf{Cat}(\Psi,-)$  are called
{\em right $\Psi$-comodule functors}.
Consider two comonads $(S_l,d_l,e_l)$ and $(S_r,d_r,e_r)$ on the same category
$\M$ and two comonads $(S'_l,d'_l,e'_l)$ and $(S'_r,d'_r,e'_r)$ on
$\M'$. Let $\Psi:S_l S_r \to S_r S_l$ and $\Psi':S'_l S'_r \to S'_r
S'_l$ be comonad distributive laws. A {\em morphism}
 from a right $\Psi'$-comodule functor $(\sqcap':\M' \to {\mathcal C}',
i':\sqcap' S'_r \to \sqcap' S'_l)$ to a right $\Psi$-module functor
$(\sqcap:\M\to {\mathcal C},i:\sqcap S_r \to \sqcap S_l)$
is a quintuple $(G,q_l,q_r,\L,\pi)$, where $\L:{\mathcal C}' \to {\mathcal C}$
is a functor,  $\pi: \sqcap G \to \L \sqcap'$ is a natural transformations,
$(G,q_l): (S'_l,d'_l,e'_l)\to (S_l,d_l,e_l)$
and
$(G,q_r): (S'_r,d'_r,e'_r)\to (S_r,d_r,e_r)$
are comonad morphisms. These data are, in addition, subject to the following
two conditions.
\begin{equation}\label{eq:com_d-law_comp}
\includegraphics{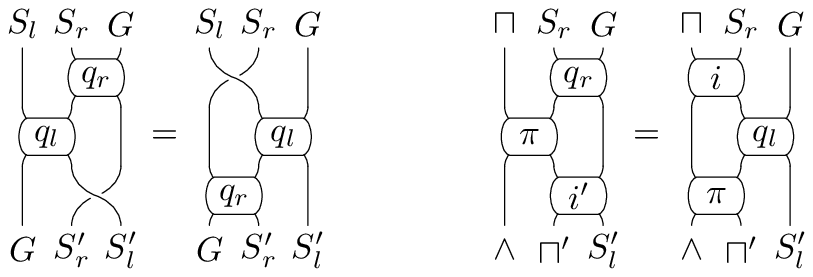}
\end{equation}
Right $\Psi$-comodule functors and their morphisms constitute a category
which is isomorphic to  the
category of 2-functors and
lax natural transformations, from the vertical opposite of the 2-category
${\mathcal R}$ in Section \ref{sec:prelims} to $\mathsf{Cat}$.

Symmetrically, a {\em left $\Psi$-comodule functor} is a coalgebra for the
comonad distributive law $\mathsf{Cat}(-,\Psi)$. A {\em morphism}
 from a left $\Psi'$-module functor $(\sqcup':{\mathcal D}' \to \M', w':S'_l
 \sqcup'\to S'_r \sqcup')$ to a left $\Psi$-module functor
$(\sqcup:{\mathcal D} \to \M,w:S_l \sqcup\to S_r \sqcup)$
is a quintuple
$(G,q_l,q_r,\V,\omega)$, where
$(G,q_l):  (S'_l,d'_l,e'_l)\to (S_l,d_l,e_l)$ and
$(G,q_r):  (S'_r,d'_r,e'_r)\to (S_r,d_r,e_r)$ are comonad morphisms
such that the first condition in \eqref{eq:com_d-law_comp} holds,
 $\V:{\mathcal D}'\to {\mathcal D}$ is a functor and  $\omega:\sqcup \V \to G
\sqcup'$
is a natural transformation such that the following relation holds.
\begin{equation}\label{eq:comod.omega}
\includegraphics{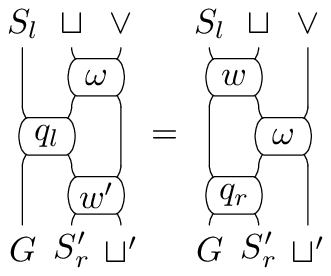}
\end{equation}
Left $\Psi$-comodule functors and their morphisms constitute a category
which is isomorphic to the
category of 2-functors and
lax natural transformations, from the horizontal and vertical opposite of the
2-category ${\mathcal R}$ in Section \ref{sec:prelims} to
$\mathsf{Cat}$.

We can define a category ${\mathcal B}$ as the opposite of the category
of 2-functors and lax natural transformations, from the vertical
opposite of the 2-category ${\mathcal S}$ in Section \ref{sec:prelims} to
$\mathsf{Cat}$:

\begin{definition} \label{def:cat_B}
The category ${\mathcal B}$ is defined to have {\em objects}
$(S_l,S_r,\Psi,\sqcap,i,\sqcup,w)$, where
\begin{itemize}
\item $S_l$ and $S_r$ are comonads on the same category $\M$;
\item $\Psi:S_l S_r \to S_r S_l$ is a comonad distributive law;
\item $(\sqcap:\M\to {\mathcal C}, i:\sqcap S_r \to \sqcap S_l)$ is
  a right $\Psi$-comodule functor;
\item $(\sqcup:{\mathcal D} \to \M,w:S_l \sqcup \to S_r \sqcup)$ is a
  left $\Psi$-comodule functor.
\end{itemize}
A {\em morphism}
$(S_l,S_r,\Psi,\sqcap,i,\sqcup,w) \to (S'_l,S'_r,\Psi',\sqcap',i',\sqcup',w')$
is a datum $(G,q_l,q_r,\L,\pi,$ $\V,\omega)$, such that
\begin{itemize}
\item $(G,q_l,q_r,\L,\pi)$ is a morphism from the right
 $\Psi'$-comodule functor $(\sqcap', i')$
to the right
 $\Psi$-comodule functor $(\sqcap,i)$;
\item $(G,q_l,q_r,\V,\omega)$ is a morphism from the left
 $\Psi'$-comodule functor $(\sqcup', w')$
to the left
 $\Psi$-comodule functor $(\sqcup,w)$.
\end{itemize}
\end{definition}

Recall that the {\em opposite} ${\mathcal C}^{op}$ of a category ${\mathcal
  C}$ has the same objects and morphisms as ${\mathcal C}$, but composition of
morphisms is opposite to that in ${\mathcal C}$. A para-cyclic object in a
category ${\mathcal C}$ is, by definition, a para-cocyclic object in
${\mathcal C}^{op}$.

\begin{definition}
{\em Objects} of
the category ${\underline{\mathcal P}}$ are para-cyclic objects in the
category of functors.
{\em Morphisms} from $(Z_*:{\mathcal D}\to {\mathcal C},d_*,s_*,t_*)$ to
$(Z'_*:{\mathcal D}'\to {\mathcal C}' ,d'_{*},s'_{*},t_{*}')$
are triples $(\L,\V,\xi_*)$, where
$\L:{\mathcal C}'\to {\mathcal C}$ and $\V:{\mathcal D}'\to {\mathcal D}$
are functors and $\xi_* :(Z_* \V,d_* \V,s_* \V,t_* \V) \to
({\L} Z'_*,\L d_{*}',\L s_{*}',\L t_{*}')$ is a morphism
of para-cyclic objects.
\end{definition}

Dually to Theorem \ref{thm:Z^*}, the following holds.

\begin{theorem}\label{thm:Z_*}
There is a functor ${\mathcal Z}_*:{\mathcal B} \to {\underline {\mathcal P}}$,
with object map
$$
(S_l,S_r,\Psi,\sqcap:\M\to {\mathcal C},i,\sqcup:{\mathcal D}\to {\mathcal
  M},w)\mapsto  \big(\ \sqcap S_l^{*+1} \sqcup:{\mathcal D}\to {\mathcal C},
d_*,s_*,t_* \ \big).
$$
The functor ${\mathcal Z}_*$ takes a morphism
$$
(G,q_l,q_r,\L,\pi,\V,\omega):(S_l,S_r,\Psi,\sqcap,i,\sqcup,w)\to
(S'_l,S'_r,\Psi',\sqcap',i',\sqcup',w')
$$
to the triple $(\L,\V,\xi_*)$. At every degree $n\geq 0$ and for $0\leq k\leq
n$,  the face morphisms $d_k$, the degeneracy morphisms $s_k$, the
para-cyclic morphism $t_n$ and the morphism $\xi_n$ are given by the natural
transformations below.
\[
\includegraphics{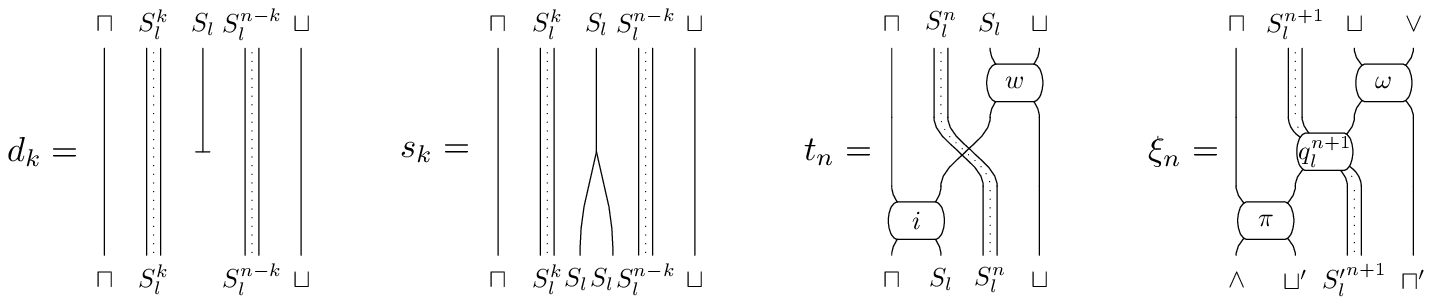}
\]
\end{theorem}

\begin{corollary}
Any object $(S_l,S_r,\Psi,\sqcap: \M\to {\mathcal
  C},i,\sqcup:{\mathcal D}\to \M,w)$ of the category ${\mathcal B}$
determines a functor from ${\mathcal D}$ to the category of para-cyclic
objects in ${\mathcal C}$.
The objects of ${\mathcal D}$ play the role of coefficients for the resulting
para-cyclic object in ${\mathcal C}$.
\end{corollary}

\section{Examples from Hopf cyclic theory}\label{sec:ex.B}

In this section we present several examples of objects in the
category ${\mathcal B}$ in Definition \ref{def:cat_B}, similar to
those we have seen in Section \ref{sec:ex.A}. Throughout this
section, the same notational conventions are used as in Section
\ref{sec:ex.A}.

The first example is obtained from \cite[Theorem 2.9]{BS} (for the
case of a bialgebra see \cite[Sec 5.3]{Kay:UniHCyc}).

\begin{example}\label{ex:5}
Let $B$ be a left bialgebroid over a $k$-algebra $L$ and $(C,\Delta,\epsilon)$
be a left $B$-comodule coring (hence in particular an $L$-coring), with
$B$-coaction $c \mapsto c_{[-1]}\ot_L c_{[0]}$. An object in the category
${\mathcal B}$ is given by the following data.
\begin{itemize}
\item The comonads $S_l= C \ot_L (-)$ and $S_r=(-)\ot_L C$ on
  $L$-$\mathsf{Mod}$-$L$, with comonad structures
\begin{eqnarray*}
&\Delta\ot_L (-): S_l \to S_l^2 \qquad
&\epsilon\ot_L (-):S_l \to L\textrm{-}\mathsf{Mod}\textrm{-}L
\qquad \textrm{and}\\
&(-)\ot_L \Delta: S_r \to S_r^2\qquad
&(-)\ot_L \epsilon:S_r \to L\textrm{-}\mathsf{Mod}\textrm{-}L,
\end{eqnarray*}
respectively;
\item The trivial comonad distributive law $\Psi= C \ot_L (-) \ot_L C$;
\item The right $\Psi$-comodule functor $(\sqcap,i)$, where $\sqcap=L
  \ot_{L^e} (-):L$-$\mathsf{Mod}$-$L\to \mathsf{Mod}$-$k$ and for any
  $L$-bimodule $P$,
$$
i_P: P {\widehat \ot}_L C \to C {\widehat \ot}_L P,\qquad
p{\widehat \ot}_L c \mapsto c {\widehat \ot}_L p;
$$
\item The left $\Psi$-comodule functor $(\sqcup,w)$, where
  $\sqcup:B$-$\mathsf{Mod}\to L$-$\mathsf{Mod}$-$L$ is the forgetful functor
  and for any left $B$-module $N$,
$$
w_N:C \ot_L N \to N\ot_L C,\qquad
c\ot_L m \mapsto c_{[-1]} m \ot_L c_{[0]}.
$$
\end{itemize}
For any left $B$-module $N$, the corresponding para-cyclic object in
$\mathsf{Mod}$-$k$ is given by $C^{{\widehat \ot}_L\, n+1} {\widehat \ot}_L
N$ at degree $n$. Face and degeneracy maps are,
for $k=0,\dots,n$,
\begin{eqnarray*}
&& d_k(c_0 \cten L c_1 \cten L \dots \cten L c_n
\cten L m) = c_0 \cten L c_1 \cten L \dots
\cten L c_{k-1} \cten L \epsilon(c_k) \cten L c_{k+1}
\cten L  \dots \cten L  c_n \cten L m \\
&&s_k(c_0 \cten L c_1 \cten L \dots \cten L c_n
\cten L m)= c_0 \cten L c_1 \cten L \dots
\cten L c_{k-1} \cten L \Delta (c_k) \cten L c_{k+1}
\cten L  \dots \cten L  c_n \cten L m.
\end{eqnarray*}
The para-cyclic map comes out as
$$
t_n(c_0 {\widehat \ot}_L c_1 {\widehat \ot}_L \dots {\widehat \ot}_L c_n
{\widehat \ot}_L m)= c_{n[0]} {\widehat \ot}_L c_0 {\widehat \ot}_L \dots
{\widehat \ot}_L c_{n-1} {\widehat \ot}_L c_{n[-1]} m.
$$
\end{example}

A next example arises from \cite[Theorem 2.11]{BS}. Restricting it to
the case of a bialgebra, it yields a symmetrical version of the para-cocyclic
module in \cite[Sec 5.1]{Kay:UniHCyc}.

\begin{example}\label{ex:6}
Let $B$ be a right bialgebroid over a $k$-algebra $R$ and
$(C,\Delta,\epsilon)$ be a right $B$-module coring (so in particular an
$R$-coring).
An object in ${\mathcal B}$ is given by the following data.
\begin{itemize}
\item The same comonads $S_l= C \ot_R (-)$ and $S_r=(-)\ot_R C$ on
  $R$-$\mathsf{Mod}$-$R$ introduced in Example \ref{ex:5} (replacing $L$ by
  $R$);
\item The trivial comonad distributive law $\Psi= C \ot_R (-) \ot_R C$;
\item The same right $\Psi$-comodule functor $(\sqcap,i)$ introduced in
  Example \ref{ex:5} (replacing $L$ by $R$);
\item The left $\Psi$-comodule functor $(\sqcup,w)$, where
  $\sqcup:\mathsf{Comod}$-$B\to R$-$\mathsf{Mod}$-$R$ is the forgetful functor
  and for any right $B$-comodule $M$, with coaction $m\mapsto m^{[0]}\ot_R
  m^{[1]}$,
$$
w_M:C \ot_R M \to M\ot_R C, \qquad c\ot_R m \mapsto m^{[0]}\ot_R c m^{[1]}.
$$
\end{itemize}
For any right $B$-comodule $M$, the corresponding para-cyclic module has the
same simplicial structure in Example \ref{ex:5} (replacing $L$ by $R$ and $N$
by $M$). The para-cyclic map is
$$
t_n(c_0 {\widehat \ot}_R c_1 {\widehat \ot}_R \dots {\widehat \ot}_R c_n
{\widehat \ot}_R m)= c_n m^{[1]}{\widehat \ot}_R c_0 {\widehat \ot}_R c_1
{\widehat \ot}_R \dots {\widehat \ot}_R c_{n-1} {\widehat \ot}_R m^{[0]}.
$$
\end{example}

\begin{example}\label{ex:7}
Let $B$ be a left bialgebroid over a $k$-algebra $L$ and
$(A,\mu,\eta)$ be a left $B$-module algebra (so in particular an
$L$-ring). An object in $\mathcal B$ is given by the
following data.
\begin{itemize}
\item The comonads $S_l=\Hom_{-,L}(A,-)$ and $S_r=\Hom_{L,-}(A,-)$ on
  $L$-$\mathsf{Mod}$-$L$, with comonad structures
\begin{eqnarray*}
&\Hom_{-,L}(\mu,-): S_l \to S_l^2\qquad\
&\Hom_{-,L}(\eta,-):S_l \to L\textrm{-}\mathsf{Mod}\textrm{-}L
\qquad \textrm{and}\\
&\Hom_{L,-}(\mu,-): S_r \to S_r^2\qquad
&\Hom_{L,-}(\eta,-):S_r \to L\textrm{-}\mathsf{Mod}\textrm{-}L,
\end{eqnarray*}
respectively;
\item The comonad distributive law
$$
\Psi: \Hom_{-,L}(A,\Hom_{L,-}(A,-)) \cong
\Hom_{L,L}(A\ot_k A,-)\cong
\Hom_{L,-}(A,\Hom_{-,L}(A,-)),
$$
given by switching the arguments;
\item The right $\Psi$-comodule functor $(\sqcap,i)$, where
  $\sqcap=\Hom_{L,L}(L,-): L$-$\mathsf{Mod}$-$L\to \mathsf{Mod}$-$k$ and
the natural transformation
$$
i:\Hom_{L,L}(L,\Hom_{-,L}(A,-)) \cong \Hom_{L,L}(A,-)\cong \Hom_{L,L}(L,
\Hom_{L,-}(A,-))
$$
is given by the hom-tensor adjunction isomorphisms;
\item The left $\Psi$-comodule functor $(\sqcup,w)$, where
  $\sqcup:\mathsf{Ctrmod}$-$B\to L$-$\mathsf{Mod}$-$L$ is the forgetful
  functor and for any right $B$-contramodule $(Q,\alpha)$,
$$
w_Q: \Hom_{-,L}(A,Q) \to \Hom_{L,-}(A,Q),\qquad
f \mapsto \big(\ a\mapsto \alpha(f(-a))\ \big).
$$
\end{itemize}
The para-cyclic module corresponding to a right $B$-contramodule $(Q,\alpha)$
is given at degree $n$ by the $k$-module $\Hom_{L,L}(A^{\ot_L\,
  n+1},Q)$. Face and degeneracy maps are, for $k=0,\dots,n$,
\begin{eqnarray*}
&&(d_k\varphi^{(n)})(a_0\stac L a_1\stac L \dots \stac L
  a_{n-1})=\varphi^{(n)}(a_0\stac L \dots \stac L a_{k-1}\stac L 1_A \stac L
  a_k \stac L \dots \stac L a_{n-1})\\
&&(s_k\varphi^{(n)})(a_0\stac L a_1\stac L \dots \stac L
  a_{n+1})=\varphi^{(n)}(a_0\stac L \dots \stac L a_{k-1}\stac L a_k a_{k+1}
  \stac L a_{k+2} \stac L \dots \stac L a_{n+1}).
\end{eqnarray*}
The para-cyclic map is
$$
(t_n\varphi^{(n)})(a_0\ot_L a_1\ot_L \dots \ot_L a_n)=
\alpha \big(\varphi^{(n)}(a_1\ot_L \dots \ot_L a_n\ot_L (-)a_0) \big),
$$
for $\varphi^{(n)} \in \Hom_{L,L}(A^{\ot_L\,  n+1},Q)$.
\end{example}

\begin{example}\label{ex:8}
Let $B$ be a right bialgebroid over a $k$-algebra $R$ and $(A,\mu,\eta)$ be a
right $B$-comodule algebra (so in particular an $R$-ring), with $B$-coaction
$a\mapsto a^{[0]}\ot_R a^{[1]}$.
An object in ${\mathcal B}$ is given by the following data.
\begin{itemize}
\item The same comonads $S_l=\Hom_{-,R}(A,-)$ and $S_r=\Hom_{R,-}(A,-)$ on
  $R$-$\mathsf{Mod}$-$R$ constructed in Example \ref{ex:7} (replacing $L$ by
  $R$);
\item The same comonad distributive law $\Psi$ constructed in Example
  \ref{ex:7} (replacing $L$ by $R$);
\item The same right $\Psi$-comodule functor $(\sqcap,i)$ constructed in
  Example \ref{ex:7} (replacing $L$ by $R$);
\item The left $\Psi$-comodule functor $(\sqcup,w)$, where
  $\sqcup:B$-$\mathsf{Mod} \to R$-$\mathsf{Mod}$-$R$ is the forgetful functor
  and for any left $B$-module $N$,
$$
w_N: \Hom_{-,R}(A,N)\to \Hom_{R,-}(A,N),\qquad
f\mapsto \big( a\mapsto a^{[1]}f(a^{[0]})\big).
$$
\end{itemize}
The para-cyclic module corresponding to a left $B$-module $N$ has the same
simplicial structure in Example \ref{ex:7} (replacing $L$ by $R$ and $Q$ by
$N$). The para-cyclic map is
$$
(t_n\varphi^{(n)})(a_0\ot_R a_1\ot_R \dots \ot_R a_n)=
a_0^{[1]}\varphi^{(n)}(a_1\ot_R \dots \ot_R a_n\ot_R a_0^{[0]}).
$$
\end{example}

\begin{example}\label{ex:1.D}
Let $B$ be a left bialgebroid over a $k$-algebra $L$ and $(A,\mu,\eta)$ be a
left $B$-module algebra (so in particular an $L$-ring). An object in
${\mathcal B}$ is given by the following data.
\begin{itemize}
\item The comonads $S_l=A\ot_L (-)$ and $S_r=(-)\ot_L A$ on
  $A\textrm{-}\mathsf{Mod} \textrm{-} A$. For any $A$-bimodule $X$,
  $ S_l X= A\ot_L X$ is an $A$-bimodule via
$$
a_1(a\ot_L x)a_2=a_1 a\ot_L xa_2,
\qquad \textrm{for } a,a_1,a_2\in A,\ x\in X.
$$
The comonad structure is given by
$$
A\ot_L \eta\ot_L X: S_l X \to S_l^2 X
\qquad \textrm{and} \qquad
S_l X \to X,\quad a\ot_L x \mapsto ax.
$$
Symmetrically, $S_r X= X\ot_L A$ is an $A$-bimodule via
$$
a_1(x\ot_L a)a_2=a_1 x\ot_L a a_2,
\qquad \textrm{for } a,a_1,a_2\in A,\ x\in X.
$$
The comonad structure is given by
$$
X \ot_L \eta\ot_L A: S_r X \to S_r^2 X
\qquad \textrm{and} \qquad
S_r X \to X,\quad x\ot_L a \mapsto xa.
$$
\item The trivial comonad distributive law $\Psi=A\ot_L (-)\ot_L
A$; \item The right $\Psi$-comodule functor $\sqcap:A{\widehat
\otimes}_A (-):
  A\textrm{-}\mathsf{Mod} \textrm{-} A \to \mathsf{Mod} \textrm{-}k$. The
  natural transformation $i$ is given by the isomorphism
$$
i: A {\widehat \otimes}_A ((-)\ot_L A) \cong
L  {\widehat \otimes}_L (-) \cong A {\widehat \otimes}_A(A \ot_L (-));
$$
\item The left $\Psi$-comodule functor $\sqcup=(-)\ot_L
  A:\mathsf{Comod}\textrm{-}B\to A\textrm{-}\mathsf{Mod} \textrm{-} A$. For
  any right $B$-comodule $M$, with coaction $m\mapsto m_{[0]}\ot_L m_{[1]}$,
  $\sqcup M  = M\ot_L A$ is an $A$-bimodule via
$$
a_1(m\ot_L a)a_2= m_{[0]} \ot_L (m_{[1]} a_1) a a_2,
\qquad \textrm{for } a,a_1,a_2\in A,\ m\in M.
$$
The natural transformation $w: A \ot_L (-)\ot_L A \to (-)\ot_L A \ot_L A$ is
given by
$$
w_M(a\ot_L m \ot_L b)= m_{[0]}\ot_L m_{[1]} a \ot_L b.
$$
\end{itemize}
For any right $B$-comodule $M$, this determines a para-cyclic
object in $\mathsf{Mod} \textrm{-}k$. At degree $n$ it is given by
$A {\widehat \otimes}_A \big(A^{\ot_L n+1} \ot_L M \ot_L A \big)
\cong A^{\ot_L n+1} {\widehat \otimes}_L M$. For every $0\leq
k\leq n-1$, the corresponding face map is
\begin{eqnarray*}
&&d_k(a_0 {\widehat \otimes}_L \dots {\widehat \otimes}_L a_n {\widehat
    \otimes}_L m)= a_0 {\widehat \otimes}_L \dots {\widehat \otimes}_L a_{k-1}
  {\widehat \otimes}_L a_k a_{k+1} {\widehat \otimes}_L a_{k+2} {\widehat
    \otimes}_L \dots {\widehat \otimes}_L a_n {\widehat \otimes}_L m,
\end{eqnarray*}
while
\begin{eqnarray*}
&&d_n(a_0 {\widehat \otimes}_L \dots {\widehat \otimes}_L a_n
{\widehat
    \otimes}_L m)=(m_{[1]} a_n) a_0 {\widehat \otimes}_L a_1 {\widehat
  \otimes}_L \dots {\widehat \otimes}_L a_{n-1} {\widehat \otimes}_L m_{[0]}.
\end{eqnarray*}
If $0\leq k\leq n$, then the degeneracy map $s_k$ is given by
\begin{eqnarray*}
&&s_k(a_0 {\widehat \otimes}_L \dots {\widehat \otimes}_L a_n {\widehat
    \otimes}_L m)= a_0 {\widehat \otimes}_L \dots {\widehat \otimes}_L a_{k}
  {\widehat \otimes}_L 1_A {\widehat \otimes}_L a_{k+1} {\widehat \otimes}_L
  \dots {\widehat \otimes}_L a_n {\widehat \otimes}_L m.
\end{eqnarray*}
The para-cyclic map is
$$
t_n(a_0 {\widehat \otimes}_L \dots {\widehat \otimes}_L a_n {\widehat
    \otimes}_L m)= m_{[1]} a_n {\widehat \otimes}_L a_0 {\widehat \otimes}_L
\dots {\widehat \otimes}_L a_{n-1} {\widehat \otimes}_L m_{[0]}.
$$
If the left bialgebroid $B$ above is a constituent left bialgebroid in a Hopf
algebroid with a bijective antipode $S$, then there is a bijective
correspondence between right $B$-comodules $M$, with coaction $m\mapsto
m_{[0]} \ot_L m_{[1]}$, and left coactions of the constituent right
bialgebroid of the Hopf algebroid on $M$, $m\mapsto S(m_{[1]}) \ot_{L^{op}}
m_{[0]}$, cf. \ref{app:comod}. Expressing in the above formulae the right
$B$-coaction on $M$ in terms of this left coaction, we obtain a para-cyclic
module. Restricting to the case of a Hopf algebra (instead of a Hopf
algebroid) and applying to the obtained para-cyclic module the functor
$\Hom_k(-,k)$, we obtain the para-cyclic module in
\cite[(2.6)--(2.9)]{HaKhRaSo2}. (Note that in \cite{HaKhRaSo2} further
properties of the comodule $M$ are assumed which ensure that the associated
para-cyclic module has a truly cyclic subobject.)
\end{example}

\begin{example}\label{ex:2.D}
Let $B$ be a left bialgebroid over a $k$-algebra $L$ and $(A,\mu,\eta)$ be a
right $B$-comodule algebra with $B$-coaction $a\mapsto a_{[0]}\ot_L
a_{[1]}$. Then $(A,\mu,\eta)$ is in particular an $R:=L^{op}$-ring. An object
in ${\mathcal B}$ is given by the following data.
\begin{itemize}
\item The same comonads $S_l=A \ot_R(-)$ and $S_r=(-)\ot_R A$ on
  $A\textrm{-}\mathsf{Mod} \textrm{-} A$ constructed in Example \ref{ex:1.D}
  (replacing $L$ by $R$);
\item The same comonad distributive law $\Psi$ in Example
\ref{ex:1.D} (replacing $L$
  by $R$);
\item The same right $\Psi$-comodule functor $(\sqcap,i)$ in  Example
  \ref{ex:1.D} (replacing $L$ by $R$);
\item The left $\Psi$-comodule functor $\sqcup=(-)\ot_R A:B
  \textrm{-}\mathsf{Mod} \to A\textrm{-}\mathsf{Mod} \textrm{-} A$.
For any left $B$-module $N$, $\sqcup N = N \ot_R A$ is an $A$-bimodule via
$$
a_1(m\ot_R a)a_2=a_{1[1]}m\ot_R a_{1[0]} aa_2,
\qquad \textrm{for }a,a_1,a_2\in A,\ m\in N.
$$
The natural transformation $w:A \ot_R (-) \ot_R A \to (-)\ot_R A \ot_R A$ is
given by
$$
w_N(a\ot_R m \ot_R b) = a_{[1]} m \ot_R a_{[0]} \ot_R b.
$$
\end{itemize}
For any left $B$-module $N$, this determines a para-cyclic object
in $\mathsf{Mod} \textrm{-}k$. At degree $n$ it is given by $A
{\widehat \otimes}_A \big(A^{\ot_R n+1} \ot_R N \ot_R A \big)
\cong A^{\ot_R n+1} {\widehat \otimes}_R N$. For every $0\leq
k\leq n-1$, the corresponding face map is
\begin{eqnarray*}
&&d_k(a_0 {\widehat \otimes}_R \dots {\widehat \otimes}_R a_n {\widehat
    \otimes}_R m)= a_0 {\widehat \otimes}_R \dots {\widehat \otimes}_R a_{k-1}
  {\widehat \otimes}_R a_k a_{k+1} {\widehat \otimes}_R a_{k+2} {\widehat
    \otimes}_R \dots {\widehat \otimes}_R a_n {\widehat \otimes}_R
    m,
\end{eqnarray*}
while
\begin{eqnarray*}
&&d_n(a_0 {\widehat \otimes}_R \dots {\widehat \otimes}_R a_n
{\widehat
    \otimes}_R m)=a_{n[0]} a_0 {\widehat \otimes}_R a_1 {\widehat
  \otimes}_R \dots {\widehat \otimes}_R a_{n-1} {\widehat \otimes}_R a_{n[1]}m.
\end{eqnarray*}
If $0\leq k\leq n$, then the degeneracy map $s_k$ is given by
\begin{eqnarray*}
&&s_k(a_0 {\widehat \otimes}_R \dots {\widehat \otimes}_R a_n {\widehat
    \otimes}_R m)= a_0 {\widehat \otimes}_R \dots {\widehat \otimes}_R a_{k}
  {\widehat \otimes}_R 1_A {\widehat \otimes}_R a_{k+1} {\widehat \otimes}_R
  \dots {\widehat \otimes}_R a_n {\widehat \otimes}_R m.
\end{eqnarray*}
The para-cyclic map is
$$
t_n(a_0 {\widehat \otimes}_R \dots {\widehat \otimes}_R a_n {\widehat
    \otimes}_R m)=a_{n[0]}{\widehat \otimes}_R a_0 {\widehat \otimes}_R
\dots {\widehat \otimes}_R a_{n-1} {\widehat \otimes}_R a_{n[1]} m.
$$
This is a non-commutative base version of the para-cyclic module in
\cite[(3.5)--(3.8)]{HaKhRaSo2}. Note that in \cite{HaKhRaSo2} further
properties of the left module $N$ are assumed which ensure that the associated
para-cyclic module has a truly cyclic subobject.
\end{example}

A {\em bicontramodule} of an $R$-coring $C$ is an $R$-bimodule $Y$, together
with a right $C$-contramodule structure $\beta_r:\Hom_{-,R}(C,Y)\to Y$ and a
left $C$-contramodule structure $\beta_l:\Hom_{R,-}(C,Y)\to Y$, such that
$\beta_r$ is a left $R$-module map, $\beta_l$ is a right $R$-module map
and
$$
\beta_l \circ \Hom_{R,-}(C,\beta_r)=\beta_r \circ \Hom_{-,R}(C,\beta_l),
$$
up to the (suppressed) canonical isomorphism
$$
\Hom_{R,-}(C,\Hom_{-,R}(C,Y))\cong \Hom_{R,R}(C\ot_k C,Y) \cong
\Hom_{-,R}(C,\Hom_{R,-}(C,Y)).
$$
A morphism of bicontramodules is a right contramodule map as well as a left
contramodule map. The category of $C$-bicontramodules is denoted by
$C\textrm{-}\mathsf{Ctrmod} \textrm{-}C$.

\begin{example}\label{ex:3.D}
Let $B$ be a right bialgebroid over a $k$-algebra $R$ and
$(C,\Delta,\epsilon)$ be a right
$B$-module coring (hence in particular an $R$-coring). An object in ${\mathcal
B}$ is given by the following data.
\begin{itemize}
\item The comonads $S_l= \Hom_{-,R}(C,-)$ and $S_r=
\Hom_{R,-}(C,-)$ on $C\textrm{-}\mathsf{Ctrmod} \textrm{-}C$. For
any $C$-bicontramodule $(Y,\beta_l,\beta_r)$, $S_l Y =
\Hom_{-,R}(C,Y)$ is a bicontramodule, via the structure maps
\begin{eqnarray*}
&&\Hom_{R,-}(C,\Hom_{-,R}(C,Y))\cong \Hom_{-,R}(C,\Hom_{R,-}(C,Y))
\stackrel{\Hom_{-,R}(C, \beta_l)}{\longrightarrow} \Hom_{-,R}(C,Y)\\
&&\Hom_{-,R}(C,\Hom_{-,R}(C,Y))\cong \Hom_{-,R}(C\ot_R C,Y))
\stackrel{\Hom_{-,R}(\Delta,Y)}{\longrightarrow} \Hom_{-,R}(C,Y).
\end{eqnarray*}
The comonad structure is given by
$$
\Hom_{-,R}(C \ot_R \epsilon,Y): S_l Y\to S_l^2 Y
\qquad \textrm{and}\qquad
\beta_r: S_l Y \to Y.
$$
Symmetrically, $S_r Y = \Hom_{R,-}(C,Y)$ is a bicontramodule, via
the structure maps
\begin{eqnarray*}
&&\Hom_{R,-}(C,\Hom_{R,-}(C,Y))\cong \Hom_{R,-}(C\ot_R C,Y))
\stackrel{\Hom_{R,-}(\Delta,Y)}{\longrightarrow} \Hom_{R,-}(C,Y)\\
&&\Hom_{-,R}(C,\Hom_{R,-}(C,Y))\cong \Hom_{R,-}(C,\Hom_{-,R}(C,Y))
\stackrel{\Hom_{R,-}(C, \beta_r)}{\longrightarrow} \Hom_{R,-}(C,Y).
\end{eqnarray*}
The comonad structure is given by
$$
\Hom_{R,-}(\epsilon \ot_R C,Y): S_r Y\to S_r^2 Y
\qquad \textrm{and}\qquad
\beta_l: S_r Y \to Y.
$$
\item The comonad distributive law
$$
\Psi:\Hom_{-,R}(C,\Hom_{R,-}(C,-))\cong \Hom_{R,R}(C\ot_k C,-)\cong
\Hom_{R,-}(C,\Hom_{-,R}(C,-)),
$$
given by switching the arguments;
\item The right $\Psi$-comodule functor $\sqcap:C\textrm{-}\mathsf{Ctrmod}
  \textrm{-}C \to \mathsf{Mod}\textrm{-}k$, given by the coequalizer
$$
\xymatrix{
\Hom_{R,R}(C,Y)\ar@<2pt>[rrr]^-{\Hom_{R,R}(R,\alpha_l)}
\ar@<-2pt>[rrr]_-{\Hom_{R,R}(R,\alpha_r)} &&&
\Hom_{R,R}(R,Y) \ar[r]& \sqcap Y.
}
$$
The natural transformation $i$ is given by the isomorphism
$$
i: \sqcap \Hom_{R,-}(C,-)\cong \Hom_{R,R}(R,-) \cong \sqcap \Hom_{-,R}(C,-).
$$
\item The left $\Psi$-comodule functor
  $\sqcup=\Hom_{R,-}(C,-):\mathsf{Ctrmod}\textrm{-} B \to
  C\textrm{-}\mathsf{Ctrmod} \textrm{-}C$. For any right $B$-contramodule
  $(Q,\alpha)$, $\sqcup Q =\Hom_{R,-}(C,Q)$ is a $C$-bicontramodule, via the
structure maps
\begin{align*} & \Hom_{R,-}(C,\Hom_{R,-}(C,Q))
\stackrel{\cong}{\longrightarrow} \Hom_{R,-}(C\ot_R C,Q)
\stackrel{\Hom_{R,-}(\Delta,Q)}{\longrightarrow}
\Hom_{R,-}(C,Q);\\
& \Hom_{-,R}(C,\Hom_{R,-}(C,Q))\stackrel{w_Q}{\longrightarrow}
\Hom_{R,-}(C,\Hom_{R,-}(C,Q))\stackrel{\cong}{\longrightarrow}
\Hom_{R,-}(C\ot_R C,Q)) \\
 & \hspace*{9.75cm}\stackrel{\Hom_{R,-}(\Delta,Q)}{\longrightarrow}
 \Hom_{R,-}(C,Q),
\end{align*}
where $w:\Hom_{-,R}(C,\Hom_{R,-}(C,-))\to
\Hom_{R,-}(C,\Hom_{R,-}(C,-))$ is given by
$$
\big(w_Q(h)\big)(c)(d)=\alpha\big(h(d-)(c)\big).
$$
\end{itemize}
For any right $B$-contramodule $(Q,\alpha)$, this determines a para-cyclic
object in $\mathsf{Mod}\textrm{-}k$. At degree $n$, it is given by
$\sqcap \Hom_{-,R}(C^{\ot_R n+1},\Hom_{R,-}(C,Q))
\cong \Hom_{R,R}(C^{\ot_R
    n+1},Q)$. Denote $\Delta(c)=c\1 \ot_R c\2$. For every $0\leq k\leq n-1$,
the corresponding face map is
\begin{eqnarray*}
&&(d_k \varphi^{(n)})(c_0\stac R \dots \stac R c_{n-1})=\varphi^{(n)}(c_0\stac R
  \dots \stac R c_{k-1}\stac R \Delta(c_k)\stac R c_{k+1}\stac R \dots \stac R
  c_{n-1}),
\end{eqnarray*}
while
\begin{eqnarray*}
&&(d_n \varphi^{(n)})(c_0\stac R \dots \stac R c_{n-1})=
\alpha\big(\varphi^{(n)}({c_0}\2 \stac R c_1 \stac R \dots \stac R
c_{n-1} \stac R {c_0}\1 -)\big).
\end{eqnarray*}
If $0\leq k\leq n$, then the degeneracy map $s_k$ is given by
\begin{eqnarray*}
&&(s_k \varphi^{(n)})(c_0\stac R \dots \stac R c_{n+1})=\varphi^{(n)}(c_0\stac R
  \dots \stac R c_{k-1}\stac R c_k\epsilon(c_{k+1})\stac R c_{k+2}\stac R
  \dots \stac R c_{n+1}).
\end{eqnarray*}
The para-cyclic map is
$$
(t_n \varphi^{(n)})(c_0\ot_R \dots \ot_R c_{n})=
\alpha\big(\varphi^{(n)}(c_1\ot_R \dots \ot_R c_{n} \ot_R c_0-)\big),
$$
for $\varphi^{(n)}\in \Hom_{R,R}(C^{\ot_R n+1},Q)$.

This is a non-commutative base version of the para-cyclic module in
\cite[p 4]{Brz:coef} (though note the minor difference of using a left or a
right module coalgebra $C$). In \cite{Brz:coef} additional properties of the
contramodule $Q$ are assumed so that the associated para-cyclic module has a
cyclic subobject.
\end{example}

\begin{example}\label{ex:4.D}
Let $B$ be a right bialgebroid over a $k$-algebra $R$ and
$(C,\Delta,\epsilon)$ be a left $B$-comodule coring with coaction $c\mapsto
c^{[-1]}\ot_R c^{[0]}$. Then
$(C,\Delta,\epsilon)$ is in particular an $L:=R^{op}$-coring. An object in
${\mathcal B}$ is given by the following data.
\begin{itemize}
\item The same comonads $S_l=\Hom_{-,L}(C,-)$ and $S_r=\Hom_{L,-}(C,-)$ on
$C\textrm{-}\mathsf{Ctrmod} \textrm{-}C$ constructed in Example
  \ref{ex:3.D} (replacing $R$ by $L$);
\item The same distributive law $\Psi$ in Example \ref{ex:3.D} (replacing $R$
  by $L$);
\item The same right $\Psi$-comodule functor $(\sqcap,i)$ in Example
  \ref{ex:3.D} (replacing $R$ by $L$);
\item The left $\Psi$-comodule functor $\sqcup=\Hom_{L,-}(C,-):B\textrm{-}
  \mathsf{Mod} \to C\textrm{-}\mathsf{Ctrmod} \textrm{-}C$. For any left
  $B$-module $N$, $\sqcup N= \Hom_{L,-}(C,N)$ is a bicontramodule via the
structure maps
\begin{align*}
&\Hom_{L,-}(C,\Hom_{L,-}(C,N))\stackrel{\cong}{\longrightarrow}
\Hom_{L,-}(C\ot_L C,N)
\stackrel{\Hom_{L,-}(\Delta,N)}{\longrightarrow}
\Hom_{L,-}(C,N)\\
&\Hom_{-,L}(C,\Hom_{L,-}(C,N))\stackrel{w_N}{\longrightarrow}
\Hom_{L,-}(C,\Hom_{L,-}(C,N))\stackrel{\cong}{\longrightarrow}
\Hom_{L,-}(C\ot_L C,N)\\
&\hspace*{9.8cm} \stackrel{\Hom_{L,-}(\Delta,N)}{\longrightarrow}
\Hom_{L,-}(C,N),
\end{align*}
where $w:\Hom_{-,L}(C,\Hom_{L,-}(C,-)) \to
\Hom_{L,-}(C,\Hom_{L,-}(C,-))$ is given by
$$
\big(w_N(h)\big)(c)(d)=d^{[-1]}h(d^{[0]})(c).
$$
\end{itemize}
For any left $B$-module $N$, this determines a para-cyclic object
in $\mathsf{Mod}\textrm{-}k$. At degree $n$, it is given by
$\sqcap \Hom_{-,L}(C^{\ot_L n+1}, \Hom_{L,-}(C,N)) \cong
\Hom_{L,L}(C^{\ot_L n+1},N)$. Denote $\Delta(c)=c\1 \ot_L c\2$.
For every $0\leq k\leq n-1$, the corresponding face map is
\begin{eqnarray*}
&&(d_k \varphi^{(n)})(c_0\stac L \dots \stac L c_{n-1})=\varphi^{(n)}(c_0\stac L
  \dots \stac L c_{k-1}\stac L \Delta(c_k)\stac L c_{k+1}\stac L \dots \stac L
  c_{n-1}),
\end{eqnarray*}
while
\begin{eqnarray*}
&&(d_n \varphi^{(n)})(c_0\stac L \dots \stac L c_{n-1})=
{{c_0}\1}^{[-1]} \varphi^{(n)}({c_0}\2 \stac L c_1 \stac L \dots
\stac L c_{n-1} \stac L
         {{c_0}\1}^{[0]} ).
\end{eqnarray*}
If $0\leq k\leq n$, then the degeneracy map $s_k$ is given by
\begin{eqnarray*}
&&(s_k \varphi^{(n)})(c_0\stac L \dots \stac L c_{n+1})=\varphi^{(n)}(c_0\stac L
  \dots \stac L c_{k-1}\stac L c_k\epsilon(c_{k+1})\stac L c_{k+2}\stac L
  \dots \stac L c_{n+1}).
\end{eqnarray*}
The para-cyclic map is
$$
(t_n \varphi^{(n)})(c_0\ot_L \dots \ot_L c_{n})=
c_0^{[-1]}\varphi^{(n)}(c_1\ot_L \dots \ot_L c_{n} \ot_L c_0^{[0]}),
$$
for $\varphi^{(n)}\in \Hom_{L,L}(C^{\ot_L n+1},N)$.
\end{example}


\section{The cyclic duality functor}

\label{sec:duality}

The functor recently known as the \emph{cyclic duality} functor,
appeared first in Connes' work \cite{Co:CycDual}.
In its original form, it is an isomorphism between the category of cyclic
objects and the category of cocyclic objects in a given category. It was
extended in \cite{KR} to an isomorphism between certain full subcategories of
the categories of para-cyclic, and of para-cocyclic objects. The objects of
these full subcategories are those para-(co)cyclic objects whose
para-(co)cyclic morphisms are isomorphisms at all degrees. The aim of the
current section is to extend cyclic duality to a functor between appropriate
subcategories of ${\mathcal{A}}$ and ${\mathcal{B}}$ in Definitions
\ref{def:cat_A} and \ref{def:cat_B}, respectively.

Connes's cyclic duality functor
(in the extended form in \cite{KR})
and also its dual version (from a subcategory of the category of para-cyclic
objects to a subcategory of the category of para-cocyclic objects) both will
be denoted by $\widehat{(-)}$.

Denote by ${\mathcal{A}}^\times$ the full subcategory of ${\mathcal{A}}$ in
Definition \ref{def:cat_A}, whose objects $(T_l,T_r,\Phi,\sqcap,i,$ $
\sqcup,w)$ obey the property that $\Phi$, $i$ and $w$ are natural
isomorphisms. In the category ${\overline {\mathcal{P}}}$ in Definition \ref
{def:cat_cocyc}, introduce the full subcategory ${\overline {\mathcal{P}}}
^\times$ whose objects have para-cocyclic morphisms which are
natural isomorphisms at all degrees. Clearly, the functor
${\mathcal{Z}}^*$ in Theorem \ref{thm:Z^*} induces a functor
${\mathcal{Z}}^{*\times}:{\mathcal{A}}^\times \to { \overline
  {\mathcal{P}}}^\times$.
Symmetrically, introduce the full subcategory ${\mathcal{B}}^\times$ of the
category ${\mathcal{B}}$ in Definition \ref{def:cat_B}, for whose objects $
(S_l,S_r,\Psi,\sqcap,i,\sqcup,w)$ the natural transformations $\Psi$, $i$
and $w$ are isomorphisms. By Theorem \ref{thm:Z_*}, there is an induced
functor ${\mathcal{Z}}_*^\times: {\mathcal{B}}^\times \to {\underline {
\mathcal{P}}}^\times$, where ${\underline {\mathcal{P}}}^\times$ is the full
subcategory of ${\underline {\mathcal{P}}}$, for whose objects the
para-cyclic morphisms are natural isomorphisms at all degrees.
Finally, denote by ${\mathcal{A}}^\times_c$ the full subcategory of ${
\mathcal{A}}^\times$, for whose objects $(T_l,T_r,\Phi,\sqcap,i,\sqcup,w)$
the codomain category ${\mathcal{C}}$ of the right $\Phi$-module functor $
\sqcap$ possesses coequalizers. Symmetrically, denote by ${\mathcal{B}}
^\times_e$ the full subcategory of ${\mathcal{B}}^\times$, for whose objects
$(S_l,S_r,\Psi,\sqcap,i,\sqcup,w)$ the codomain category ${\mathcal{C}}$ of
the right $\Psi$-comodule functor $\sqcap$ possesses equalizers.
Restrictions of the functor ${\mathcal{Z}}^{*\times}$ to ${\mathcal{A}}
^\times_c$ and ${\mathcal{Z}}_*^{\times}$ to ${\mathcal{B}}^\times_e$ are
denoted by the same symbols ${\mathcal{Z}}^{*\times}$ and ${\mathcal{Z}}
_*^{\times}$.

\begin{theorem}
\label{thm:cyclic_dual} Using the notation in the paragraph preceding the
theorem, there exists a functor $(\widetriangle{-}):{\mathcal{A}}^\times_c
\to {\mathcal{B}}^\times$, such that the diagram
\begin{equation*}
\xymatrix{
{\mathcal A}^\times_c \ar[rr]^-{(\wt{\,-\,})}
\ar[d]_-{{\mathcal Z}^{*\times}}&&
{\mathcal B}^\times \ar[d]^-{{\mathcal Z}_*^\times}\\
{\overline {\mathcal P}}^\times \ar[rr]_-{(\widehat{\,-\,})}&&
{\underline {\mathcal P}}^\times }
\end{equation*}
commutes up to a natural isomorphism.

Symmetrically, there exists a functor $(\widetriangle{-}):{\mathcal{B}}
^\times_e \to {\mathcal{A}}^\times$, such that the functors ${\mathcal{Z}}
^{*\times} (\widetriangle{-})$ and $\widehat{{\mathcal{Z}}_*^{\times}(-)}
$ are naturally isomorphic.
\end{theorem}

The proof of Theorem \ref{thm:cyclic_dual} goes through a series of lemmata.

Recall that for the Eilenberg-Moore category $\mathcal{M}^T$ of a monad $T:
\mathcal{M }\to \mathcal{M}$, there is a \emph{forgetful functor} $U:
\mathcal{M}^T \to \mathcal{M}$, with object map $(M,\varrho) \mapsto M$ and
acting on the morphisms as the identity map. The forgetful functor $U$ has a
left adjoint $F$, with object map $M \mapsto (TM,mM)$ and morphism map $f
\mapsto Tf$.

\begin{lemma}
\label{lem:Phi.hat} Let $(T_l,m_l,u_l)$ and $(T_r,m_r,u_r)$ be monads on the
same category $\mathcal{M}$ and $\Phi:T_r T_l\to T_l T_r$ be a distributive
law which is a natural isomorphism. Consider the induced monad
\eqref{eq:composite_monad} and the forgetful functor $U:\mathcal{M}^{T_l
T_r}\to \mathcal{M}$.

(1) There is a comonad $({\widetriangle T}_l, d_l,e_l)$ on $\mathcal{M}^{T_l
T_r}$ such that $U {\widetriangle T}_l = T_l U$, $U d_l = T_l u_l U$ and $U
e_l (M,\varrho)= \varrho \circ T_l u_r M$, for any $T_l T_r$-algebra $
(M,\varrho)$.

(2) There is a comonad $({\widetriangle T}_r,d_r,e_r)$ on $\mathcal{M}^{T_l
T_r}$ such that $U {\widetriangle T}_r = T_r U$, $U d_r = T_r u_r U$ and $U
e_r (M,\varrho)= \varrho \circ u_l T_r M$, for any $T_l T_r$-algebra $
(M,\varrho)$.

(3) There is a comonad distributive law ${\widetriangle \Phi}:{
\widetriangle
T}_l {\widetriangle T}_r \to {\widetriangle T}_r {\widetriangle T}_l$, such
that $U {\widetriangle \Phi} = \Phi^{-1} U$.
\end{lemma}

\begin{proof}
(1) By Beck's classical theorem \cite[p 122]{Be}, the distributive law $\Phi$
induces a monad $({\widetilde T}_l,{\tilde m}_l,{\tilde u}_l)$ on the
category of $T_r$-algebras, such that the forgetful functor $U_r:\mathcal{M}
^{T_r}\to \mathcal{M}$ satisfies
\begin{equation*}
U_r {\widetilde T}_l= T_l U_r,\qquad U_r {\tilde{m}}_l= m_l U_r,\qquad U_r {
\tilde{u}}_l= u_l U_r.
\end{equation*}
Moreover, the category $(\mathcal{M}^{T_r})^{{\tilde T}_l} $ of ${\widetilde
T}_l$-algebras is isomorphic to the category of $T_l T_r$-algebras. Consider
the forgetful functor ${\widetilde U}_l:(\mathcal{M}^{T_r})^{{\tilde T}_l}
\to \mathcal{M}^{T_r}$ and its left adjoint ${\widetilde F}_l$. The
composite functor ${\widetilde U}_l{\widetilde F}_l$ is equal to ${
\widetilde T}_l$ while ${\widetriangle T}_l:={\widetilde F}_l{\widetilde U}
_l $ is a comonad on $(\mathcal{M}^{T_r})^{{\tilde T}_l}\cong \mathcal{M}
^{T_l T_r}$. Its comultiplication is given by ${\widetilde F}_l {\tilde u}_l
{\widetilde U}_l:{\widetilde F}_l{\widetilde U}_l \to {\widetilde F}_l{
\widetilde T}_l {\widetilde U}_l = {\widetilde F}_l{\widetilde U}_l{
\widetilde F}_l{\widetilde U}_l$ and the counit is ${\tilde \varrho}_l:{
\widetilde F}_l {\widetilde U}_l(M,{\tilde \varrho}_l)=({\widetilde T}_l M,{
\tilde m}_l M) \to (M,{\tilde \varrho}_l)$, for any ${\widetilde T}_l$
-algebra $(M,{\tilde \varrho}_l)$. Since the composite functor $U_r {
\widetilde U}_l$ differs from the forgetful functor $U:\mathcal{M}^{T_l
T_r}\to \mathcal{M}$ by the isomorphism $(\mathcal{M}^{T_r})^{{\tilde T}
_l}\cong \mathcal{M}^{T_l T_r}$, the comonad ${\widetriangle T}_l$ obeys the
required properties.

Part (2) follows by applying the same reasoning as in part (1) to the
distributive law $\Phi^{-1}$.

For any $T_{l}T_{r}$-algebra $(M,\varrho )$, the $T_{l}T_{r}$-actions on ${
\widetriangle T}_{l}(M,\varrho )$ and ${\widetriangle T}_{r}(M,\varrho )$
are given by the respective morphisms $\wt{\rho}_l:T_{l}T_{r}T_{l}M\rightarrow
T_{l}M$ and $\wt{\rho}_r:T_{l}T_{r}T_{r}M\rightarrow T_{r}M$, where
\begin{align}
&\wt{\rho}_l:=T_{l}\varrho \circ T_{l}u_{l}T_{r}M\circ m_{l}T_{r}M\circ
  T_{l}\Phi M\
\label{eq:ac_on_T^l} \\
&\wt{\rho}_r:=T_{r}\varrho \circ T_{r}T_{l}u_{r}M\circ \Phi ^{-1}M\circ
  T_{l}m_{r}M.
\label{eq:ac_on_T^r}
\end{align}

(3) The composite monad $T_lT_r$ in \eqref{eq:composite_monad}
induces a monad $\mathsf{Cat}(\M,T_lT_r)$ on the category
$\mathsf{Cat}(\M,\M)$ of functors $\M\to \M$. The natural transformation
$\Phi^{-1}$ yields a morphism of its algebras $(T_l T_r,T_l m_r \circ m_l T_r
T_r\circ T_l \Phi T_r) \to (T_r T_l, T_r m_l \circ \Phi^{-1} T_l\circ T_l m_r
T_l)$. Indeed, using string computation, we have
\[
\includegraphics{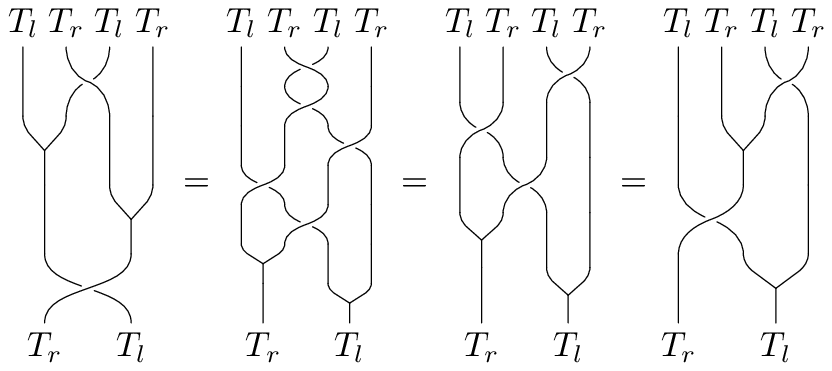}
\]
where for the first and the third equations one uses the definition of
distributive laws and the second relation follows by
$\Phi^{-1}\circ\Phi=T_rT_l$.
Comparing these actions with \eqref{eq:ac_on_T^l} and \eqref{eq:ac_on_T^r},
this implies the existence of a natural transformation ${\widetriangle\Phi}:
{\widetriangle T}_l{\widetriangle T}_r\to {\widetriangle T}_r {\widetriangle
  T}_l$ such that $U {\widetriangle \Phi}= \Phi^{-1}U$. Since the forgetful
functor $U$ reflects isomorphisms, $\wt \Phi$ is a natural
isomorphism.
Using that $U$ is faithful, it is easy to check that since $\Phi$ is a monad
distributive law, ${\widetriangle \Phi}$ is a comonad distributive law.
\end{proof}

For an object $(T_{l},T_{r},\Phi ,\sqcap :\mathcal{M}\rightarrow {\mathcal{C}
},i,\sqcup :{\mathcal{D}}\rightarrow \mathcal{M},w)$ of ${\mathcal{A}}
_{c}^{\times }$, consider the forgetful functor $U:\mathcal{M}
^{T_{l}T_{r}}\rightarrow \mathcal{M}$ and the natural transformations $\xi
_{l}:T_{l}U\rightarrow U$ and $\xi _{r}:T_{r}U\rightarrow U$, given for any $
T_{l}T_{r}$-algebra $(M,\varrho )$ by the morphisms
\begin{equation} \label{eq:xi_l_r}
\xi _{l}(M,\varrho ):=\varrho \circ T_{l}u_{r}M\qquad \text{and}\qquad \xi
_{r}(M,\varrho ):=\varrho \circ u_{l}T_{r}M.
\end{equation}
Since coequalizers in ${\mathcal{C}}$ exist by assumption, we can define a
functor ${\widetriangle\sqcap }:\mathcal{M}^{T_{l}T_{r}}\rightarrow {
\mathcal{C}}$ via the coequalizer
\begin{equation}
\xymatrix{ \sqcap T_l U \ar@<2pt>[rr]^-{\sqcap \xi_l}
\ar@<-2pt>[rr]_-{\sqcap \xi_r \circ i U}&& \sqcap U \ar[rr]^-{p}&&
{\widetriangle \sqcap} }  \label{eq:Pi.hat}
\end{equation}
in the category of functors. For any $T_{l}T_{r}$-algebra $(M,\varrho )$, $
p(M,\varrho ):\sqcap M\rightarrow \widetriangle{\sqcap}(M,\varrho )$ is the
coequalizer of $\sqcap \varrho \circ \sqcap T_{l}u_{r}M$ and $\sqcap \varrho
\circ \sqcap u_{l}T_{r}M\circ iM$. For any morphism $f:(M,\varrho
)\rightarrow (M^{\prime },\varrho ^{\prime })$ in $\mathcal{M}^{T_{l}T_{r}}$,
the composite $p(M^{\prime },\varrho ^{\prime })\circ \sqcap f$
coequalizes the parallel morphisms in \eqref{eq:Pi.hat} (evaluated at $
(M,\varrho )$). Hence we can define $\widetriangle{\sqcap}f$ as the unique
morphism for which $\widetriangle{\sqcap}f\circ p(M,\varrho )=p(M^{\prime
},\varrho ^{\prime })\circ \sqcap f$.

\begin{lemma}
\label{lem:i.hat} Consider an object $(T_{l},T_{r},\Phi ,\sqcap,i,\sqcup,w)$
of ${\mathcal{A}}_{c}^{\times }$. For the forgetful functor $U:\mathcal{M}
^{T_{l}T_{r}}\rightarrow \mathcal{M}$, the monads ${\widetriangle T}_{l}$
and ${\widetriangle T}_{r}$ in Lemma \ref{lem:Phi.hat} and the functor ${
\widetriangle\sqcap }$ in \eqref{eq:Pi.hat},
there are natural isomorphisms $\theta _{l}:{\widetriangle\sqcap
}{\widetriangle T}_{l}\rightarrow \sqcap U$ and $\theta
_{r}:{\widetriangle\sqcap }{\widetriangle T}_{r}\rightarrow \sqcap U$ such
that
\begin{equation*}
\theta _{l}\circ p{\widetriangle T}_{l}\circ i^{-1}U=\sqcap \xi_{r}
\qquad \textrm{and}\qquad
\theta _{r}\circ p{\widetriangle T}_{r}\circ iU=\sqcap \xi _{l}.
\end{equation*}
\end{lemma}

\begin{proof}
By definition, ${\widetriangle\sqcap }{\widetriangle T}_{r}$ is the
coequalizer of the natural transformations $\sqcap T_{r}\xi _{l}\circ \sqcap
\Phi ^{-1}U$ and $\sqcap m_{r}U\circ iT_{r}U$. Since $\Phi $ is an
isomorphism, this is equivalent to the coequalizer of $\sqcap T_{r}\xi _{l}$
and $\sqcap m_{r}U\circ iT_{r}U\circ \sqcap \Phi U=iU\circ \sqcap
m_{l}U\circ i^{-1}T_{l}U$, where we used that $(\sqcap ,i)$ is a right $\Phi
$-module functor, cf. first condition in \eqref{eq:right_t-module}. Using
that $i$ is an isomorphism, we conclude that ${\widetriangle\sqcap }{
\widetriangle T}_{r}$ is the coequalizer of $i^{-1}U\circ \sqcap T_{r}\xi
_{l}\circ iT_{l}U=\sqcap T_{l}\xi _{l}$ and $\sqcap m_{l}U$. The coequalizer
\begin{equation*}
\xymatrix{ T_lT_lU \ar@<2pt>[rr]^-{T_l \xi_l}\ar@<-2pt>[rr]_-{ m_l U}&& T_l
U \ar[rr]^-{\xi_l}&& U }
\end{equation*}
is split, as $(M,\xi _{l}(M,\varrho ))$ is a $T_{l}$-algebra for any
$T_{l}T_{r}$-module $(M,\varrho )$. Hence it is preserved by composing with
$\sqcap $ on the left. Thus, by the universal property of coequalizers,
there is a unique natural isomorphism $\theta _{r}:{
\widetriangle\sqcap }{\widetriangle T}_{r}\rightarrow \sqcap U$ such that
$\theta _{r} \circ p{\widetriangle T}_{r}\circ iU=\sqcap \xi _{l}$. The
existence of the isomorphism  $\theta _{l}:{\widetriangle\sqcap }{
\widetriangle T}_{l}\rightarrow \sqcap U$ is proven similarly.
\end{proof}

\begin{lemma}
\label{lem:U.hat} Consider an object $(T_{l},T_{r},\Phi ,\sqcap,i,\sqcup,w)$
of ${\mathcal{A}}_{c}^{\times }$. Then $T_{r}\sqcup $ is a left module
functor for the composite monad $T_{l}T_{r}$ in \eqref{eq:composite_monad}.
Hence there is a functor ${\widetriangle\sqcup }:{\mathcal{D}}
\rightarrow {\mathcal{M}}^{T_{l}T_{r}}$, with object map and morphism map
\begin{equation*}
\wt{\sqcup}Y:=(T_{r}\sqcup Y,w^{-1}Y\circ m_{l}\sqcup Y\circ T_{l}wY\circ
T_{l}m_{r}\sqcup Y)
\quad \text{and}\quad
\wt{\sqcup}f:= T_{r}\sqcup f.
\end{equation*}
\end{lemma}

\begin{proof}
Unitality of the given $T_{l}T_{r}$-action is immediate by unitality of the
multiplications of $T_{l}$ and $T_{r}$. Associativity is checked as follows.
\begin{equation*}
\includegraphics{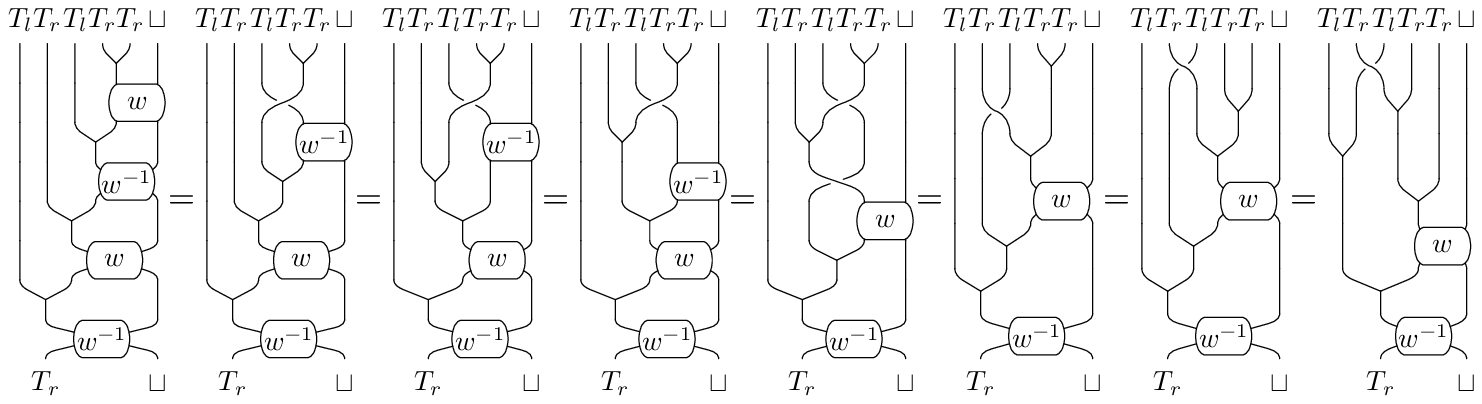}
\end{equation*}
The first and the fourth equalities follow by using that $(\sqcup ,w)$ is a
left $\Phi $-module functor, cf. first identity in \eqref{eq:left_t-module}.
The second equality follows by associativity of $m_{r}$. Naturality is used
in the third equality, and also in the last two ones, in the case of the last
equality together with the associativity of $m_{l}$ and of $m_{r}$.
The fifth equality is a consequence of the fact that $\Phi $ is a distributive
law.
\end{proof}

\begin{proof}[Proof of Theorem \protect\ref{thm:cyclic_dual}]
Let $(T_l,T_r,\Phi,\sqcap,i,\sqcup,w)$ be an object in ${\mathcal{A}}
^\times_c$, with $\sqcap:\mathcal{M}\to \mathcal{C}$ and $\sqcup:{\mathcal{D}
}\to \mathcal{M}$. First we show that the septuple $({\widetriangle T}_l,{
\widetriangle T}_r,{\widetriangle \Phi}, {\widetriangle\sqcap}, {
\widetriangle i}, {\widetriangle \sqcup}, {\widetriangle w})$ is an object
of ${\mathcal{B}}$, where $U:\mathcal{M}^{T_l T_r} \to {\mathcal{M}}$ is the
forgetful functor, the comonads ${\widetriangle T}_l$ and ${\widetriangle T}
_r$ and the comonad distributive law ${\widetriangle \Phi}$ are constructed
in Lemma \ref{lem:Phi.hat}, the functor ${\widetriangle
\sqcap}$ is defined by \eqref{eq:Pi.hat}, the functor ${\widetriangle
\sqcup}$ is constructed in Lemma \ref{lem:U.hat}, while $\widetriangle i$
and $\widetriangle w$ are given by the relations
\begin{equation}  \label{eq:ihat-what}
{\widetriangle i}:={\theta_l}^{-1}\circ\theta_r
\qquad\textrm{and}\qquad
U {\widetriangle w}:= T_r w^{-1} \circ \Phi^{-1} \sqcup.
\end{equation}
The following string computation proves that ${\widetriangle w}Y:{\wt
  T}_l {\wt \sqcup}Y\to {\wt T}_r {\wt \sqcup}Y$ is a
morphism in
$\M^{T_lT_r}$, for any object $Y$ in ${\mathcal D}$.
\begin{equation*}
\includegraphics{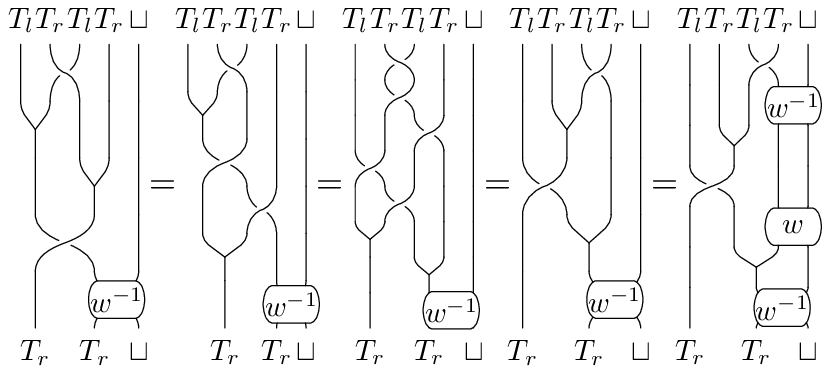}
\end{equation*}
Here we used the form of the $T_lT_r$-actions on ${\widetriangle
T}_l(M,\varrho)$, ${\widetriangle T}_r(M,\varrho)$ and ${\widetriangle
\sqcup}Y$, for any $T_lT_r$-algebra $(M,\varrho)$ and for any object $Y$ of
${\mathcal D}$; see relations \eqref{eq:ac_on_T^l}, \eqref{eq:ac_on_T^r}
and Lemma \ref{lem:U.hat}. For the first, second and third equalities
we also used that $\Phi$ is a distributive law and that the multiplication of
a monad is a natural transformation. Hence ${\widetriangle w}$ can be regarded
as a natural transformation ${\widetriangle T}_l {\widetriangle \sqcup} \to
{\widetriangle T}_r {\widetriangle \sqcup}$.

By naturality of $p$ in \eqref{eq:Pi.hat} and by Lemma
\ref{lem:Phi.hat},
$$
p {\wt T}_r^2\circ \sqcap T_r u_r U = {\wt \sqcap}d_r\circ p {\wt T}_r,\qquad
p {\wt T}_l^2\circ \sqcap T_l u_l U = {\wt \sqcap}d_l\circ p {\wt T}_l,\qquad
p {\wt T}_r {\wt T}_l \circ  \sqcap \Phi^{-1}U ={\wt \sqcap}{\wt \Phi} \circ p
{\wt T}_l {\wt T}_r .
$$
Using Lemma \ref{lem:i.hat} (in the first and third equalities), the second
relation in \eqref{eq:right_t-module} on $i$ and unitality of $\xi_r$ (in the
second equality) and the first equality in \eqref{eq:ihat-what} defining $\wt
i$ (in the fourth equality), we see that
\begin{eqnarray}\label{eq:i_hat_form}
p {\wt T}_l \circ \sqcap u_l U \circ \sqcap \xi_l \circ i^{-1} U &=&
\theta_l^{-1}\circ \sqcap \xi_r \circ iU\circ \sqcap u_l U \circ \sqcap \xi_l
\circ i^{-1} U  \nonumber\\
&=&\theta_l^{-1}\circ \sqcap \xi_l \circ i^{-1} U =
\theta_l^{-1}\circ \theta_r \circ p {\wt T}_r =
{\wt i}\circ p {\wt T}_r.
\end{eqnarray}
The first condition in \eqref{eq:comod.i} for $({\widetriangle\sqcap },
{\widetriangle i})$ can be expressed as commutativity of the inner square in
the diagram below. Since $p {\wt T}_r$ is a natural epimorphism, the above
considerations imply that it holds true if and only if the outer square in
$$
\xymatrix@C=.75pc{
\sqcap T_r U \ar[rrr]^-{i^{-1}U}\ar[rrd]^-{p {\wt T}_r}\ar[ddd]_-{\sqcap T_r
  u_r U}&&&
\sqcap T_lU \ar[rrr]^-{\sqcap \xi_l}&&&
\sqcap U \ar[rrr]^-{\sqcap u_lU}&&&
\sqcap T_lU\ar[lld]_-{p {\wt T}_l}\ar[dd]^-{\sqcap T_l u_l U}\\
&&{\wt \sqcap}{\wt T}_r \ar[rrrrr]^-{\wt i}\ar[d]_-{{\wt \sqcap}d_r}&&&&&
{\wt \sqcap}{\wt T}_l \ar[d]^-{{\wt \sqcap}d_l}&&\\
&&{\wt \sqcap}{\wt T}_r^2\ar[rr]^-{{\wt i}{\wt T}_r}&&
{\wt \sqcap} {\wt T}_l {\wt T}_r \ar[r]^-{{\wt \sqcap}{\wt\Phi}}&
{\wt \sqcap} {\wt T}_r {\wt T}_l \ar[rr]^-{{\wt i}{\wt T}_l}&&
{\wt \sqcap} {\wt T}_l^2\ar@{=}[rrd]&&
\sqcap T_l^2U\ar[d]^-{p{\wt T}_l^2}\\
\sqcap T_r^2 U\ar[rrrr]_-{
\sqcap u_lT_rU \circ
\sqcap T_r\xi_l\circ
\sqcap \Phi^{-1}U\circ
i^{-1}T_rU}
\ar[urr]^-{p {\wt T}_r^2}&&&&
\,\,\sqcap T_lT_rU\ar[r]_-{\sqcap \Phi^{-1}U}\ar[u]^-{p{\wt T}_l {\wt T}_r}&
\sqcap T_rT_lU\ar[rrr]_-{
\sqcap u_lT_lU\circ
\sqcap m_l U\circ
i^{-1}T_lU}
\ar[u]_-{p{\wt T}_r {\wt T}_l}&&&
\sqcap T_l^2 U\ar[r]_-{p{\wt T}_l^2}&
{\wt \sqcap}{\wt T}_l^2
}
$$
 commutes, what follows by the compatibility of $\Phi$ and $i$ with the
units $u_l$ and $u_r$ of both monads and unitality of the multiplication
$m_l$, cf. the following string computation.
\begin{equation*}
\includegraphics{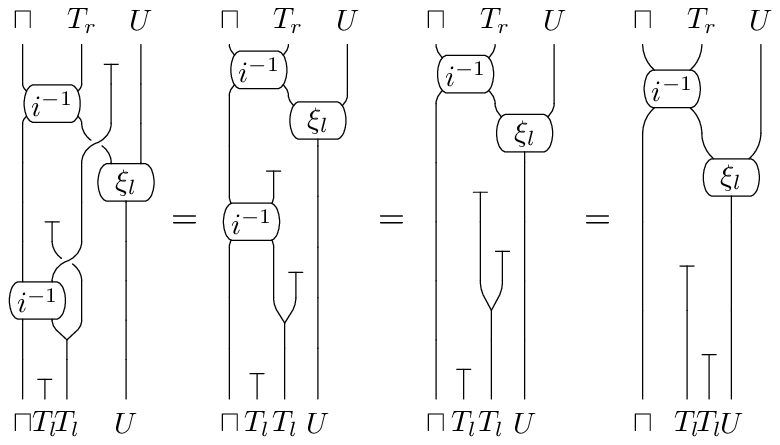}
\end{equation*}
Similarly, by \eqref{eq:i_hat_form}, by naturality of $p$ and Lemma
\ref{lem:Phi.hat}, the second condition in \eqref{eq:comod.i} for
$({\widetriangle\sqcap },{\widetriangle i})$, i.e. the identity ${\wt
  \sqcap}e_l\circ {\wt i}= {\wt \sqcap} e_r$, holds true if and only of the
outer square in
$$
\xymatrix{
\sqcap T_r U \ar[r]^-{i^{-1}U}\ar[rd]^-{p{\wt T}_r}\ar[d]_-{\sqcap \xi_r}&
\sqcap T_lU\ar[r]^-{\sqcap \xi_l}&
\sqcap U\ar[r]^-{\sqcap u_lU}&
\sqcap T_lU\ar[d]^-{\sqcap \xi_l}\ar[ld]_-{p{\wt T}_l}\\
\sqcap U\ar[d]_-{p}&
{\wt \sqcap}{\wt T}_r\ar[r]^-{\wt i}\ar[dl]^-{{\wt \sqcap}e_r}&
{\wt \sqcap}{\wt T}_l\ar[dr]_-{{\wt \sqcap}e_l}&
\sqcap U\ar[d]^-p\\
{\wt \sqcap}\ar@{=}[rrr]&&&
{\wt \sqcap}
}
$$
commutes, what follows by the unitality of $\xi_l$ and the definition of $p$
via the coequalizer in \eqref{eq:Pi.hat}.

Next we prove that $({\wt\sqcup},\widetriangle{w})$
satisfies the conditions in \eqref{eq:comod.i}. Since $U$ is faithful,
applying it to the first relation in \eqref{eq:comod.i} we obtain an
equivalent condition. In view of Lemma \ref{lem:Phi.hat} and the construction
of ${\wt w}$ via the second equality in \eqref{eq:ihat-what}, it takes the
form
$$
T_r^{ 2} w^{-1}\circ
T_r\Phi^{ -1} \sqcup \circ
\Phi^{ -1} T_r\sqcup \circ
T_l T_r w^{-1} \circ
T_l\Phi^{-1}\sqcup \circ
T_l u_l T_r \sqcup
= T_r u_r T_r \sqcup \circ T_r w^{-1} \circ \Phi^{-1} \sqcup.
$$
This holds true by the computation below, where we use the
compatibility between the unit of $T_l$ with $\Phi$ and $w$, and the fact that
$u_l$ is a natural transformation.
\begin{equation*}
\includegraphics{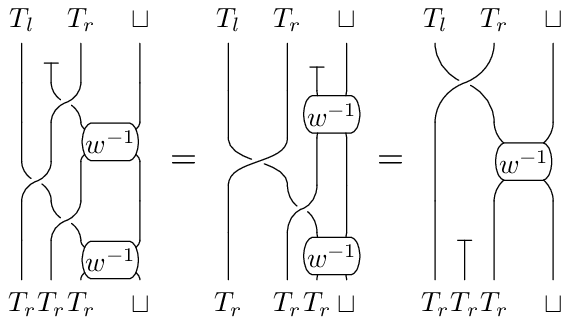}
\end{equation*}
The second relation in \eqref{eq:comod.i} may be proved analogously, showing
that it is equivalent to the fact that $w$ is compatible with the
multiplications of $T_l$ and $T_r$. In conclusion, $({\widetriangle T}_l,{
\widetriangle T}_r,{\widetriangle \Phi}, {\widetriangle \sqcap}, {
\widetriangle i}, {\widetriangle w})$ is an object of ${\mathcal{B}}$.

Next we define the functor $(\widetriangle{-})$ on morphisms. Let $
(G,q_{l},q_{r},\wedge ,\pi ,\vee ,\omega )$ be a morphism in ${\mathcal{A}}
_{c}^{\times }$. Its image in ${\mathcal{B}}^{\times}$ under
$(\widetriangle{-})$ will be denoted by $({\widetriangle
  G},{\widetriangle q}_{l},{\widetriangle q}_{r},\wedge,{\widetriangle\pi
},\vee , {\widetriangle\omega })$. The functor ${\widetriangle G}:
{\mathcal{M}^{\prime }\,}^{T_{l}^{\prime}T_{r}^{\prime }}\rightarrow
\mathcal{M}^{T_{l}T_{r}}$ is defined, for $(M,\rho )$ in ${\mathcal{M}^{\prime
}}^{\,T_{l}^{\prime }T_{r}^{\prime }}$, by $\widetriangle{G}(M,\rho
)=(GM,\widetriangle{\rho})$, where
\begin{equation}
\widetriangle{\rho}=G\rho \circ q_{l}T_{r}^{\prime }M\circ T_{l}q_{r}M.
\label{eq:G-hat}
\end{equation}
On a morphism $f$ in ${\mathcal{M}^{\prime}}^{\,T_{l}^{\prime}
T_{r}^{\prime }}$, we put ${\wt G}f=Gf$.
In fact, $\widetriangle{G}$ is the lifting of $G$ (in the sense that $U{
\widetriangle G}=GU^{\prime }$), which is induced by the monad morphism $
q_{l}T_{r}^{\prime }\circ T_{l}q_{r}$, cf. Lemma \ref{lem:comp_d_law} and
\cite[Lemma 1]{Johnst}. Since
\begin{align*}
\widetriangle{T}_{l}\widetriangle{G}(M,\rho )& =(T_{l}GM,T_{l}G\rho \circ
T_{l}q_{l}T_{r}^{\prime }M\circ T_{l}^2q_{r}M\circ T_{l}u_{l}T_{r}GM\circ
m_{l}T_{r}GM\circ T_{l}\Phi GM), \\
\widetriangle{G}\widetriangle{T}_{l}^{\prime }(M,\rho )& =(GT_{l}^{\prime
}M,GT_{l}^{\prime }\rho \circ GT_{l}^{\prime }u_{l}^{\prime }T_{r}^{\prime
}M\circ Gm_{l}^{\prime }T_{r}^{\prime }M\circ GT_{l}^{\prime }\Phi ^{\prime
}M\circ q_{l}T_{r}^{\prime }T_{l}^{\prime }M\circ T_{l}q_{r}T_{l}^{\prime
}M),
\end{align*}
we can check easily that $q_{l}M:T_{l}GM\rightarrow GT_{l}^{\prime }M$ is a
morphism of $T_{l}T_{r}$-algebras with respect to the
above actions. Hence one may define ${\widetriangle q}_{l}:{\widetriangle T}
_{l}{\widetriangle G\rightarrow \widetriangle G\widetriangle T'_{l}}
$ by $ U{\widetriangle q}_{l}(M,\rho ):=q_{l}M$ and, proceeding
similarly, one may take ${\widetriangle q}_{r}:{\widetriangle
  T}_{r}{\widetriangle G\rightarrow \widetriangle G\widetriangle T_{r}'}$ to
be defined by $ U{\widetriangle q}_{r}(M,\rho):=
q_{r}M$.

The left square in the following diagram is commutative by
naturality, if choosing the upper ones of the parallel arrows. It is
commutative also choosing the lower ones of the parallel arrows, by
\eqref{eq:right_t-mod_morphism} and naturality.
Therefore, universality of the coequalizer in the top row implies the
existence of a unique natural transformation ${\widetriangle\pi }$,
rendering commutative the diagram.
$$
\xymatrix{
\sqcap T_l G U'\ar@<2pt>[rrr]^-{\sqcap G \xi'_l\circ \sqcap q_l U'}
\ar@<-2pt>[rrr]_-{\sqcap G \xi'_r \circ \sqcap q_r U' \circ i G U'}
\ar[d]_-{\pi T'_l U'\circ \sqcap q_l U'} &&&
\sqcap G U'\ar[rr]^-{p {\widetriangle G}} \ar[d]^-{\pi U'}&&
{\widetriangle \sqcap}{\widetriangle G}
\ar@{-->}[d]^-{{\widetriangle \pi}} \\
\L \sqcap'T'_l U'\ar@<2pt>[rrr]^-{\L \sqcap'
  \xi'_l}\ar@<-2pt>[rrr]_-{\L \sqcap' \xi'_r \circ \L i'
  U'}&&&
\L \sqcap' U' \ar[rr]_-{\L p'}&& \L {\widetriangle
  \sqcap} '
}
$$
Finally, we put $U{\widetriangle\omega }:=q_{r}\sqcup' \circ
T_{r}\omega$.
We need to show that this defines indeed a natural transformation
 ${\widetriangle\omega }:{\widetriangle\sqcup } \vee \rightarrow
{\widetriangle G}{\wt \sqcup'} $; that is, that $q_{r}\sqcup' \circ
T_{r}\omega$ commutes with the $T_{l}T_{r}$-actions on $U{\wt
  \sqcup}\V=T_{r}\sqcup\vee$ and $U{\wt G}{\wt \sqcup'}=GT'_{r}\sqcup'$.
For a diagrammatic proof see the following computation, where for the first and
third equations one uses that $(G,q_l)$ and $(G,q_r)$ are morphisms of
monads. The second and the fourth relations follow by
\eqref{eq:left_t-mod_morphism}. In the first three equalities we also use that
all maps are natural transformations.
\begin{equation*}
\includegraphics{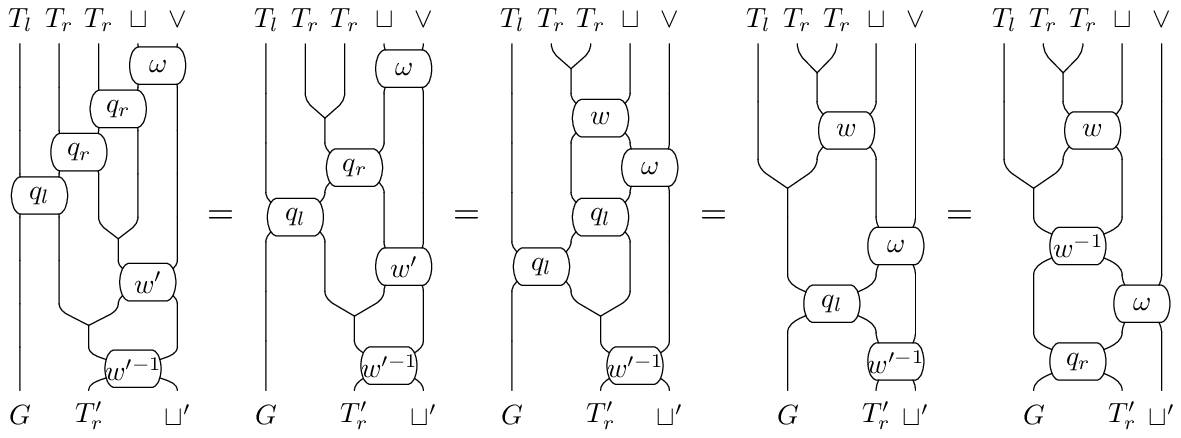}
\end{equation*}

We turn to proving that the so constructed datum $({\widetriangle
G}, {\widetriangle q}_l,{\widetriangle q}_r, \L, {\widetriangle \pi},\V,
{\widetriangle \omega})$ yields a morphism in ${\mathcal B}$.
In order to check that $(G,{\widetriangle q}_r)$ is a comonad
morphism, it suffices to show that
$$
U {\widetriangle q}_r {\widetriangle T}'_r \circ U {\widetriangle
T}_r {\widetriangle q}_r \circ U d_r {\widetriangle G} =
U{\widetriangle G} d'_r \circ U {\widetriangle q}_r
\qquad \textrm{and}\qquad
U e_r {\widetriangle G}=U {\widetriangle G} e'_r \circ U {\widetriangle q}_r ,
$$
since the forgetful functor $U$ is faithful. By Lemma \ref{lem:Phi.hat},
the first condition is equivalent to
$$
q_r T'_r U'\circ T_r q_r U' \circ T_r u_r G U' = G T'_r u'_r U' \circ q_r U',
$$
which holds true since $(G,q_r)$ is a monad morphism and by
naturality. The second condition holds true by construction of the functor
${\widetriangle G}$ (cf. \eqref{eq:G-hat}) and the
relations in Lemma \ref{lem:Phi.hat} on $e_r$ and $e'_r$.
Symmetrically, also $(G,\widetriangle{q}_{l})$ is a comonad
morphism. By faithfulness of $U$, the first condition in
\eqref{eq:com_d-law_comp} is equivalent to
$$
q_r T'_l U'\circ T_r q_l U' \circ \Phi^{-1} G U' = G \Phi^{\prime -1} U'\circ
q_l T'_r U' \circ T_l q_r U',
$$
which holds true by \eqref{eq:d-law_comp}.
The second condition in  \eqref{eq:com_d-law_comp} is equivalent to
commutativity of the inner square in the following diagram.
Since $p{\widetriangle T}_r {\widetriangle G}$ is a natural epimorphism, it
follows by the constructions of the morphisms $\wt q_r$, $\wt q_l$ and $\wt
\pi$ any by the equality \eqref{eq:i_hat_form} that the second condition in
\eqref{eq:com_d-law_comp} holds true if and only if the outer square in
$$
\xymatrix{
\sqcap T_rGU'\ar[rr]^-{\sqcap q_r U'}\ar[rd]^-{p{\wt T}_r{\wt G}}
\ar[d]_-{i^{-1}GU'}&&
\sqcap GT'_r U' \ar[rr]^-{\pi T'_rU'}\ar[d]^-{p{\wt G}{\wt T'}_r}&&
\L\sqcap' T'_rU'\ar[d]^-{\L i^{\prime -1} U'}\ar[dl]_-{\L p'{\wt T'}_r}\\
\sqcap T_lGU'\ar[d]_-{\sqcap q_lU'}&
{\wt \sqcap}{\wt T}_r{\wt G}\ar[r]^-{{\wt \sqcap}{\wt q}_r}
\ar[dd]^-{{\wt i}{\wt G}}&
{\wt \sqcap}{\wt G}{\wt T'}_r\ar[r]^-{{\wt \pi}{\wt T'}_r}&
\L {\wt \sqcap'}{\wt T'}_r\ar[dd]^-{\L{\wt i'}}&
\L\sqcap'T'_lU'\ar[d]^-{\L \sqcap'\xi'_l}\\
\sqcap G T'_lU'\ar[d]_-{\sqcap G \xi'_l}&&&&
\L \sqcap'U'\ar[d]^-{\L \sqcap'u'_lU'}\\
\sqcap G U'\ar[d]_-{\sqcap u_lGU'}&
{\wt \sqcap}{\wt T}_l{\wt G}\ar[r]^-{{\wt \sqcap}{\wt q}_l}&
{\wt \sqcap}{\wt G}{\wt T'}_l \ar[r]^-{{\wt \pi}{\wt T'}_l}&
\L {\wt \sqcap'}{\wt T'}_l&
\L \sqcap' T'_l U' \ar[d]^-{\L p'{\wt T'}_l}
\\
\sqcap T_l GU'\ar[rr]_-{\sqcap q_l U'}
\ar[ur]_(.65){\!\!\!\!\!\!p{\wt T}_l{\wt G}}&&
\sqcap G T'_lU'\ar[r]_-{\pi T'_lU'}\ar[u]_-{p{\wt G}{\wt T'}_l}&
\L \sqcap'T'_l U'\ar[r]_-{\L p'{\wt T'}_l}&
\L {\wt \sqcap'}{\wt T'}_l\ar@{=}[lu]
}
$$
commutes. This follows by the compatibility of $q_l$ with $u_l$, naturality
and \eqref{eq:right_t-mod_morphism}.
Finally, by faithfulness of $U$, \eqref{eq:comod.omega} is
equivalent to
$$
G T'_r w^{\prime -1}\circ G \Phi^{\prime -1} \sqcup' \circ q_l T'_r \sqcup '
\circ T_l q_r \sqcup' \circ T_l T_r \omega =
q_r T'_r \sqcup' \circ T_r q_r \sqcup' \circ T_r^2 \omega \circ T_r w^{-1} \V
\circ \Phi^{-1} \sqcup \V.
$$
This holds true by naturality, \eqref{eq:d-law_comp} and
\eqref{eq:left_t-mod_morphism}.
This finishes the construction of the functor $(\wt{-})$.
 It is straightforward to see that it is a functor indeed, i.e. it
preserves identity morphisms and composition.

It remains to construct a natural isomorphism
$\tau:{\mathcal Z}^{\times}_{*} (\widetriangle{-}) \to
\widehat{{\mathcal Z}^{*\times}(-)}$. For a given object
$({T}_l,{T}_r,{\Phi},{\sqcap},{i}, { \sqcup},{w})$ in $\mathcal{A}^\times_c$,
let $({\widetriangle T}_l,{\widetriangle T}_r,{\widetriangle\Phi},
{\widetriangle \sqcap},{\widetriangle i},{\widetriangle \sqcup},{\widetriangle
w})$ denote the object in $\mathcal{B}^\times$, constructed as in the first
part of the proof.
The para-cyclic object $\wt{\mathcal
Z}_*^{\times}:={\mathcal Z}_*^{\times}({\widetriangle T}_l,{\widetriangle
T}_r,{\widetriangle\Phi},{\widetriangle \sqcap},{\widetriangle i},
{\widetriangle \sqcup}, {\widetriangle w})$, associated to the latter
object as in Theorem \ref{thm:Z_*}, is given, at any non-negative degree
$n$, by the functor ${\widetriangle \sqcap}{\widetriangle T}_l^{n+1}
{\widetriangle\sqcup}$.
The para-cocyclic object ${\mathcal Z}^{*{\times}}:={\mathcal
  Z}^{*\times}(T_l,T_r,\Phi,\sqcap,i,\sqcup,w)$ in Theorem \ref{thm:Z^*}
(and thus also its cyclic dual) is given at degree $n$ by $\sqcap
T_l^{n+1}\sqcup$.
The desired natural transformation $\tau_n$ is defined as the composition of
the following two morphisms
\begin{equation}\label{Tau}
\xymatrix{
{\widetriangle \sqcap}{\widetriangle T}_l^{n+1} {\widetriangle\sqcup}
\ar[r]^-{\theta_l\wt{T}_l^n\wt{\sqcup}}&
\sqcap U {\wt T}_l^{n}{\wt \sqcup}=\sqcap T_l^nT_r\sqcup
\ar[r]^-{\sqcap T_l^nw}&
\sqcap T_l^{n+1}\sqcup}.
\end{equation}
Clearly $\tau_n$ is a natural isomorphism.
We claim that $\tau_\ast$ is also an isomorphism of para-cyclic objects
between $\wt{\mathcal Z}_*^{\times}$  and the cyclic dual $\widehat{{\mathcal
Z}^{*{\times}}}$ of ${\mathcal Z}^{*{\times}}$. Applying  Connes's cyclic
duality functor (in the form it can be found in \cite{KR}) to
${\mathcal Z}^{*{\times}}$, for every degree $n$ the para-cyclic morphism
$\hat{t}_n:\sqcap T_l^{n+1}\sqcup \to \sqcap T_l^{n+1}\sqcup$ of
$\widehat{{\mathcal Z}^{*{\times}}}$ comes out as
$$
\hat{t}_n = i^{-1} T_l^{n}\sqcup \circ \sqcap \Phi^{-1} T_l^{n-1}\sqcup \circ
\dots \circ \sqcap T_l^{n-1} \Phi^{-1} \sqcup \circ \sqcap T_l^n
w^{-1}.
$$
On the other hand, by Theorem \ref{thm:Z_*}, the
para-cyclic morphism
$\wt{t}_n:
{\widetriangle \sqcap}{\widetriangle T}_l^{n+1} {\widetriangle\sqcup}\to
{\widetriangle \sqcap}{\widetriangle T}_l^{n+1} {\widetriangle\sqcup}$
of $\wt{\mathcal Z}_*^{\times}$ is given by
$$
{\widetriangle t}_n= {\widetriangle i}{\widetriangle T}_l^n{\widetriangle
\sqcup} \circ {\widetriangle \sqcap} {\widetriangle \Phi}{\widetriangle
T}_l^{n-1} {\widetriangle \sqcup} \circ {\widetriangle \sqcap}{\widetriangle
T}_l{\widetriangle \Phi}{\widetriangle T}_l^{n-2} {\widetriangle \sqcup}
\circ \dots \circ {\widetriangle \sqcap}
{\widetriangle T}_l^{n-1} {\widetriangle \Phi}{\widetriangle \sqcup} \circ
{\widetriangle \sqcap} {\widetriangle T}_l^{n}{\widetriangle w},
$$
for every $n$. Since $p\wt{T}^{n+1}\wt{\sqcup}$ is an epimorphism,
$\tau_n$ commutes with these para-cyclic morphisms if and only if
\begin{equation}\label{CycOpCommute}
\hat{t}_n\circ\tau_n\circ
p\wt{T}^{n+1}\wt{\sqcup}=\tau_n\circ\wt{t}_n\circ p\wt{T}^{n+1}\wt{\sqcup}.
\end{equation}
For the proof of this relation, using string computations, see the following
sequence of equations.
\[
\includegraphics{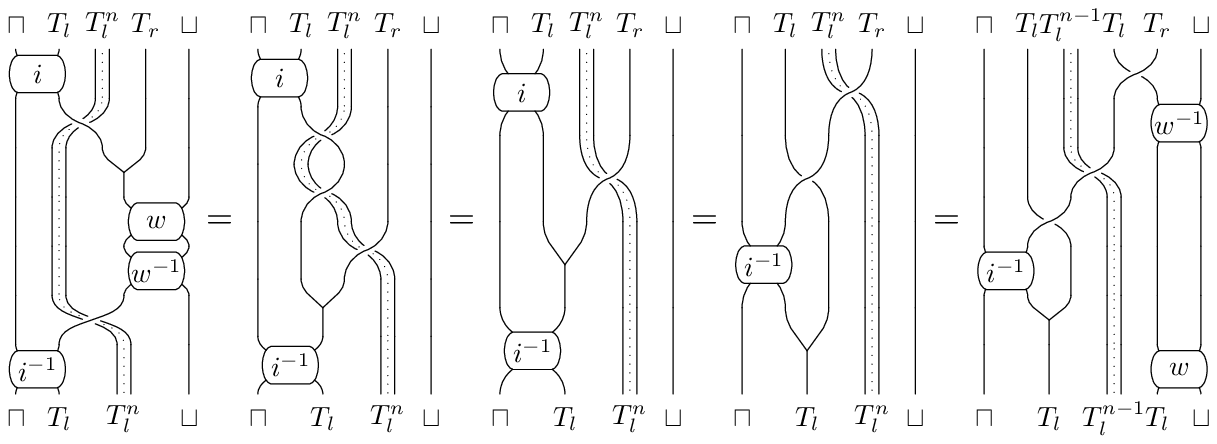}
\]
Taking into account the formulae that give $\hat{t}_n$ and $\tau_n$, on one
hand, and the defining relation  of $\theta_l$ in Lemma \ref{lem:i.hat}
together with the definition of $\xi_r$ in \eqref{eq:xi_l_r},
on the other hand, one can see that the first diagram represents
the left hand side of \eqref{CycOpCommute}. In view of the
definition of the maps $\wt{t}_n$, $\wt{i}_n$, $\wt{w}_n$, $\wt{\Phi}$ and
$\theta_r$, the last diagram corresponds to the right hand side
of \eqref{CycOpCommute}. To deduce the first equality we
use $n$ times the compatibility relation between the multiplication of $T_r$
and $\Phi$. The second relation is a consequence of the fact that the string
and the stripe may be unknotted, by using repeatedly the relation
$\Phi^{-1}\circ\Phi=T_rT_l$. For the third equation one uses that
$(\sqcap,i)$ is right $\Phi$-module functor, i.e. the first condition in
\eqref{eq:right_t-module}. The last relation follows by using that $w$
is natural.

For $0\leq k\leq n$, let ${\widetriangle d}_k:={\widetriangle
\sqcap}{\widetriangle T}_l^k e_l{\widetriangle T}_l^{n-k} {\widetriangle
\sqcup}$ be the face maps of $\wt{\mathcal Z}_\ast^\times$.
For $0< k\leq n$, the face maps of $\widehat{{\mathcal Z}^{*{\times}}}$
are $\hat{d}_k:=\sqcap T_l^{k-1} m_l T_l^{n-k}\sqcup$ while
$\hat{d}_0:=\sqcap T_l^{n-1}m_l \sqcup \circ {\hat t}_n^{-1}$.
We prove next that the operators $\tau_\ast$ are compatible with these face
maps, i.e.  for  $0\leq k\leq n$,
\begin{equation}\label{TauCompFace}
\hat{d}_k\circ\tau_n\circ
p\wt{T}_l^{n+1}\wt{\sqcup}=\tau_{n-1}\circ\wt{d}_k\circ
p\wt{T}_l^{n+1}\wt{\sqcup}.
\end{equation}
If $1\leq k\leq n$, then proceeding as in the case of the para-cyclic operator,
one can see that the string representation of \eqref{TauCompFace} is the
following.
\[
\includegraphics{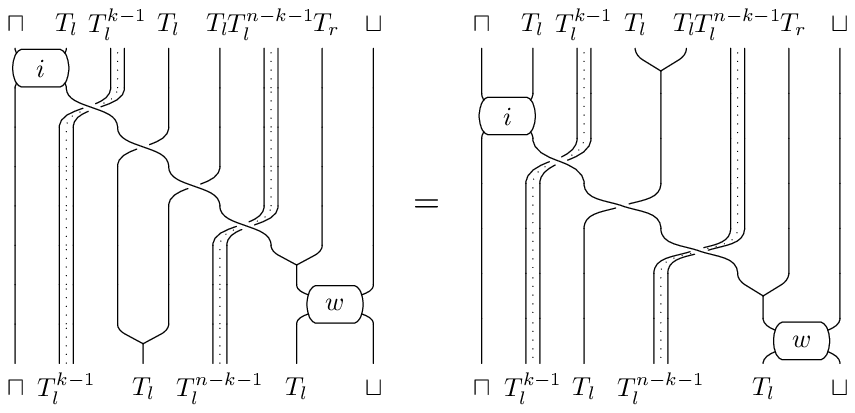}
\]
Visibly this equality holds, as the multiplication of a monad is a natural
transformation and $\Phi$ is a distributive law. Finally, by the definition of
para-cyclic objects, $\hat{d}_0=\hat{d}_n\circ \hat{t}_n^{-1}$ and
$\wt{d}_0=\wt{d}_n\circ \wt{t}_n^{-1}$. Thus, in view of relation
\eqref{CycOpCommute} and of the fact that \eqref{TauCompFace} holds for $k=n$,
we conclude that \eqref{TauCompFace} also holds for $k=0$.
One proves in a similar way that the operators $\tau_\ast$ are compatible with
the codegeneracy maps.

Our final task is to prove naturality of this isomorphism $\tau$. For an
arbitrary morphism
$$
\zeta:=(G,q_l,q_r,\L,\pi,\V,\omega):(T_l,T_r,\Phi,\sqcap,i,\sqcup,w)\to
(T'_l,T'_r,\Phi',\sqcap',i',\sqcup',w')$$
in $\mathcal{A}^\times_c$, let $\wt{\zeta}:=({\widetriangle G}, {\widetriangle
q}_l,{\widetriangle q}_r, \L, {\widetriangle \pi},\V,
{\widetriangle \omega})$ denote the corresponding morphism in ${\mathcal
  B}$. Since $p$ is a natural epimorphism, it is sufficient to
prove that
\begin{equation}\label{eq:nat}
\tau'_\ast\circ\wt{\zeta}_\ast\circ
p\wt{T}_l^{n+1}\wt{\sqcup}{\vee}=
\hat{\zeta}_\ast\circ\tau_\ast\circ p\wt{T}_l^{n+1}\wt{\sqcup}
    {\vee},
\end{equation}
where $\wt{\zeta}_\ast:=\mathcal{Z}^\times_\ast(\wt{\zeta})$
(cf. Theorem \ref{thm:Z_*}),
$\hat{\zeta}_\ast:=\widehat{\mathcal{Z}^{\ast\times}(\zeta)}$
(cf. Theorem \ref{thm:Z^*}); and $\tau_\ast$ and
$\tau'_\ast$ are morphisms constructed as in \eqref{Tau}, corresponding
to the domain and the codomain of $\zeta$, respectively. The string
representations of both sides of \eqref{eq:nat} are given in the first
and last terms of the following sequence of equations.
\[
\includegraphics{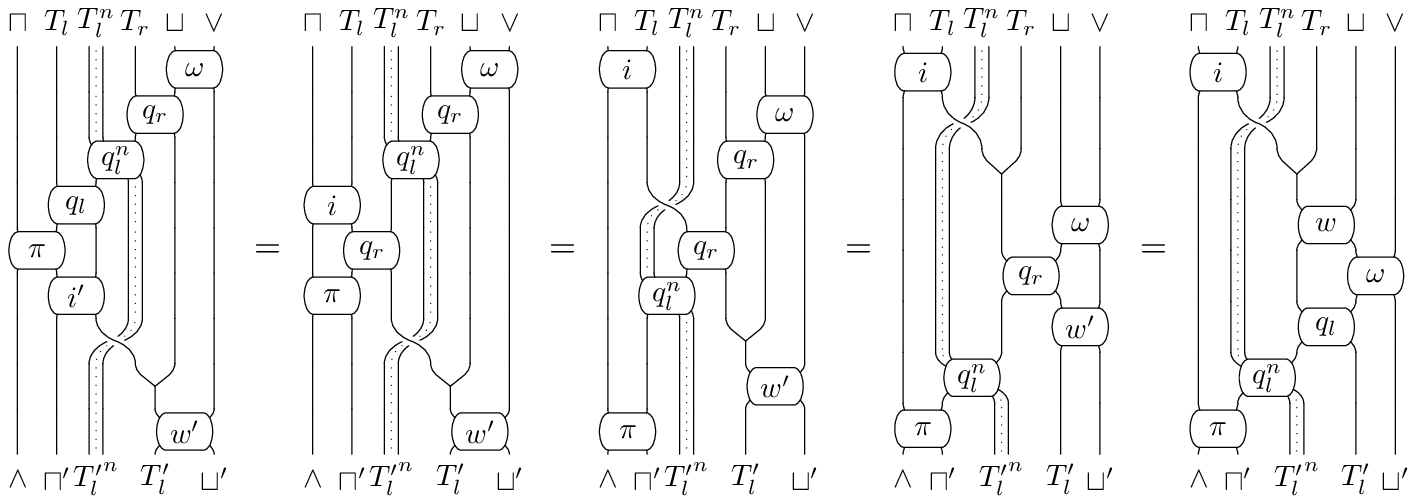}
\]
Recall that $q_l^n:T_l^nG\rightarrow G{T'_l}^n$ is defined inductively by
$q_l^1=q_l$ and $q_l^{n+1}:=q_l{T'_l}^n\circ T_lq_l^n$.
The first identity follows by applying \eqref{eq:right_t-mod_morphism}, while
the second one is obtained by using repeatedly \eqref{eq:d-law_comp} and
naturality. The third equation follows as $q_r$ is a morphism of monads
and the last equality is a consequence of
 \eqref{eq:left_t-mod_morphism}.

 The statements about the functor $(\wt{-}):{\mathcal B}^\times_e\to
    {\mathcal A}^\times$ are proven symmetrically.
\end{proof}


\section{Back to the examples}\label{sec:ex_duals}

In this section we return to the examples in Section
\ref{sec:ex.A} and Section \ref{sec:ex.B}. We show that they all
belong to ${\mathcal A}^\times_c$ and ${\mathcal B}^\times_e$,
respectively, whenever the occurring bialgebroids correspond to a
Hopf algebroid with a bijective antipode and the  occurring
(co(ntra))modules are (co(ntra))modules of the Hopf algebroid
(cf. Appendix). In light of Theorem \ref{thm:cyclic_dual}, in
this case we may apply the functor $\bigcheck$ to any of these
examples. The examples in Section \ref{sec:ex.A} and Section
\ref{sec:ex.B} turn out to be pairwise related via the functor
$\bigcheck$.

\begin{example}
Let $H$ be a Hopf algebroid over base algebras $L$ and $R$, with a bijective
antipode $S$, and $A$ be a left $H$-module algebra. This means that $A$ is a
left module algebra of the constituent left bialgebroid $H_L$, so there is a
corresponding object $(T_l,T_r,\Phi,\sqcap,i,\sqcup,w)$ of ${\mathcal A}$ in
Example \ref{ex:1}. Obviously, using the forgetful functor
$F:H\textrm{-}\mathsf{Comod} \to H_L\textrm{-}\mathsf{Comod}$, we can
construct another object
\begin{equation}\label{eq:1.d}
(T_l,T_r,\Phi,\sqcap,i,\sqcup F,wF)
\end{equation}
of ${\mathcal A}$.
We claim that \eqref{eq:1.d} belongs to ${\mathcal A}^\times$ (hence to
${\mathcal A}^\times_c$, since $\mathsf{Mod}\textrm{-}k$ possesses
coequalizers). Indeed, $\Phi$
and $i$ are obviously isomorphisms. The inverse of $w_{F M}$ is given, for a
left $H$-comodule $M$ with
coaction $m\mapsto m^{[-1]}\ot_R m^{[0]}$ of the constituent right bialgebroid,
as
$$
A \ot_L M\to M\ot_L A,\qquad a\ot_L m \mapsto m^{[0]}\ot_L S^{-1}(m^{[-1]})a.
$$
Thus we can apply to \eqref{eq:1.d} the functor $\bigcheck$. The resulting
object of ${\mathcal B}^\times$ can be obtained from the object in Example
\ref{ex:1.D}, by composing on the right the left comodule functor
$\sqcup: \mathsf{Comod}$-$H_L \to A$-$\mathsf{Mod}$-$A$ in Example
\ref{ex:1.D} with the forgetful functor
$\mathsf{Comod}\textrm{-}H \to \mathsf{Comod}\textrm{-} H_L$ and with the
isomorphism $I_S:H\textrm{-}\mathsf{Comod}\to \mathsf{Comod}\textrm{-}H$,
induced by the bijective antipode $S$, cf. \ref{app:comod}.
\end{example}

\begin{example}
Let $H$ be a Hopf algebroid over base algebras $L$ and $R$, with a bijective
antipode $S$, and $A$ be a right $H$-comodule algebra. Then $A$ is in
particular a right comodule algebra of the constituent right
bialgebroid. Hence there is a corresponding object of ${\mathcal A}$ as in
Example \ref{ex:2}. We claim that it belongs to ${\mathcal A}^\times$ (hence
to ${\mathcal A}^\times_c$). Indeed,
$\Phi$ and $i$ are obviously isomorphisms. The inverse of $w_N$ is given, for
any right $H$-module $N$, by
$$
A\ot_R N \to N \ot_R A,\qquad a \ot_R m \mapsto  m S^{-1}(a_{[1]}) \ot_R a_{[0]},
$$
where $a \mapsto a_{[0]} \ot_L  a_{[1]}$ denotes the coaction of the
constituent left bialgebroid.
Thus we can apply the functor $\bigcheck$. The resulting
object of ${\mathcal B}^\times$ can be obtained from the object in Example
\ref{ex:2.D}, by composing on the right the left comodule functor
$\sqcup: H$-$\mathsf{Mod}\equiv H_L$-$\mathsf{Mod} \to A$-$\mathsf{Mod}$-$A$
in Example \ref{ex:2.D} with the
isomorphism $I_{S^{-1}}^{-1}:\mathsf{Mod}\textrm{-}H \to H
\textrm{-}\mathsf{Mod}$,  induced by the bijective antipode $S$,
cf. \ref{app:mod}.
\end{example}

\begin{example}
Let $H$ be a Hopf algebroid over base algebras $L$ and $R$, with a bijective
antipode $S$, and $C$ be a right $H$-module coring. This means that $C$ is a
right module coring of the constituent right bialgebroid $H_R$, so there is a
corresponding object $(T_l,T_r,\Phi,\sqcap,i,\sqcup,w)$ of ${\mathcal A}$ in
Example \ref{ex:3}. Using the forgetful functor
$F:H\textrm{-}\mathsf{Ctrmod} \to H_R\textrm{-}\mathsf{Ctrmod}$, we can
construct another object $(T_l,T_r,\Phi,\sqcap,i,\sqcup F,wF)$ of ${\mathcal
A}$.
We claim that the modified object belongs to ${\mathcal A}^\times$ (hence
to ${\mathcal A}^\times_c$). Indeed,
$\Phi$ and $i$ are obviously isomorphisms. The inverse of $w_{F Q}$ is given,
for a left $H$-contramodule $Q$ with structure maps
$\alpha_L:\Hom_{L,-}(H,Q)\to Q$ and $\alpha_R:\Hom_{R,-}(H,Q)\to Q$, by
$$
\Hom_{-,R}(C,Q) \to \Hom_{R,-}(C,Q),\qquad f \mapsto \big(c\mapsto \alpha_L(f(c
S^{-1}(-)))\big).
$$
Thus we can apply the functor $\bigcheck$. The resulting object of
${\mathcal B}^\times$ can be obtained from the object in Example
\ref{ex:3.D}, by composing on the right the left comodule functor
$\sqcup:\mathsf{Ctrmod}$-$H_R\to C$-$\mathsf{Ctrmod}$-$C$ in
Example \ref{ex:3.D} with the forgetful functor
$\mathsf{Ctrmod}\textrm{-}H \to \mathsf{Ctrmod}\textrm{-} H_R$ and
with the isomorphism $I_{S}:H \textrm{-}\mathsf{Ctrmod}\to
\mathsf{Ctrmod}\textrm{-} H$, induced by the bijective antipode
$S$, cf. \ref{app:contra}.
\end{example}

\begin{example}
Let $H$ be a Hopf algebroid over base algebras $L$ and $R$, with a bijective
antipode $S$, and $C$ be a left $H$-comodule coring. Then $C$ is in
particular a left comodule coring of the constituent left
bialgebroid. Hence there is a corresponding object of ${\mathcal A}$ as in
Example \ref{ex:4}. We claim that it belongs to ${\mathcal A}^\times$ (hence
to ${\mathcal A}^\times_c$). Indeed,
$\Phi$ and $i$ are obviously isomorphisms. The inverse of $w_N$ is given, for
any right $H$-module $N$, by
$$
\Hom_{-,L}(C,N)\to \Hom_{L,-}(C,N),\qquad
f \mapsto \big(c\mapsto f(c^{[0]})S^{-1}(c^{[-1]})\big),
$$
where $c\mapsto c^{[-1]} \ot_R c^{[0]}$ denotes the coaction of the constituent
right bialgebroid.
Thus we can apply the functor $\bigcheck$. The resulting
object of ${\mathcal B}^\times$ can be obtained from the object in Example
\ref{ex:4.D}, by composing on the right the left comodule functor
$\sqcup: H$-$\mathsf{Mod}\equiv H_R$-$\mathsf{Mod} \to C$-$\mathsf{Comod}$-$C$
in Example \ref{ex:4.D} with the
isomorphism $I_{S^{-1}}^{-1}:\mathsf{Mod}\textrm{-}H \to H
\textrm{-}\mathsf{Mod}$,  induced by the bijective antipode $S$,
cf. \ref{app:mod}.
\end{example}

\begin{example}
Let $H$ be a Hopf algebroid over base algebras $L$ and $R$, with a bijective
antipode $S$, and $C$ be a left $H$-comodule coring. Then $C$ is in
particular a left comodule coring of the constituent left
bialgebroid. Hence there is a corresponding object of ${\mathcal B}$ as in
Example \ref{ex:5}. We claim that it belongs to ${\mathcal B}^\times$ (hence
to ${\mathcal B}^\times_e$). Indeed,
$\Phi$ and $i$ are obviously isomorphisms. The inverse of $w_N$ is given, for
any left $H$-module $N$, by
$$
N\ot_L C \to C\ot_L N,\qquad
m\ot_L c \mapsto c^{[0]}\ot_L S^{-1}(c^{[-1]})m,
$$
where $c\mapsto c^{[-1]} \ot_R c^{[0]}$ denotes the coaction of the constituent
right bialgebroid. Thus we can apply the functor $\bigcheck$. The resulting
object of ${\mathcal A}^\times$ can be obtained from the object in Example
\ref{ex:5.D}, by composing on the right the left comodule functor
$\sqcup: \mathsf{Mod}$-$H \equiv \mathsf{Mod}$-$H_R \to
C$-$\mathsf{Comod}$-$C$ in Example \ref{ex:5.D} with the
isomorphism $I_{S}:H\textrm{-}\mathsf{Mod}\to \mathsf{Mod}\textrm{-}H$,
induced by the bijective antipode $S$, cf. \ref{app:mod}.
\end{example}

\begin{example}
Let $H$ be a Hopf algebroid over base algebras $L$ and $R$, with a bijective
antipode $S$, and $C$ be a right $H$-module coring. This means that $C$ is a
right module coring of the constituent right bialgebroid $H_R$, so there is a
corresponding object $(T_l,T_r,\Phi,\sqcap,i,\sqcup,w)$ of ${\mathcal B}$ in
Example \ref{ex:6}. Using the forgetful functor
$F:\mathsf{Comod}\textrm{-}H \to \mathsf{Comod}\textrm{-}H_R$, we can
construct another object $(T_l,T_r,\Phi,\sqcap,i,\sqcup F,wF)$ of ${\mathcal
B}$.
We claim that the modified object belongs to ${\mathcal B}^\times$ (hence
to ${\mathcal B}^\times_e$). Indeed,
$\Phi$ and $i$ are obviously isomorphisms. The inverse of $w_{F M}$ is given,
for a right $H$-comodule $M$ with coaction $m\mapsto m_{[0]}\ot_L
m_{[1]}$ of the constituent left bialgebroid, by
$$
M\ot_R C \to C \ot_R M,\qquad m\ot_R c \mapsto c S^{-1}(m_{[1]}) \ot_R m_{[0]}.
$$
Thus we can apply the functor $\bigcheck$. The resulting object of
${\mathcal A}^\times$ can be obtained from the object in Example
\ref{ex:6.D}, by composing on the right the left comodule functor
$\sqcup:H_R$-$\mathsf{Comod}\to C$-$\mathsf{Comod}$-$C$ in Example
\ref{ex:6.D} with the forgetful functor $H \textrm{-}
\mathsf{Comod}\to H_R \textrm{-}\mathsf{Comod}$ and with
the isomorphism $I_{S^{-1}}^{-1}:\mathsf{Comod}\textrm{-} H\to H
\textrm{-}\mathsf{Comod}$, induced by the bijective antipode $S$,
cf. \ref{app:comod}.
\end{example}

\begin{example}
Let $H$ be a Hopf algebroid over base algebras $L$ and $R$, with a bijective
antipode $S$, and $A$ be a left $H$-module algebra. This means that $A$ is a
left module algebra of the constituent left bialgebroid $H_L$, so there is a
corresponding object $(T_l,T_r,\Phi,\sqcap,i,\sqcup,w)$ of ${\mathcal B}$ in
Example \ref{ex:7}. Using the forgetful functor
$F:\mathsf{Ctrmod}\textrm{-}H \to \mathsf{Ctrmod}\textrm{-}H_L$, we can
construct another object $(T_l,T_r,\Phi,\sqcap,i,\sqcup F,wF)$ of ${\mathcal
B}$.
We claim that the modified object belongs to ${\mathcal B}^\times$ (hence
to ${\mathcal B}^\times_e$). Indeed,
$\Phi$ and $i$ are obviously isomorphisms. The inverse of $w_{F Q}$ is given,
for a right $H$-contramodule $(Q,\alpha_L,\alpha_R)$ by
$$
\Hom_{L,-}(A,Q)\to \Hom_{-,L}(A,Q), \qquad
g \mapsto \big(a\mapsto \alpha_R(g(S^{-1}(-)a))\big).
$$
Thus we can apply the functor $\bigcheck$. The resulting object of
${\mathcal A}^\times$ can be obtained from the object in Example
\ref{ex:7.D}, by composing on the right the left comodule functor
$\sqcup:H_L$-$\mathsf{Ctrmod} \to A$-$\mathsf{Mod}$-$A$ in Example
\ref{ex:7.D} with the forgetful functor $H \textrm{-}
\mathsf{Ctrmod}\to H_L \textrm{-}\mathsf{Ctrmod}$ and with the isomorphism
$I_{S^{-1}}^{-1}:\mathsf{Ctrmod}\textrm{-} H\to H \textrm{-}\mathsf{Ctrmod}$,
induced by the bijective antipode $S$, cf. \ref{app:contra}.
\end{example}

\begin{example}
Let $H$ be a Hopf algebroid over base algebras $L$ and $R$, with a bijective
antipode $S$, and $A$ be a right $H$-comodule algebra. Then $A$ is in
particular a right comodule algebra of the constituent right
bialgebroid. Hence there is a corresponding object of ${\mathcal B}$ as in
Example \ref{ex:8}. We claim that it belongs to ${\mathcal B}^\times$ (hence
to ${\mathcal B}^\times_e$). Indeed,
$\Phi$ and $i$ are obviously isomorphisms. The inverse of $w_N$ is given, for
any left $H$-module $N$, by
$$
\Hom_{R,-}(A,N)\to \Hom_{-,R}(A,N),\qquad
g\mapsto \big(a\mapsto S^{-1}(a_{[1]})g(a_{[0]})\big),
$$
where $a \mapsto a_{[0]} \ot_L  a_{[1]}$ denotes the coaction of the
constituent left bialgebroid. Thus we can apply the functor $\bigcheck$. The
resulting object of ${\mathcal A}^\times$ can be obtained from the object in
Example \ref{ex:8.D}, by composing on the right the left comodule functor
$\sqcup:\mathsf{Mod}$-$H\equiv \mathsf{Mod}$-$H_L \to A$-$\mathsf{Mod}$-$A$ in
Example \ref{ex:8.D} with the isomorphism $I_{S}:H \textrm{-}\mathsf{Mod}\to
\mathsf{Mod}\textrm{-} H$, induced by the bijective antipode $S$,
cf. \ref{app:mod}.
\end{example}
\appendix

\section{Modules, comodules and contramodules of Hopf algebroids}

In this appendix we shortly review algebraic structures over non-commutative
base algebras,  which  are used to construct the examples in the
paper. For more information on them we refer to \cite{Bohm:HoA}. Structures as
$R$-rings, $R$-corings, bialgebroids and Hopf algebroids below, generalize the
notions of an algebra, a coalgebra, a bialgebra and a Hopf algebra over a
commutative ring, respectively.

Throughout, let $k$ be a commutative, associative and unital ring.
By an {\em algebra} $R$ we mean an associative and unital algebra over
$k$. The {\em enveloping algebra} $R\ot_k R^{op}$ is denoted by $R^e$. We
tacitly identify $R^e$-modules with $R$-bimodules.

\begin{claim}
An {\em $R$-ring} is a monoid in the monoidal category of $R$-bimodules. In
fact, an $R$-ring $A$ is equivalent to a $k$-algebra $A$, together
with an algebra map $\iota:R\to A$. Denoting the multiplication in an $R$-ring
$A$ by $\mu:A\ot_R A\to A$, there is an induced monad
\begin{equation}\label{eq:right_tensor_monad}
((-)\ot_R A, (-)\ot_R \mu, (-)\ot_R \iota)
\end{equation}
on  the category $\mathsf{Mod}$-$R$  of right $R$-modules.
Algebras of this monad are equivalent to right modules of the
$k$-algebra $A$. Symmetrically, algebras for the monad
\begin{equation}\label{eq:left_tensor_monad}
(A\ot_R (-),\mu\ot_R (-), \iota\ot_R (-))
\end{equation}
on  the category $R$-$\mathsf{Mod}$  of left $R$-modules
are equivalent to left modules of the $k$-algebra $A$. Note that
the same formulae \eqref{eq:right_tensor_monad} and
\eqref{eq:left_tensor_monad} define monads also on  the category
$R$-$\mathsf{Mod}$-$R$  of $R$-bimodules, with respect to
the $R$-actions
\begin{equation}\label{eq:R-actions_tensor}
r(p\ot_R a)r'=rp\ot_R ar'\qquad \textrm{and}\qquad r(a\ot_R
p)r'=r a\ot_R pr' ,
\end{equation}
for $r,r'\in R$, $a\in A$, $p\in P$ and any $R$-bimodule $P$.
\end{claim}

\begin{claim}
For our considerations, $R^e$-rings are of special interest. Note that an
$R^e$-ring is equivalent to an algebra $B$,
together with algebra maps $s:R\to B$ and $t:R^{op}\to B$, such that
$s(r)t(r')=t(r')s(r)$, for all $r,r'\in R$. The maps $s$ and $t$ are known as
the {\em source} and {\em target} maps, respectively. An immediate example of
an $R^e$-ring is the algebra $\End_k(R)$ of $k$-linear endomorphisms of
$R$. It is an algebra via composition of endomorphisms and source and target
maps are
$$
R\to \End_k(R),\quad  r\mapsto r(-)\qquad
R^{op} \to \End_k(R),\quad r \mapsto (-)r.
$$
Any $R^e$-ring $B$ carries four commuting $R$-actions:
$$
r \blacktriangleright b = s(r) b,\qquad
b \blacktriangleleft r = t(r) b,\qquad
r \triangleright b = b t(r),\qquad
b \triangleleft r = b s(r).
$$
In terms of these actions, the following construction can be performed. Take
first the $R$-module tensor product
$$
B \ot_R B := B \ot_k B / \{\ b \blacktriangleleft r - r\blacktriangleright
b, \quad \forall \ b\in B, \ r\in R\ \}
$$
and then the $k$-submodule
$$
B \times_R B := \{\ \Sigma \,  b_i \ot_R b'_i \in B \ot_R B\ |\ \Sigma\,
r\triangleright b_i \ot_R b'_i = \Sigma\, b_i \ot_R b'_i \triangleleft r,
\ \forall \ r\in R\ \}.
$$
$B\times_R B$ is known as the {\em Takeuchi product}, and it is easily checked
to be an $R^e$-ring with factorwise multiplication and source and target maps
$$
R \to B \times_R B,\quad  r\mapsto s(r)\ot_R 1_B, \qquad
R^{op}\to B \times_R B,\quad  r\mapsto 1_B \ot_R t(r).
$$
\end{claim}

\begin{claim}
An {\em $R$-coring} is a comonoid in the monoidal category of
$R$-bimodules. That is, an $R$-coring is an $R$-bimodule $C$,
equipped with an $R$-bilinear coassociative comultiplication
$\Delta:C \to C \ot_R C$, possessing an $R$-bilinear counit
$\epsilon:C \to R$.  For the comultiplication we use the index
notation $c \mapsto c\1 \ot_R c\2$, where implicit summation is
understood.

Any $R$-coring $C$ induces a comonad
\begin{equation}\label{eq:right_tensor_comonad}
((-)\ot_R C, (-)\ot_R \Delta,(-)\ot_R \epsilon)
\end{equation}
on $\mathsf{Mod}$-$R$. Coalgebras of this comonad are called {\em
right $C$-comodules}. Explicitly, this means right $R$-modules
$M$, equipped with a right $R$-linear coaction $M\to M \ot_R C$,
subject to coassociativity and counitality constraints.  For a
right coaction, the index notation $m \mapsto m_{[0]}\otimes_R
m_{[1]}$ is used (with lower, or with upper indeces), where
implicit summation is understood. Symmetrically, coalgebras
for the comonad
\begin{equation}\label{eq:left_tensor_comonad}
(C\ot_R (-),\Delta \ot_R (-), \epsilon \ot_R (-))
\end{equation}
on $R$-$\mathsf{Mod}$ are called {\em left $C$-comodules}.  For
the coaction on a left $C$-comodule, we use the index notation $m
\mapsto m_{[-1]}\ot_R m_{[0]}$ (with lower, or with upper
indeces), where implicit summation is understood. Morphisms of
(right or left) $C$-comodules are morphisms of coalgebras for the
appropriate comonad (\eqref{eq:right_tensor_comonad} or
\eqref{eq:left_tensor_comonad}). That is, (right or left)
$R$-module maps which are compatible with the $C$-coaction.

Note that the same formulae \eqref{eq:right_tensor_comonad} and
\eqref{eq:left_tensor_comonad} define comonads also on $R$-$\mathsf{Mod}$-$R$,
with respect to the $R$-bimodule structures as in \eqref{eq:R-actions_tensor}.

To an $R$-coring $C$, one can associate also monads. The triple
\begin{equation}\label{eq:right_hom_monad}
(\Hom_{-,R}(C,-),\Hom_{-,R}(\Delta,-),\Hom_{-,R}(\epsilon,-))
\end{equation}
is a monad on $\mathsf{Mod}$-$R$ (where we used standard hom-tensor identities
to identify \break
$\Hom_{-,R}(C\ot_R
C,-)\cong \Hom_{-,R}(C,\Hom_{-,R}(C,-))$ and $\Hom_{-,R}(R,-)\cong
\mathsf{Mod}$-$R$). Algebras of this monad are called {\em right
  $C$-contramodules}, cf.  \cite{EM}, \cite{BBW}. Symmetrically, left
$C$-contramodules are algebras of the monad
\begin{equation}\label{eq:left_hom_monad}
(\Hom_{R,-}(C,-),\Hom_{R,-}(\Delta,-),\Hom_{R,-}(\epsilon,-))
\end{equation}
on $R$-$\mathsf{Mod}$.
 Morphisms of (right or left) $C$-contramodules are morphisms of algebras
for the appropriate monad (\eqref{eq:right_hom_monad} or
\eqref{eq:left_hom_monad}). That is, (right or left) $R$-module maps
which are compatible with the contramodule structure.

Note that the same formulae \eqref{eq:right_hom_monad}
and \eqref{eq:left_hom_monad} define monads also on $R$-$\mathsf{Mod}$-$R$,
with respect to the $R$-bimodule structures
$$
(rfr')(c)=rf(r'c)\qquad \textrm{and} \qquad (rgr')(c)=g(cr)r',
$$
for $r,r'\in R$, $c\in C$, $f\in \Hom_{-,R}(C,P)$, $g\in \Hom_{R,-}(C,P)$ and
any $R$-bimodule $P$.
\end{claim}

\begin{claim}
A {\em left $R$-bialgebroid} \cite{Tak:bgd}, \cite{Lu:hgd}
is an $R^e$-ring $(B,s,t)$, possessing also an
$R$-coring structure $(B,\blacktriangleright, \blacktriangleleft, \Delta,
\epsilon)$, subject to the following compatibility axioms.
\begin{itemize}
\item The comultiplication $\Delta:B \to B \ot_R B$ factorizes through $B
  \times_R B$;
\item corestriction of $\Delta$ is a homomorphism of $R^e$-rings $B \to B
  \times_R B$;
\item the map $B \to \End_k(R)$, $ b\mapsto \epsilon(b s(-))$ is a
  homomorphism of $R^e$-rings.
\end{itemize}
Some equivalent forms of the definition can be found e.g. in \cite{BrzMil:bgd}.
The notion of a {\em right $R$-bialgebroid} is obtained symmetrically, by
interchanging the roles of the $R$-actions
$(\blacktriangleright,\blacktriangleleft)$ and
$(\triangleright,\triangleleft)$ in an $R^e$-ring, given by multiplication on
the right, and on the left, respectively. For more details we refer to
\cite{KadSzl:D2bgd} or \cite{Bohm:HoA}.
\end{claim}

\begin{claim}
{\em Modules} of an $R$-bialgebroid $B$ are modules of the
underlying $k$-algebra $B$. Since $B$ is an $R^e$-ring, there is a
forgetful functor from the category of (left or right) $B$-modules
to the category of (left or right) $R^e$-modules, equivalently, to
the category of $R$-bimodules. By \cite[Theorem 5.1]{Scha:bia_nc},
the category of left (resp. right) modules of a left (resp. right)
bialgebroid is a monoidal category, with monoidal product given by
the $R$-module tensor product. Left (resp. right) {\em module
algebras} of a left (resp. right) bialgebroid $B$ are defined as
monoids in the monoidal category of left (resp. right)
$B$-modules. $B$-module algebras are thus in particular $R$-rings.
By the same principle, left (resp. right) {\em $B$-module corings}
are comonoids in the monoidal category of left (resp. right)
$B$-modules.  They are in particular $R$-corings.
\end{claim}

\begin{claim}
{\em Comodules} of a (left or right) $R$-bialgebroid are comodules
of the constituent $R$-coring. As a consequence of the bialgebroid
axioms, any right comodule $M$ of a right $R$-bialgebroid $B$ can
be equipped also with a unique left $R$-action, such that the
range of the coaction $m \mapsto m^{[0]}\ot_R m^{[1]}$ lies within
the center of the $R$-bimodule $M \ot_R B$. That is, for any $m\in
M$ and $r\in R$, $rm^{[0]} \ot_R m^{[1]}= m^{[0]} \ot_R t(r)
m^{[1]}$ (where $t:R^{op} \to B$ is the target map). This equips
any $B$-comodule with an $R$-bimodule structure and the category
of right $B$-comodules becomes monoidal with respect to the
$R$-module tensor product (cf. \cite[Theorem 3.18]{Bohm:HoA}). In
other words, there is a strict monoidal `forgetful' functor from
the category $\mathsf{Comod}$-$B$  of right $B$-comodules
to $R$-$\mathsf{Mod}$-$R$. Symmetrically, also the category of
left comodules of a right $R$-bialgebroid is monoidal, via
$\ot_{R^{op}}$. In the same way, categories of left and right
comodules of a left $R$-bialgebroid are monoidal, with respect to
the $R$-module tensor product and the $R^{op}$-module tensor
product, respectively. Left (resp. right) {\em comodule algebras}
of a left or right bialgebroid $B$ are defined as monoids in the
monoidal category of left (resp. right) $B$-comodules.
$B$-comodule algebras are in particular $R$-rings or
$R^{op}$-rings (depending on the monoidal product of the
appropriate comodule category). By the same principle, left (resp.
right) {\em $B$-comodule
  corings} are comonoids in the monoidal category of left (resp. right)
$B$-comodules (hence they are $R$, or $R^{op}$-corings).
\end{claim}

\begin{claim}
{\em Contramodules} of a (left or right) $R$-bialgebroid are contramodules of
the constituent $R$-coring.
As a consequence of the bialgebroid axioms, any left contramodule
$(Q,\alpha:\Hom_{R,-}(B,Q)\to Q)$ of a right $R$-bialgebroid $B$ can be
equipped also with a right $R$-action
$$
qr:= \alpha \big( \epsilon(s(r)-) q \big),
$$
for $q\in Q$, $r\in R$, such that  $Q$ becomes an
$R$-bimodule. This construction yields a `forgetful' functor from
the category $B$-$\mathsf{Ctrmod}$  of left $B$-contramodules
to $R$-$\mathsf{Mod}$-$R$. Symmetrically, also right
contramodules of a right $R$-bialgebroid and left and right
contramodules of a left $R$-bialgebroid possess canonical
$R$-bimodule structures. Note, however, that the category of
contramodules of an arbitrary bialgebroid is not known to be
monoidal.
\end{claim}

\begin{claim} A {\em Hopf algebroid} $H$ consists of a left
bialgebroid structure $(H,s_L,t_L,\Delta_L,\epsilon_L)$ over a base algebra
$L$, and a right bialgebroid structure $(H,s_R,t_R,\Delta_R,\epsilon_R)$ over
a base algebra $R$, on the same $k$-algebra $H$, together with
a $k$-module map $S:H\to H$, called the {\em antipode}. These structures are
subject to the following axioms.
\begin{itemize}
\item $s_R \circ \epsilon_R \circ t_L =t_L$,
$t_R \circ \epsilon_R \circ s_L =s_L$,
$s_L \circ \epsilon_L \circ t_R =t_R$,
$t_L \circ \epsilon_L \circ s_R =s_R$.
\item $(\Delta_R \ot_L H)\circ \Delta_L = (H \ot_R \Delta_L )\circ \Delta_R$
  and
$(\Delta_L \ot_R H)\circ \Delta_R = (H \ot_L \Delta_R)\circ \Delta_L$.
\item For all $h\in H$, $l\in L$ and $r\in R$, $S(t_L(l)ht_R(r))=s_R(r) S(h)
  s_L(l)$.
\item $\mu_R \circ (H \ot_R S)\circ \Delta_R=s_L \circ \epsilon_L$ and
$\mu_L \circ (S \ot_L H)\circ \Delta_L=s_R \circ \epsilon_R$,
\end{itemize}
where $\mu_R:H\ot_R H \to H$ denotes multiplication in the $R$-ring $s_R:R\to
H$ and $\mu_L:H\ot_L H \to H$ denotes multiplication in the $L$-ring $s_L:L\to
H$.
Note that the second axiom is meaningful because of the first axiom and the
fourth axiom is meaningful because of the third one.

These axioms imply that the algebras $L$ and $R$ are anti-isomorphic, and the
antipode is a bialgebroid morphism from the constituent left bialgebroid to the
opposite-coopposite of the right bialgebroid, and also from the constituent
right bialgebroid to the opposite-coopposite of the left bialgebroid.
\end{claim}

\begin{claim} \label{app:mod}
{\em Modules} of a Hopf algebroid $H$ are by definition
modules of the underlying $k$-algebra. In this way the category of (left or
right) $H$-modules coincides with the (left or right) module category of any
of the  constituent bialgebroids. Hence both categories of left and right
$H$-modules are monoidal. A (left or right) {\em module algebra} of a Hopf
algebroid $H$ is defined as a monoid in the monoidal category of (left or
right) $H$-modules. Similarly, a (left or right) {\em module
coring} of a Hopf algebroid $H$ is defined as a comonoid in the monoidal
category of (left or right) $H$-modules.

If the antipode $S$ of a Hopf algebroid $H$ is bijective, then it
induces an isomorphism $I_S:H\textrm{-}\mathsf{Mod}\to
\mathsf{Mod}\textrm{-}H$ between the categories of left and right
$H$-modules. This isomorphism takes a left $H$-module $N$ to
$N$ as a right $H$-module with action $n \triangleleft h =
S^{-1}(h)n$. On the morphisms $I_S$ acts as the identity map. A
similar isomorphism $I_{S^{-1}}:H\textrm{-}\mathsf{Mod}\to
\mathsf{Mod}\textrm{-}H$ is obtained by replacing $S$ by $S^{-1}$.
\end{claim}

\begin{claim} \label{app:comod}
Right {\em comodules} of a Hopf algebroid $H$ over base
algebras $L$ and $R$ are triples $(M,\varrho_L,\varrho_R)$, where $M$ is a
right $L$-module and a right $R$-module, $(M,\varrho_L)$ is a right comodule
of the constituent left bialgebroid, $(M,\varrho_R)$ is a right comodule
of the constituent right bialgebroid, such that both coactions are comodule
maps for the other bialgebroid as well. That is, $\varrho_L$ is a right
$R$-module map, $\varrho_R$ is a right $L$-module map and
the compatibility conditions
$$
(M\ot_L \Delta_R) \circ \varrho_L = (\varrho_L \ot_R H) \circ \varrho_R
\qquad \textrm{and} \qquad
(M\ot_R \Delta_L) \circ \varrho_R = (\varrho_R \ot_L H) \circ \varrho_L
$$
hold. It follows that the right $R$-, and $L$-actions on $M$
commute, i.e. $M$ is a right $R\ot L$-module. Morphisms of
$H$-comodules are defined as comodule maps for both constituent
bialgebroids.  Right comodules of a Hopf algebroid $H$ and their
morphisms constitute the category $\mathsf{Comod}\textrm{-}H$. The
category $H\textrm{-}\mathsf{Comod}$ of left $H$-comodules is
defined symmetrically.

The category of (left or right) comodules of a Hopf algebroid is monoidal and
the forgetful functors to the comodule categories of the constituent
bialgebroids are strict monoidal \cite[Theorem 4.9]{Bohm:HoA}. A (left or
right) {\em comodule algebra} of a Hopf algebroid $H$ is defined as a
monoid in the monoidal category of (left or right)
$H$-comodules. Similarly, a (left or right) {\em comodule coring} of a
Hopf algebroid $H$ is defined as a comonoid in the monoidal category of
(left or right) $H$-comodules.

If the antipode $S$ of a Hopf algebroid $H$ is bijective, then it
induces an isomorphism $I_S:H\textrm{-}\mathsf{Comod}\to
\mathsf{Comod}\textrm{-}H$. Take a left $H$-comodule $M$, with
coaction $m\mapsto m_{[-1]} \ot_L m_{[0]}$ of the constituent left
$L$-bialgebroid and coaction $m\mapsto m^{[-1]} \ot_R m^{[0]}$ of
the constituent right $R$-bialgebroid. The isomorphism $I_S$ takes
it to $M$ as a right $H$-comodule,  with right $R$, and
$L$-actions induced by the algebra isomorphism $R \cong L^{op}$
and with coaction $m\mapsto m^{[0]} \ot_L S^{-1}(m^{[-1]})$ of
the constituent left $L$-bialgebroid and coaction $m\mapsto
m_{[0]} \ot_R S^{-1}(m_{[-1]})$ of the constituent right
$R$-bialgebroid. On the morphisms $I_S$ acts as the identity map.
A similar isomorphism $I_{S^{-1}}:H\textrm{-}\mathsf{Comod}\to
\mathsf{Comod}\textrm{-}H$ is obtained by replacing $S$ by
$S^{-1}$.
\end{claim}

\begin{claim} \label{app:contra}
Right {\em contramodules} of a  Hopf algebroid $H$ over
base algebras $L$ and $R$ are triples $(Q,\alpha_L,\alpha_R)$, where $Q$ is a
right $L$-module and a right $R$-module, $(M,\alpha_L)$ is a right
contramodule of the constituent left bialgebroid, $(M,\alpha_R)$ is a right
contramodule of the constituent right bialgebroid, such that both structure
maps $\alpha_L$ and $\alpha_R$ are contramodule maps for the other bialgebroid
as well. That is, $\alpha_L$ is a right $R$-module map, $\alpha_R$ is a
right $L$-module map and the compatibility conditions
\begin{eqnarray*}
&&\alpha_L \circ \mathsf{Hom}_{-,L}(H,\alpha_R)= \alpha_R \circ
\mathsf{Hom}_{-,R}(\Delta_L,H)
\qquad \textrm{and} \\
&&\alpha_R \circ \mathsf{Hom}_{-,R}(H,\alpha_L)= \alpha_L \circ
\mathsf{Hom}_{-,L}(\Delta_R,H)
\end{eqnarray*}
hold. It follows that the right $R$-, and $L$-actions on $Q$ commute, i.e. $Q$
is a right $R\ot L$-module.
Morphisms of $H$-contramodules are defined as contramodule maps
for both constituent bialgebroids. Right contramodules of a Hopf algebroid $H$
and their morphisms constitute the category $\mathsf{Ctrmod}\textrm{-}H$. The
category $H\textrm{-}\mathsf{Ctrmod}$ of left $H$-contramodules is defined
symmetrically.

If the antipode $S$ of a Hopf algebroid $H$ is bijective, then it
induces an isomorphism $I_S:H\textrm{-}\mathsf{Ctrmod}\to
\mathsf{Ctrmod}\textrm{-}H$. Take a left $H$-contramodule $Q$,
with structure map $\alpha_L:\Hom_{L,-}(H,Q)\to Q$ as a
contramodule of the constituent left $L$-bialgebroid and structure
map $\alpha_R:\Hom_{R,-}(H,Q)\to Q$ as a contramodule of the
constituent right $R$-bialgebroid. The isomorphism $I_S$ takes it
to $Q$ as a right $H$-contramodule,  with right $R$, and
$L$-actions induced by the algebra isomorphism $R \cong L^{op}$
and with structure maps
$$
\Hom_{-,L}(H,Q)\to Q, \ f \mapsto \alpha_R(f\circ S^{-1})
\quad \textrm{and} \quad
\Hom_{-,R}(H,Q)\to Q, \ g \mapsto \alpha_L(g\circ S^{-1}).
$$
On the morphisms it acts as the identity map.
A similar isomorphism $I_{S^{-1}}:H\textrm{-}\mathsf{Ctrmod}\to
\mathsf{Ctrmod}\textrm{-}H$ is obtained by
replacing $S$ by $S^{-1}$.
\end{claim}

\section*{Acknowledgments}
The first named author was financially supported by the Hungarian
Scientific Research Fund OTKA K 68195. The second named author
was financially supported by CNCSIS, Contract 560/2009 (CNCSIS
code ID\_69).

\end{document}